\newif\ifjournalversion \journalversionfalse
\newif\ifmarked \markedfalse
\documentclass[11pt,a4paper]{article}

\newcommand{\MDPIabstract}{}
\newcommand{\MDPIkeyword}{}
\newcommand{\MDPImsc}{}
\newcommand{\MDPIaddress}{}
\newcommand{\Title}[1]{\title{#1}}
\newcommand{\Author}[1]{\author{#1}}
\newcommand{\AuthorNames}[1]{}
\providecommand{\abstract}[1]{}
\renewcommand{\abstract}[1]{\gdef\MDPIabstract{#1}}
\providecommand{\address}[1]{}
\renewcommand{\address}[1]{\gdef\MDPIaddress{#1}}
\newcommand{\keyword}[1]{\gdef\MDPIkeyword{#1}}
\newcommand{\MSC}[1]{\gdef\MDPImsc{#1}}
\newcommand{\MDPIfrontmatter}{%
  \maketitle
  \begin{center}\begin{minipage}{0.92\linewidth}\small
  \noindent{\bf Abstract.} \MDPIabstract\par\medskip
  \noindent{\bf Keywords:} \MDPIkeyword\par\medskip
  \noindent{\bf MSC:} \MDPImsc
  \end{minipage}\end{center}\medskip}
\newcommand{\MDPIbackmatter}[2]{\par\medskip\noindent{\bf #1} #2\par}
\newcommand{\funding}[1]{\MDPIbackmatter{Funding.}{#1}}
\newcommand{\dataavailability}[1]{\MDPIbackmatter{Data availability.}{#1}}
\newcommand{\conflictsofinterest}[1]{%
  \MDPIbackmatter{Conflicts of interest.}{#1}}
\newcommand{\useofartificialintelligence}[1]{%
  \MDPIbackmatter{Use of artificial intelligence.}{#1}}
\newcommand{\reftitle}[1]{}
\newcommand{\PublishersNote}{}
\usepackage{amsthm}
\usepackage{changepage}
\newlength{\extralength}
\newtheoremstyle{mdpi}{12pt}{12pt}{\itshape}{}{\bfseries}{.}{.5em}{}
\theoremstyle{mdpi}
\newcounter{theorem}     \newtheorem{Theorem}[theorem]{Theorem}
\newcounter{lemma}       \newtheorem{Lemma}[lemma]{Lemma}
\newcounter{corollary}   \newtheorem{Corollary}[corollary]{Corollary}
\newcounter{proposition} \newtheorem{Proposition}[proposition]{Proposition}
\newcounter{example}     \newtheorem{Example}[example]{Example}
\newcounter{remark}      \newtheorem{Remark}[remark]{Remark}
\newcounter{definition}  \newtheorem{Definition}[definition]{Definition}
\newcounter{hypothesis}  \newtheorem{Hypothesis}[hypothesis]{Hypothesis}
\newcounter{assumption}  \newtheorem{Assumption}[assumption]{Assumption}
     
\theoremstyle{plain}

\usepackage[sort&compress]{natbib}
\providecommand{\articlenumber}[1]{\expandafter\gdef\csname @articlenumber\endcsname{#1}}
\providecommand{\copyrightyear}[1]{\expandafter\gdef\csname @copyrightyear\endcsname{#1}}
\providecommand{\dateaccepted}[1]{\expandafter\gdef\csname @dateaccepted\endcsname{#1}}
\providecommand{\datepublished}[1]{\expandafter\gdef\csname @datepublished\endcsname{#1}}
\providecommand{\datereceived}[1]{\expandafter\gdef\csname @datereceived\endcsname{#1}}
\providecommand{\daterevised}[1]{\expandafter\gdef\csname @daterevised\endcsname{#1}}
\providecommand{\firstpage}[1]{\expandafter\gdef\csname @firstpage\endcsname{#1}}
\providecommand{\hreflink}[1]{\expandafter\gdef\csname @hreflink\endcsname{#1}}
\providecommand{\issuenum}[1]{\expandafter\gdef\csname @issuenum\endcsname{#1}}
\providecommand{\pubvolume}[1]{\expandafter\gdef\csname @pubvolume\endcsname{#1}}
\providecommand{\pubyear}[1]{\expandafter\gdef\csname @pubyear\endcsname{#1}}
\firstpage{1}
\pubvolume{1}
\issuenum{1}
\articlenumber{0}
\pubyear{2026}
\copyrightyear{2026}
\datereceived{ }
\daterevised{ }
\dateaccepted{ }
\datepublished{ }
\hreflink{https://doi.org/}

\NewDocumentEnvironment{theorem}{o}
  {\IfValueTF{#1}{\begin{Theorem}[#1]}{\begin{Theorem}}}{\end{Theorem}}
\NewDocumentEnvironment{lemma}{o}
  {\IfValueTF{#1}{\begin{Lemma}[#1]}{\begin{Lemma}}}{\end{Lemma}}
\NewDocumentEnvironment{corollary}{o}
  {\IfValueTF{#1}{\begin{Corollary}[#1]}{\begin{Corollary}}}{\end{Corollary}}
\NewDocumentEnvironment{proposition}{o}
  {\IfValueTF{#1}{\begin{Proposition}[#1]}{\begin{Proposition}}}{\end{Proposition}}
\NewDocumentEnvironment{definition}{o}
  {\IfValueTF{#1}{\begin{Definition}[#1]}{\begin{Definition}}}{\end{Definition}}
\NewDocumentEnvironment{example}{o}
  {\IfValueTF{#1}{\begin{Example}[#1]}{\begin{Example}}}{\end{Example}}
\NewDocumentEnvironment{remark}{o}
  {\IfValueTF{#1}{\begin{Remark}[#1]}{\begin{Remark}}}{\end{Remark}}
\NewDocumentEnvironment{hypothesis}{o}
  {\IfValueTF{#1}{\begin{Hypothesis}[#1]}{\begin{Hypothesis}}}{\end{Hypothesis}}
\NewDocumentEnvironment{assumption}{o}
  {\IfValueTF{#1}{\begin{Assumption}[#1]}{\begin{Assumption}}}{\end{Assumption}}

\usepackage[utf8]{inputenc}
\usepackage[novbox]{pdfsync}   

\makeatletter
\@ifundefined{AAdefNnogeometry}{\usepackage[a4paper,margin=2cm]{geometry}}{}
\makeatother

\usepackage{amsmath,amssymb,amsthm,mathtools,bm,mathrsfs,bbm}
\usepackage{centernot}

\usepackage{graphicx,float}
\makeatletter
\@ifpackageloaded{subfig}{}{\usepackage{subcaption}}
\makeatother
\usepackage{booktabs,longtable,threeparttable,array,multirow}
\usepackage{enumitem}

\usepackage{xcolor,listings,tcolorbox}
\makeatletter
\@ifundefined{c@algorithm}{\usepackage{algorithm}}{}
\makeatother
\usepackage{algorithmicx,algpseudocode}

\usepackage{tikz,tikz-cd,pgfplots}
\usetikzlibrary{graphs,decorations.pathmorphing,decorations.markings,
arrows.meta,automata,positioning}
\tikzcdset{every label/.append style = {font = \small}}
\pgfplotsset{compat=1.18}

\usepackage[version=4]{mhchem}

\usepackage{microtype}
\usepackage{lineno}
\usepackage{nomencl}
\makenomenclature

\makeatletter
\@ifundefined{hyperrefpackage}{%
  \@ifundefined{AAdefNnohyperref}{\usepackage{hyperref}}{}%
}{}
\@ifundefined{hyperrefpackage}{\usepackage{bookmark}}{}
\makeatother
\usepackage[toc,sort=use]{glossaries}
\makeglossaries

\newif\ifbookbuild
\bookbuildfalse

\newif\ifmadec     \madectrue
\newif\ifslowfast  \slowfasttrue

\newif\ifhnetpart  \hnetparttrue

\RequirePackage{etoolbox}
\usepackage{makecell}
\usepackage[novbox]{pdfsync}
\usepackage{amssymb, amsfonts,mathtools,bbm,eurosym,graphicx,latexsym, amsmath,amsthm,times,nccmath,url}
\usepackage{listings,verbatim,cases,tabularx}
\makeatletter
\@ifpackageloaded{hyperref}{}{%
  \@ifundefined{hyperrefpackage}{%
    \@ifundefined{AAdefNnohyperref}{\usepackage{hyperref}}{}%
  }{}%
}
\makeatother
\makeatletter
\AtBeginDocument{\@ifundefined{hypersetup}{}{\hypersetup{linktoc=all}}}
\makeatother
\usepackage{graphicx,float,fancybox}
\usepackage{epsfig,epstopdf,tikz}
\usepackage[sans]{dsfont}
\usetikzlibrary{positioning,arrows,arrows.meta,calc}
\lstdefinestyle{mathematica}{
 language=Mathematica,
 basicstyle=\small\ttfamily,
 backgroundcolor=\color{gray!10},
 frame=single,
 breaklines=true,
 commentstyle=\color{green!60!black},
 keywordstyle=\color{blue},
 stringstyle=\color{red},
 showstringspaces=false,
 tabsize=2
}
\usepackage{tcolorbox}
\newtcolorbox{keyresult}[1][]{
 colback=blue!5!white,
 colframe=blue!75!black,
 title=#1,
 fonttitle=\bfseries
}
\newtcolorbox{examplebox}[1][]{
 colback=green!5!white,
 colframe=green!75!black,
 title=#1,
 fonttitle=\bfseries
}

\makeatletter
\@ifundefined{theorem}   {\newtheorem{theorem}   {Theorem}}   {}
\@ifundefined{definition}{\newtheorem{definition}{Definition}} {}
\@ifundefined{open}      {}   {}
\@ifundefined{example}   {\newtheorem{example}   {Example}}   {}
\makeatother

\newcounter{tabcounter}
\newcounter{defcounter}

\newcounter{thmcounter}
\newcounter{excounter}
\newcounter{ascounter}
\newcounter{cjcounter}

\newcounter{figcounter}

\newglossarystyle{unifiedstyle}{%
  \setglossarystyle{list}%
  \renewcommand{\glossentry}[2]{%
    \edef\glscursec{\glsentryuseriii{##1}}%
    \ifx\glscursec\glslastsec\else
      \item[]\textbf{\S\,\glscursec\enspace\glsentryuseriv{##1}}%
      \global\let\glslastsec\glscursec
    \fi
    \item[\glsentryuserii{##1} \glsentryuseri{##1}:] ``\glossentryname{##1}'' \glsentrydesc{##1}%
  }%
}

\makeatletter
\def\@firstwordof#1 #2\relax{#1}
\def\@takethree#1#2#3#4\@nil{#1#2#3}

\NewDocumentCommand{\beXa}{s O{} O{} o}{%
  \def\@temparg{#2}%
  \ifx\@temparg\@empty
    \begin{example}%
  \else
    \let\@savedSecTitle\@currentsectionname%
    \refstepcounter{excounter}%
    \IfValueTF{#4}{%
      \edef\@exlblkey{#4}%
    }{%
      \edef\@exfirstword{\@firstwordof#2 \relax}%
      \edef\@exlblkey{%
        \expandafter\@takethree\@exfirstword\@empty\@empty\@empty\@nil}%
    }%
    \def\@tempshort{#3}%
    \ifx\@tempshort\@empty
      \expandafter\xdef\csname @xshort@\@exlblkey\endcsname{\@exlblkey}%
    \else
      \expandafter\def\csname @xshort@\@exlblkey\endcsname{#3}%
    \fi
    \IfBooleanTF{#1}{\let\@extag\@empty}{%
      \edef\@extag{\noexpand\ \noexpand\textup{[\noexpand\texttt{e:\@exlblkey}]}}}%
    \begin{example}[#2\ifshowkeys\@extag\fi]%
    \label{e:\@exlblkey}%
    \protected@write\@auxout{}{%
      \string\@thmtitle{e:\@exlblkey}{\unexpanded{#2}}}%
    \expanded{%
      \noexpand\newglossaryentry{e\theexcounter}{%
        name={\unexpanded{#2}%
          \IfBooleanF{#1}{ [\noexpand\texttt{e:\@exlblkey}]}},%
        sort={\theexample},%
        description={p.~\noexpand\pageref{e:\@exlblkey}},%
        user1={\theexample},%
        user2={Example},%
        user3={\thesection},%
        user4={\expandafter\unexpanded\expandafter{\@savedSecTitle}}%
      }%
    }%
    \protected@write\@auxout{}{\string\glsliveentry{e\theexcounter}}%
  \fi
}
\newcommand{\eeXa}{\end{example}}

\NewDocumentCommand{\beD}{s O{} O{} o}{%
  \def\@temparg{#2}%
  \ifx\@temparg\@empty
    \begin{definition}%
  \else
    \let\@savedSecTitle\@currentsectionname%
    \refstepcounter{defcounter}%
    \IfValueTF{#4}{%
      \edef\@deflblkey{#4}%
    }{%
      \edef\@deffirstword{\@firstwordof#2 \relax}%
      \edef\@deflblkey{%
        \expandafter\@takethree\@deffirstword\@empty\@empty\@empty\@nil}%
    }%
    \def\@tempshort{#3}%
    \ifx\@tempshort\@empty
      \expandafter\xdef\csname @dshort@\@deflblkey\endcsname{\@deflblkey}%
    \else
      \expandafter\def\csname @dshort@\@deflblkey\endcsname{#3}%
    \fi
    \IfBooleanTF{#1}{\let\@deftag\@empty}{%
      \edef\@deftag{\noexpand\ \noexpand\textup{[\noexpand\texttt{d:\@deflblkey}]}}}%
    \begin{definition}[#2\ifshowkeys\@deftag\fi]%
    \label{d:\@deflblkey}%
    \protected@write\@auxout{}{%
      \string\@thmtitle{d:\@deflblkey}{\unexpanded{#2}}}%
    \expanded{%
      \noexpand\newglossaryentry{d\thedefcounter}{%
        name={\unexpanded{#2}%
          \IfBooleanF{#1}{ [\noexpand\texttt{d:\@deflblkey}]}},%
        sort={\thedefinition},%
        description={p.~\noexpand\pageref{d:\@deflblkey}},%
        user1={\thedefinition},%
        user2={Definition},%
        user3={\thesection},%
        user4={\expandafter\unexpanded\expandafter{\@savedSecTitle}}%
      }%
    }%
    \protected@write\@auxout{}{\string\glsliveentry{d\thedefcounter}}%
  \fi
}
\newcommand{\eeD}{\end{definition}}

\newif\ifshowkeys \showkeystrue
\newif\ifcommonred \commonredfalse
\NewDocumentCommand{\beT}{s O{} O{} o}{%
  \def\@temparg{#2}%
  \ifx\@temparg\@empty
    \begin{theorem}%
  \else
    \let\@savedSecTitle\@currentsectionname%
    \refstepcounter{thmcounter}%
    \IfValueTF{#4}{%
      \edef\@thmlblkey{#4}%
    }{%
      \edef\@thmfirstword{\@firstwordof#2 \relax}%
      \edef\@thmlblkey{%
        \expandafter\@takethree\@thmfirstword\@empty\@empty\@empty\@nil}%
    }%
    \IfBooleanTF{#1}{\let\@thmtag\@empty}{%
      \edef\@thmtag{\noexpand\ \noexpand\textup{[\noexpand\texttt{t:\@thmlblkey}]}}}%
    \def\@tempshort{#3}%
    \ifx\@tempshort\@empty
      \expandafter\xdef\csname @tshort@\@thmlblkey\endcsname{\@thmlblkey}%
    \else
      \expandafter\def\csname @tshort@\@thmlblkey\endcsname{#3}%
    \fi
    \begin{theorem}[#2\ifshowkeys\@thmtag\fi]%
    \label{t:\@thmlblkey}%
    \protected@write\@auxout{}{%
      \string\@thmtitle{t:\@thmlblkey}{\unexpanded{#2}}}%
    \expanded{%
      \noexpand\newglossaryentry{t\thethmcounter}{%
        name={\unexpanded{#2}%
          \IfBooleanF{#1}{ [\noexpand\texttt{t:\@thmlblkey}]}},%
        sort={\thetheorem},%
        description={p.~\noexpand\pageref{t:\@thmlblkey}},%
        user1={\thetheorem},%
        user2={Theorem},%
        user3={\thesection},%
        user4={\expandafter\unexpanded\expandafter{\@savedSecTitle}}%
      }%
    }%
    \protected@write\@auxout{}{\string\glsliveentry{t\thethmcounter}}%
  \fi
}
\newcommand{\eeT}{\end{theorem}}

\newcommand{\xpart}[1]{\@ifundefined{r@part:glv}{}{\paragraph{Relation with Parts~I--II.}#1}}

\newcommand{\bookref}[3]{\@ifundefined{r@#1}{#3}{#2}}

\newcommand{\beAs}[1][]{%
  \def\@temparg{#1}%
  \ifx\@temparg\@empty
    \begin{assumption}%
  \else
    \let\@savedSecTitle\@currentsectionname%
    \refstepcounter{ascounter}%
    \begin{assumption}[#1]%
    \label{as:\theascounter}%
    \expanded{%
      \noexpand\newglossaryentry{as\theascounter}{%
        name={\unexpanded{#1}},%
        sort={\theascounter},%
        description={p.~\noexpand\pageref{as:\theascounter}},%
        user1={\theascounter},%
        user2={Assumption},%
        user3={\thesection},%
        user4={\expandafter\unexpanded\expandafter{\@savedSecTitle}}%
      }%
    }%
    \protected@write\@auxout{}{\string\glsliveentry{as\theascounter}}%
  \fi
}
\newcommand{\eeAs}{\end{assumption}}
\newcommand{\beCj}[1][]{%
  \def\@temparg{#1}%
  \ifx\@temparg\@empty
    \begin{conjecture}%
  \else
    \let\@savedSecTitle\@currentsectionname%
    \refstepcounter{cjcounter}%
    \begin{conjecture}[#1]%
    \label{cj:\thecjcounter}%
    \expanded{%
      \noexpand\newglossaryentry{cj\thecjcounter}{%
        name={\unexpanded{#1}},%
        sort={\thecjcounter},%
        description={p.~\noexpand\pageref{cj:\thecjcounter}},%
        user1={\thecjcounter},%
        user2={Conjecture},%
        user3={\thesection},%
        user4={\expandafter\unexpanded\expandafter{\@savedSecTitle}}%
      }%
    }%
    \protected@write\@auxout{}{\string\glsliveentry{cj\thecjcounter}}%
  \fi
}
\newcommand{\eeCj}{\end{conjecture}}
\newcommand{\@mkresult}[5]{%
  \expandafter\NewDocumentCommand\csname be#1\endcsname{s O{} O{} o}{%
    \def\@temparg{##2}%
    \ifx\@temparg\@empty
      \begin{#2}%
    \else
      \let\@savedSecTitle\@currentsectionname%
      \refstepcounter{#4}%
      \IfValueTF{##4}{%
        \edef\@relbl{##4}%
        \IfBooleanTF{##1}{\let\@retag\@empty}{%
          \edef\@retag{\noexpand\ \noexpand\textup{[\noexpand\texttt{#3:##4}]}}}%
      }{%
        \edef\@relbl{\arabic{#4}}%
        \let\@retag\@empty
      }%
      \def\@tempshort{##3}%
      \ifx\@tempshort\@empty
        \expandafter\xdef\csname @sshort@#3\@relbl\endcsname{\@relbl}%
      \else
        \expandafter\def\csname @sshort@#3\@relbl\endcsname{##3}%
      \fi
      \begin{#2}[##2\ifshowkeys\@retag\fi]%
      \edef\@tmplbl{#3:\@relbl}\expandafter\label\expandafter{\@tmplbl}%
      \expanded{%
        \noexpand\newglossaryentry{#3\csname the#4\endcsname}{%
          name={\unexpanded{##2}%
            \IfBooleanF{##1}{\IfValueT{##4}{ [\noexpand\texttt{#3:##4}]}}},%
          sort={\csname the#2\endcsname},%
          description={p.~\noexpand\pageref{#3:\@relbl}},%
          user1={\csname the#2\endcsname},%
          user2={#5},%
          user3={\thesection},%
          user4={\expandafter\unexpanded\expandafter{\@savedSecTitle}}%
        }%
      }%
      \protected@write\@auxout{}{\string\glsliveentry{#3\csname the#4\endcsname}}%
      \protected@write\@auxout{}{\string\@thmtitle{\@tmplbl}{\unexpanded{##2}}}%
    \fi
  }%
  \expandafter\newcommand\csname ee#1\endcsname{\end{#2}}%
}
\@mkresult{L}{lemma}{l}{lemcounter}{Lemma}
\@mkresult{R}{remark}{r}{remcounter}{Remark}
\@mkresult{C}{corollary}{c}{corcounter}{Corollary}
\@mkresult{P}{proposition}{p}{propcounter}{Proposition}
\@mkresult{O}{open}{o}{opcounter}{Problem}

\newcommand{\@currentsectionname}{}
\AtBeginDocument{%
  \let\@epidcrn@orig@refstepcounter\refstepcounter
  \renewcommand{\refstepcounter}[1]{%
    \@epidcrn@orig@refstepcounter{#1}%
    \ifstrequal{#1}{section}{\global\let\@currentsectionname\@currentlabelname}{}%
  }%
}

\newcommand{\beTa}[1][]{%
  \def\@temparg{#1}%
  \ifx\@temparg\@empty
    \begin{table}[htbp]%
  \else
    \let\@savedSecTitle\@currentsectionname%
    \refstepcounter{tabcounter}%
    \begin{table}[htbp]%
    \label{A:\thetabcounter}%
    \expanded{%
      \noexpand\newglossaryentry{A\thetabcounter}{%
        name={\unexpanded{#1}},%
        sort={\thetabcounter},%
        description={p.~\noexpand\pageref{A:\thetabcounter}},%
        user1={\thetabcounter},%
        user2={Table},%
        user3={\thesection},%
        user4={\expandafter\unexpanded\expandafter{\@savedSecTitle}}%
      }%
    }%
    \protected@write\@auxout{}{\string\glsliveentry{A\thetabcounter}}%
  \fi
}
\newcommand{\eeTa}{\end{table}}

\newcommand{\beF}[1][]{%
  \def\@temparg{#1}%
  \ifx\@temparg\@empty
    \begin{figure}[H]%
  \else
    \begin{figure}[H]%
    \let\@savedSecTitle\@currentsectionname%
    \refstepcounter{figcounter}%
    \label{fig:\thefigcounter}%
    \expanded{%
      \noexpand\newglossaryentry{Fig\thefigcounter}{%
        name={\unexpanded{#1}},%
        sort={\thefigcounter},%
        description={p.~\noexpand\pageref{fig:\thefigcounter}},%
        user1={\thefigcounter},%
        user2={Figure},%
        user3={\thesection},%
        user4={\expandafter\unexpanded\expandafter{\@savedSecTitle}}%
      }%
    }%
    \protected@write\@auxout{}{\string\glsliveentry{Fig\thefigcounter}}%
  \fi
}
\newcommand{\eeF}{\end{figure}}
\makeatother
\makeatletter
\newcommand{\@glslive}{}
\newcommand{\@glsaddlive}[1]{\ifglsentryexists{#1}{\glsadd{#1}}{}}
\newcommand{\glsliveentry}[1]{\g@addto@macro\@glslive{\@glsaddlive{#1}}}
\makeatother
\makeatletter
\newcommand{\printboth}{%
  \@glslive   
  \printglossary[style=unifiedstyle, title={0 \hspace{.22em} Definitions, Theorems, Examples and Open Problems}]%
}
\makeatother
\newcommand{\tref}[1]{%
  Thm~\ref{t:#1}(\ifcsname @tshort@#1\endcsname\csname @tshort@#1\endcsname\else#1\fi)%
}

\makeatletter
\newcommand{\@thmtitle}[2]{\global\@namedef{@tlong@#1}{#2}}
\newcommand{\nref}[1]{%
  \ref{#1}%
  \ifcsname @tlong@#1\endcsname
    \ ``\csname @tlong@#1\endcsname''%
  \fi
}
\newcommand{\nm}[1]{%
  \ifcsname @tlong@#1\endcsname
    \csname @tlong@#1\endcsname%
  \else
    \ref{#1}%
  \fi
}
\makeatother

\makeatletter
\@ifundefined{lemma}      {\newtheorem{lemma}      {Lemma}}      {}
\@ifundefined{proposition}{\newtheorem{proposition}{Proposition}} {}
\@ifundefined{corollary}  {\newtheorem{corollary}  {Corollary}}  {}
\@ifundefined{conjecture} {\newtheorem{conjecture} {Conjecture}}  {}
\@ifundefined{remark}     {\newtheorem{remark}     {Remark}}     {}
\@ifundefined{question}   {\newtheorem{question}   {Question}}   {}
\@ifundefined{assumption} {\newtheorem{assumption} {Assumption}} {}
\@ifundefined{hypothesis} {\newtheorem{hypothesis} {Hypothesis}} {}
\@ifundefined{fact}       {}       {}
\makeatother

\makeatletter
\providecommand{\citer}{%
  \@ifundefined{r@t:relayRF}%
    {\cite{AAHrel}}%
    {Theorem~\ref{t:relayRF} of \cite{AAHrel}}%
}
\providecommand{\citeR}{%
  \@ifundefined{r@t:Pstar}%
    {\cite{AAHRah}}%
    {Theorem~\ref{t:Pstar} of \cite{AAHRah}}%
}
\makeatother

\def\beQ{\begin{question}}\def\eeQ{\end{question}}
\def\beA{\begin{assumption}}\def\eeA{\end{assumption}}

\definecolor{funccolor}{RGB}{25,25,112}
\definecolor{desccolor}{RGB}{64,64,64}

\def\bep{\begin{pmatrix}}\def\eep{\end{pmatrix}}
\def\bev{\begin{vmatrix}}\def\eev{\end{vmatrix}}
\def\bea{\begin{eqnarray*}}\def\eea{\end{eqnarray*}}
\def\bc{\begin{cases}}\def\ec{\end{cases}}
\def\BEN{\begin{enumerate}}\def\EEN{\end{enumerate}}
\def\BI{\begin{itemize}}\def\EI{\end{itemize}}

\newcommand{\be}[1]{\begin{equation}\label{#1}}
\newcommand{\ee}{\end{equation}}
\newcommand{\beq}{\begin{eqnarray}}
\def\eeq{\end{eqnarray}}
\def\eqr{\eqref}\def\lbl{\label}

\newcommand{\ind}[1]{\ensuremath{\mathbbm{1}_{\left\{#1\right\}}}}

\newcommand{\R}{\mathbb{R}}  

\def\T{\widetilde}

\newcommand{\bff}[1]{{\mbox{\boldmath$#1$}}}

\def\v1{\vec {\bff 1}}

\DeclareMathOperator{\diag}{Diag}

\long\def\symbolfootnote[#1]#2{%
\begingroup
\def\thefootnote{\fnsymbol{footnote}}\footnote[#1]{#2}%
\endgroup}

\def\im{\item}

\def\ssec{\subsection}

\newcommand\CRN{chemical reaction network}

\def\DFE{disease free equilibrium}

\def\NGM{next generation matrix}

\makeatletter
\newcommand{\ideaitem}[2]{%
  \expandafter\gdef\csname @idea@#1\endcsname{#2}%
  \item\hyperlink{idea:#1}{#2}}
\newcommand{\ideahead}[2][]{%
  \item\hypertarget{idea:#2}{}%
  \def\@ideatmp{#1}%
  \ifx\@ideatmp\@empty
    \@ifundefined{@idea@#2}{%
      \PackageWarning{preamble}{Idea #2 is used in Main Results but was never
        defined by \string\ideaitem\space in the abstract}%
      \textbf{[idea #2 undefined]}%
    }{\textbf{\begingroup\def\ideaterm##1##2{##2}\def\idearest##1{}\@nameuse{@idea@#2}\endgroup}}%
  \else
    \textbf{#1}%
  \fi}
\newcommand{\ideaterm}[2]{\hyperlink{#1}{#2}}
\newcommand{\idearest}[1]{#1}
\makeatother

\newif\ifshort \shortfalse

\newcommand{\Jtr}{J^\perp}
\newcommand{\JtrS}[1]{J_{#1}^\perp}
\showkeysfalse
\usetikzlibrary{arrows.meta,positioning}

\AddToHook{begindocument/end}{\hypersetup{hidelinks}}

\definecolor{revcolor}{rgb}{0.00,0.22,0.70}
\definecolor{opencolor}{rgb}{0.62,0.07,0.07}
\definecolor{newcolor}{rgb}{0.00,0.38,0.14}
\ifmarked\else\colorlet{revcolor}{black}\colorlet{newcolor}{black}\fi
\ifmarked\long\def\rev#1{\textcolor{revcolor}{#1}}\else\long\def\rev#1{#1}\fi
\ifmarked\newenvironment{revblock}{\color{revcolor}\ignorespaces}{\par}\else\fi
\ifmarked\newenvironment{unrefblock}{\color{newcolor}\ignorespaces}{\par}\else\fi
\AtBeginDocument{\pdfstringdefDisableCommands{%
  \long\def\rev#1{#1}%
}}
\Title{Threshold Networks: Global Convergence, Arbitrary Multistability, and Cusp Bifurcations}

\Author{Florin Avram}

\AuthorNames{Florin Avram}

\address{Laboratoire de Math\'{e}matiques Appliqu\'{e}es, Universit\'{e} de Pau, Pau, France;
avramf3@gmail.com}

\abstract{We consider positive production and threshold networks with nonnegative monotone responses. For the positive Lur'e family $\dot x=-Vx+B\Phi(\beta Wx)$, with $V$ a nonsingular $M$-matrix, $B,W\ge0$, $\beta\ge0$ and a common nonnegative response with $0\le\Phi'\le L$, the threshold $\beta\rho(BWV^{-1})=1/L$ separates two response-class statements: below it every admissible response gives a unique globally exponentially attracting equilibrium, whereas above it a smooth bounded response can be constructed with at least any prescribed finite number of locally exponentially stable equilibria. The contraction direction is classical; in the construction the response depends on the requested count. Vanishing production and absence of linear influx from outside the face characterize invariant siphon faces; the first-order production terms and linear removal determine their Metzler transverse Jacobians. For fixed nonnegative production polynomials satisfying a support condition, common smooth responses with flat onset yield arbitrarily many stable positive equilibria while preserving the siphons and their transverse matrices. For odds-polynomial responses, degree four is necessary and sufficient for a stable origin to coexist with a positive nondegenerate cusp. A three-dimensional construction with overlapping minimal siphons permits independently variable invasion numbers at two nonzero boundary residents and an interior cusp coexisting with at least eight stable equilibria obtained by weak coupling of bistable coordinates. Invariant partitions provide exact reduced equations and lift nondegenerate bifurcations under stable complementary directions. Positive diagonal diffusion preserves the constructed stable constant equilibria and yields a linear-comparison upper bound on spreading; the nonlinear propagation speed requires separate analysis.}

\keyword{threshold networks; monotone systems; $M$-matrices; multistability; cusp bifurcation; siphons; spreading bounds}

\MSC{34D20; 34D23; 34C23; 37N25; 92B20; 92D30}

\providecommand{\bookref}[3]{#3}

\providecommand{\JtrS}[1]{J^{\perp}_{#1}}
\providecommand{\Jtr}{J^{\perp}}

\begin{document}
\MDPIfrontmatter
\tableofcontents

\section{Introduction to production and threshold models with response functions}
\label{s:int}

\beD[Production and threshold models][][prMod]\label{d:KB}
A \emph{production model} is a system on a forward-invariant region of the nonnegative orthant which decomposes as
\begin{equation}\label{eq:KB}
 \dot x=f(x,t)=\Pi(x,t)-Vx,
\end{equation}
where $\Pi\ge0$ and $V$ is a nonsingular $M$-matrix. The state has $n_S$ coordinates. If $\Pi$ is nondecreasing in $x$ and differentiable, the Jacobian $\partial_x\Pi-V$ is Metzler; more generally, nonnegative off-diagonal production derivatives suffice.

A \emph{threshold response} is a nonnegative nondecreasing function $\Phi$
(boundedness and smoothness assumptions are stated in each result),
which is applied to response arguments. In the polynomial production subclass, the response is applied coordinatewise:
\begin{equation}\label{polR}
 \Pi_i(x)=b_i\Phi_i\bigl(h_i+P_i(x)\bigr),\qquad b_i>0,
\end{equation}

Here $P_i(x)$ are polynomials with nonnegative coefficients (with $P_i(0)=0$ after absorbing the constant term into $h_i$); the linear case $P_i(x)=\beta(Wx)_i$, $W\ge0$, is a subclass.

\eeD

\ssec{Some subclasses of production models}\lbl{s:clas}
{\bf  Our common production--loss framework}  \eqref{eq:KB} unifies several important subclasses, distinguished by the form of $\Pi$. The classes overlap; Class 3 allows production that need not be nondecreasing. Remark~\ref{r:genM} also allows nonlinear removal.

\BEN
\im {\bf Polynomial production  models:}
With $P_i$ in \eqr{polR} a polynomial with nonnegative coefficients, monomials such as $x_1x_3$ represent higher-order interactions, and the
choice of monomials decides which siphons survive.
\emph{Studied in} Theorems~\ref{t:siphMulti} (multistability with the siphons
fixed), \ref{t:Net} and~\ref{t:KBoverlapCusp} (interior cusps over stable
siphon faces) and~\ref{t:OvlRelay} (backward relay); the scalar response
these three use is the one of Theorem~\ref{t:Qua}.
\im
{\bf (Positive) Lur'e systems:}
\begin{equation}\label{eq:LureClass}
 \dot x=-Vx+B\Phi(h+\beta Wx,t),
 \qquad B,W\ge0,
\end{equation}
where $\Phi(y,t)$ is coordinatewise nonnegative, nondecreasing
and $C^1$ in $y$. Contraction results for Lur'e systems are
found for example in
\cite{Davydov,davydov2024exponential,OfirOvseevichMargaliot2024,OfirSlotineMargaliot2025}.
\emph{Studied in} Theorem~\ref{t:KBsharp}, which, for autonomous systems
with a common scalar response, zero stimulus and a slope bound
$0\le\Phi'\le L$, establishes the sharp threshold $\bar R_0=1/L$, where
\begin{equation}\label{eq:LureR0}
 \bar R_0=\beta\rho(V^{-1}BW),
\end{equation}
and in Theorem~\ref{t:SpdBd}, which carries both conclusions to diagonal
diffusion. The threshold depends on the product $BW$, whereas the nonlinear
production generally depends on $B$ and $W$ separately. For any number $n_R$
of response coordinates one may take $B\in\mathbb R_+^{n_S\times n_R}$ and
$W\in\mathbb R_+^{n_R\times n_S}$. The following nested subclasses are of particular
interest.
\BEN
\im {\bf Linear production threshold models:}
$B=\operatorname{Diag}(b_i)$, $b_i>0$, with autonomous responses. These also
belong to the polynomial class above, since the polynomials in \eqr{polR}
are then linear. The scalar case is $\dot x=-x+\Phi(h+\beta x)$, where the
cumulative distribution function (CDF)
\begin{equation}\label{eq:response}
 \Phi(q)=\mathbb P(\Theta\le q)
\end{equation}
gives the proportion of individuals whose threshold does not exceed $q$. The
network case with $B=V=I$ and $W\mathbf1=\mathbf1$ is
$\dot x_i=-x_i+\Phi(h+\beta(Wx)_i)$.
\im {\bf Firing-rate networks:} the subclass of the previous one with
$B=I$, $V$ positive diagonal and $W\ge0$:
$\dot x=-Vx+\Phi(Wx+h)$. For $V=I$ the variable $y=Wx+h$ satisfies the
Hopfield equation $\dot y=-y+W\Phi(y)+h$, which again belongs to the general
Lur'e class. Firing-rate networks are a central model of neural population
dynamics \cite{Betteti2025,Davydov,Gokhale2026}; multiple cusp bifurcations
in neural networks were studied in \cite{Izhikevich1998}, and the Hopfield
equation is \cite{Hopfield1984}.
\EEN

\im {\bf Models with a state-dependent amplification:}
\begin{equation}\label{eq:KB-general-response}
 \Pi_i(x)=b_i(x_i)\Phi_i\bigl(h_i+P_i(x)\bigr),\qquad b_i\ge0\ \text{of class }C^1,
\end{equation}
which is \eqref{polR} with the constant amplification $b_i$ replaced by a
function of $x_i$.
The main case is a \emph{capacity factor} $1-\rho_ix_i$, $\rho_i>0$, which decreases from one to zero as $x_i$ increases to its
capacity $1/\rho_i$ (for instance a susceptible fraction), possibly multiplying a production $G_i\ge0$ of no specified form:
\begin{equation}\label{eq:KB-capacity}
 \Pi_i(x)=(1-\rho_ix_i)G_i(x),\qquad x\in[0,1/\rho],\quad V(1/\rho)\ge0.
\end{equation}
Then $\Pi\ge0$ holds only on this box, whose upper faces point inward, since $f_i(x)\le-(V(1/\rho))_i\le0$ at $x_i=1/\rho_i$, and
$\Pi_i$ need not be nondecreasing in $x_i$. If also $\partial_jG_i\ge0$ for $j\ne i$, the Jacobian is Metzler on the box.
Higher-order SIS networks and simplicial contagion \cite{IPB} are
\begin{equation}\label{eq:SIS}
 \dot x_i=-\mu_ix_i+(1-x_i)
 \Bigl(\sum_ja_{ij}x_j+\sum_{|E|\ge2}\beta_{iE}\prod_{j\in E}x_j\Bigr),
 \qquad \mu_i>0,\quad a_{ij},\beta_{iE}\ge0,
\end{equation}
on $[0,1]^{n_S}$, with the susceptible fraction $1-x_i$ as capacity factor. Their scalar form is
\begin{equation}\label{eq:simp}
 \dot u=-u+\lambda u(1-u)+\lambda_\Delta u^2(1-u).
\end{equation}
In epidemiology such models exhibit backward bifurcation \cite{ArinoMc}. The excitatory Wilson--Cowan equation
\cite{WilsonCowan1972} is $\tau\dot E=-E+(k-rE)S(cE+P)$, where distinct response functions represent differences between
excitatory and inhibitory populations; an example with nonlinear removal, covered by Remark~\ref{r:genM}, is logistic growth with saturating loss, $\dot x=rx(1-x/K)-cx^2/(x^2+h^2)$.
\emph{Studied in} Theorem~\ref{t:InvF} (the face criterion and the transverse
matrices, under $b_i(0)>0$), Proposition~\ref{p:KBshared} (shared responses,
with $b_i$ constant in the contraction test) and Remark~\ref{r:KBtransG} (the
transverse matrix of the capacity-factor and SIS forms). No multistability or
cusp construction of this paper is carried out for this class.

\EEN

Table~\ref{tab:scope} lists the principal theorems. We claim as new the constructions in Theorems~\ref{t:siphMulti}, \ref{t:KBsharp}(part 2), \ref{t:Qua}, \ref{t:Net}, \ref{t:KBoverlapCusp} and \ref{t:OvlRelay}; the remaining results are developments of established methods.

\beTa[Principal theorems]
\centering
\small
\setlength{\tabcolsep}{4pt}
\renewcommand{\arraystretch}{1.15}
\begin{tabular}{p{0.13\linewidth}p{0.58\linewidth}p{0.23\linewidth}}
\hline
Theorem & Title & Subclass\\
\hline
\multicolumn{3}{l}{\emph{Claimed new}}\\
Theorem~\ref{t:siphMulti} & \nm{t:siphMulti} & Polynomial production\\
Theorem~\ref{t:KBsharp} & \nm{t:KBsharp} & Autonomous positive Lur'e\\
Theorem~\ref{t:Qua} & \nm{t:Qua} & Scalar odds response\\
Theorem~\ref{t:Net} & \nm{t:Net} & Polynomial production, odds response\\
Theorem~\ref{t:KBoverlapCusp} & \nm{t:KBoverlapCusp} & Polynomial
production, odds response\\
Theorem~\ref{t:OvlRelay} & \nm{t:OvlRelay} & The network of
Theorem~\ref{t:KBoverlapCusp}\\
\hline
\multicolumn{3}{l}{\emph{Obtained by established methods}}\\
Theorem~\ref{t:InvF} & \nm{t:InvF} & General production models\\
Theorem~\ref{t:Spe} & \nm{t:Spe} & Any vector field with an invariant
partition\\
Theorem~\ref{t:SpdBd} & \nm{t:SpdBd} & Autonomous positive Lur'e, with
diffusion\\
\hline
\end{tabular}
\caption{Principal theorems; the third column gives the subclass of Section~\ref{s:clas} each one
is stated for.}\lbl{tab:scope}
\eeTa


\beR[Where the threshold responses come from][][histhr]
The probabilistic form \eqref{eq:response} has early roots in Thurstone's
comparative judgment \cite{Thurstone1927}; related foundations are the
all-or-none units of McCulloch and Pitts \cite{McCulloch1943}, Schelling's
models of segregation and tipping \cite{Schelling1971}, and Glauber's spin
dynamics \cite{Glauber1963}. In Granovetter's model \cite{granovetter1978}
individuals act once the active fraction exceeds a personal threshold; his
map $r(t+1)=\Phi[r(t)]$ has the equilibria of $\dot x=-x+\Phi(x)$. In Kaye's
dawn-chorus model \cite{Kaye26}, light and neighbouring songs promote
singing, and Kaye relates folds and cusps to the threshold density;
collective animal choruses are also modelled in
\cite{Aihara2014,HartbauerRomer2016}.  The network case of the linear production threshold models is studied by
Boudourides \cite{Boud26}: homogeneous branches, reductions, multistability,
heterogeneous cusps, and the persistence of families of stable equilibria
under weak interblock coupling. Network versions show how local exposure and
topology control { \bf global cascades} \cite{Watts2002}, and smooth and stochastic
variants produce saddle-node transitions and hysteresis
\cite{Wiedermann2020,Kook2021}. Scheffer et al.\ \cite{Scheffer2001} analyse
 {\bf catastrophic shifts} through $\dot x=a-bx+rx^p/(x^p+h^p)$, whose response
includes the constant term $a$; see also \cite{Scheffer2009,Kuehn2011}, and
\cite{kroenke2020} for tipping cascades between coupled subsystems.

\eeR

\beR[Family containing all of the above][][genM]
Consider
\begin{equation}\label{eq:intro-family}
 \Pi_i(x)=(1-\rho_ix_i)\Bigl[\eta_i+\sum_{k}b_{ik}\Phi_{ik}\bigl(h_{ik}+P_{ik}(x)\bigr)\Bigr],
 \qquad P_{ik}(x)=\sum_{\nu\ne0}a_{ik\nu}x^\nu,\quad a_{ik\nu}\ge0,
\end{equation}
with $\rho_i,\eta_i\ge0$, $b_{ik}>0$, finitely many monomials, and nonnegative nondecreasing $C^1$ responses (bounded distribution
functions and the identity are allowed). The removal may be nonlinear, $Vx+L(x)$, with $L_i(x_i)$ nondecreasing $C^1$ and
$L_i(0)=0$ (grazing, harvesting, saturating loss), so that $\dot x=-Vx-L(x)+\Pi(x)$. This combines threshold responses
($\rho=L=0$), pairwise or product terms, capacity factors, and saturating losses. Responses, removal coefficients and
interaction strengths may vary by coordinate; activity levels, fractions or concentrations need not all lie in $[0,1]$.
\eeR

{\bf Why folds and cusps matter for modelling.}
Past a fold the state can jump to another attractor, and reversing the parameter need not restore it. This is the mechanism
of catastrophic regime shifts \cite{Scheffer2001} and tipping cascades \cite{kroenke2020}, and of the backward bifurcation of
epidemic and social-contagion models, where an established state persists although invasion is impossible
\cite{IPB,Kook2021,ArinoMc}. A cusp, where two fold curves meet, marks the change from a gradual and reversible response to an
abrupt and hysteretic one when its outer branches are attracting, which is why cusps serve as design and diagnostic tools
\cite{Guc}. Several cusps give  several regions of bistability and several switching thresholds.

\section{Siphons of a production model.}\label{s:KB}
\beD[Siphons separate absent variables from resident variables on an invariant face]
[][Siphons]
Let $\dot x=f(x)$ be $C^1$ near a forward-invariant region
$\Omega\subseteq\mathbb R^{n_S}_+$, with $f_i(x)\ge0$ on $x_i=0$.
A nonempty index set $\Sigma$ is a \emph{siphon} if
\[
 \mathcal F_\Sigma=\{x\in\Omega:x_i=0\ (i\in\Sigma)\},
 \qquad f_\Sigma(0,y)=0
 \quad\text{for every }(0,y)\in\Omega.
\]
Equivalently, $\mathcal F_\Sigma$ is forward invariant. A minimal
siphon contains no smaller nonempty siphon. With $z=x_\Sigma$ and
$y=x_{\Sigma^c}$, the equations on the face are
$\dot y=f_{\Sigma^c}(0,y)$, called the resident system.
Minimal siphons may overlap and need not form a partition
(Example~\ref{e:KBoverlap}).
\eeD

{\bf The relevance of siphon calculus.}
Several of the methods we use below (for multistability, cusp constructions, propagation of waves,  etc), have already been used  in particular examples of epidemic, genetic,
 neural and threshold models. We are revisiting them here from the point of view  of siphon calculus, which was developed recently in papers like \cite{AAHK,AAHRah,AAHrel,AA26}, building upon the \NGM\ method of \cite{Van},  the persistence result of Angeli, De Leenheer
and Sontag \cite{AdLS},  the algebraic analysis
of siphons and boundary steady states by Shiu and Sturmfels
\cite{ShiuStu}, and on the underlying \CRN\ quest of identifying ``structures" (or the minimal assumptions), which make  the resolution of similar  particular examples work out.  The transverse Jacobian block of any siphon face is Metzler, and hence the Freidlin--G\"artner--Weinberger--Lewis--Li wave-speed formula gives a general principle: it yields the linear spreading speed of the invasion of any resident on a siphon face under Laplacian diffusion\rev{, provided the resident is spatially constant, the transverse block is irreducible and at least one absent species diffuses. It is the speed of the linearization; whether a nonlinear front travels at it is a separate question, left open here}.

Minimal siphons turn out to correspond often to epidemiologic strains \cite{AABH25}. They also distinguish classic models where minimal siphons are typically disjoint from less studied ones, where intersection is possible, like  the
three-population  network of Theorem~\ref{t:KBoverlapCusp},
$\dot x_1=-x_1+\Phi_a(x_1+\alpha x_2)$, $\dot x_2=-x_2+\Phi_a(x_2+\beta x_1x_3)$,
$\dot x_3=-x_3+\Phi_a(x_3+\gamma x_2)$, with $\Phi_a$ a quartic odds response.%

\beT[Autonomous production models admit explicit criteria for siphons and formulas for their transverse Jacobians]
[siphons and transverse invasion][InvF]\label{p:KBsiph}
\label{p:KBresponsezeros}\label{p:KBorder}

Let $\Omega$ be $\mathbb R_+^{n_S}$ or $[0,\bar x]$ with $\bar x\gg0$, and assume
the autonomous production model of Definition~\ref{d:KB} has $f$ of class $C^1$ on a neighbourhood of $\Omega$ and satisfies
\[
 f_i(x)\ge0\quad(x_i=0,\ x\in\Omega).
\]
In the box case assume also $f_i(x)\le0$ on $x_i=\bar x_i$.
The lower-face inequalities hold automatically when $\Pi\ge0$.
For $\varnothing\ne \Sigma\subseteq\{1,\ldots,n_S\}$ write $z=x_\Sigma$, $y=x_{\Sigma^c}$.
Siphons are understood as in Definition~\ref{d:Siphons} (for the lifetime of solutions in the orthant case).
Write $V_\Sigma=V_{\Sigma\Sigma}$, let $\rho$ denote spectral radius, and write $s(M)=\max_{\lambda\in\sigma(M)}\operatorname{Re}\lambda$.
\begin{enumerate}
\item \emph{Face invariance and transverse linearization.}
The set $\Sigma$ is a siphon exactly when
\begin{equation}\label{eq:KB-face-general}
 \Pi_\Sigma(0,y)=V_{\Sigma\Sigma^c}\,y\quad\text{for every }y\text{ with }(0,y)\in\Omega.
\end{equation}
Since $\Pi_\Sigma\ge0$ and $V_{\Sigma\Sigma^c}\le0$, this is equivalent to
\begin{equation}\label{eq:KB-face}
 V_{\Sigma\Sigma^c}=0,\qquad \Pi_\Sigma(0,y)=0\quad\text{for every }y\text{ with }(0,y)\in\Omega.
\end{equation}
On the face, with $J=Df$,
\[
 J(0,y)=\begin{pmatrix}J^\perp&0\\J_{\mathrm{mix}}&J_{\mathrm{tan}}\end{pmatrix},
 \quad F_\Sigma(y)=\partial_z\Pi_\Sigma(0,y),\quad
 J^\perp=F_\Sigma(y)-V_\Sigma,
\]
where $J_{\mathrm{mix}}=\partial_zf_{\Sigma^c}(0,y)$ and
$J_{\mathrm{tan}}=\partial_yf_{\Sigma^c}(0,y)$.
The matrix $J^\perp$ is Metzler, $F_\Sigma(y)\ge0$, and the next-generation matrix and invasion number on the face are $K_\Sigma(y)=F_\Sigma(y)V_\Sigma^{-1}$ and $R_\Sigma(y)=\rho(K_\Sigma(y))$. Then
\[
 \operatorname{sign}s(J^\perp)
 =\operatorname{sign}\bigl(R_\Sigma(y)-1\bigr).
\]
At a
boundary equilibrium every eigenvalue of the full Jacobian has negative
real part exactly when this holds for both $J^\perp$ and
$J_{\mathrm{tan}}$. For $\Sigma^c=\varnothing$ the tangential condition is omitted.

\item \emph{Vanishing of the response on a face.} For \eqref{eq:KB-general-response}, suppose
$b_i(0)>0$ for every $i\in \Sigma$ and put
\[
 \mathcal Q_{i,\Sigma}=\{h_i+P_i(0,y):(0,y)\in\Omega\}.
\]
Then the face criterion becomes
\begin{equation}\label{eq:KB-response-zero-test}
 V_{\Sigma\Sigma^c}=0,\qquad
 \mathcal Q_{i,\Sigma}\subseteq\{q:\Phi_i(q)=0\}\quad(i\in \Sigma).
\end{equation}
Keeping $V$, $b_i$ and $P_i$ fixed, replacing responses by
nonnegative functions with the same zero sets
on attained arguments preserves every siphon. Strictly positive responses
on those arguments exclude nonempty siphons when all $b_i(0)>0$.
For an index with $b_i(0)=0$, its production vanishes on the face
regardless of its response.

\item \emph{Polynomials with nonnegative coefficients.} In \eqref{polR}, write $P_i(x)=\sum_{\nu\ne0}a_{i\nu}x^\nu$, with multi-indices $\nu\in\mathbb Z_{\ge0}^{n_S}$ and $a_{i\nu}\ge0$, and
assume $\Phi_i(q)=0$ for $q\le0$ and $\Phi_i(q)>0$ for $q>0$.
On the full orthant, the criterion is
\[
 V_{\Sigma\Sigma^c}=0,\quad h_i\le0\ (i\in \Sigma),\quad
 \operatorname{supp}\nu\cap \Sigma\ne\varnothing
 \quad(i\in \Sigma,\ a_{i\nu}>0).
\]
On a finite box it is equivalent when $h_i=0$ for $i\in \Sigma$.
For arbitrary $h_i$, the exact box criterion is
\[
 V_{\Sigma\Sigma^c}=0,\qquad h_i+P_i(0,\bar x_{\Sigma^c})\le0\quad(i\in \Sigma).
\]
Here $\operatorname{supp}\nu=\{j:\nu_j>0\}$. For $P_i(x)=\beta(Wx)_i$, $W\ge0$, zero stimulus,
and $\beta>0$, the monomial condition is $W_{\Sigma\Sigma^c}=0$.
With diagonal $V$, minimal siphons are the source strongly connected
components of the influence graph $j\to i$ when $w_{ij}>0$---the closed
recurrent classes in the opposite row-transition convention for stochastic $W$.

\item \emph{Order in absent variables.}
On a siphon face write
\[
 \Pi_\Sigma(z,y)=\Pi_\Sigma(0,y)+F_\Sigma(y)z+\mathcal R_\Sigma(z,y).
\]
If $\mathcal R_\Sigma=O(\|z\|^2)$ at each resident state under consideration,
changing only the remainder preserves $J^\perp$.
More generally, adding a $C^1$ vector field $H$ preserves face dynamics
and the transverse matrix under
\[
 H(0,y)=0,\qquad \partial_zH_\Sigma(0,y)=0.
\]
Preservation of the entire face Jacobian additionally requires
$\partial_zH_{\Sigma^c}(0,y)=0$; otherwise $J_{\mathrm{mix}}$ may change.
These identities must hold along the face, not just at one equilibrium.
In particular, $\Pi_\Sigma(z,y)-\Pi_\Sigma(0,y)=O(\|z\|^2)$ gives
$J^\perp=-V_\Sigma$ along every resident trajectory.
For nonnegative production on a siphon face, $\Pi_\Sigma(0,y)=0$.
\end{enumerate}
\eeT
\begin{proof}
\begin{enumerate}
\item On the face the absent equations are
$\dot z=-V_{\Sigma\Sigma^c}y+\Pi_\Sigma(0,y)$.
Invariance forces this expression to vanish at every initial point
of the face. Conversely, under \eqref{eq:KB-face-general}, a solution
of the resident equations with $z\equiv0$ solves the full system;
uniqueness proves invariance.

Since the off-diagonal entries of $V$ are nonpositive, both
$-V_{\Sigma\Sigma^c}y$ and $\Pi_\Sigma(0,y)$ are nonnegative when production is
nonnegative on the face. Their sum vanishes only if both vanish.
Taking $y=t e_j$ for each resident index and an admissible $t>0$
gives $V_{\Sigma\Sigma^c}=0$. This uses the specified orthant or nondegenerate box.

The triangular form and the Metzler property of $J^\perp$
are Proposition~\ref{p:Tri}. Its one-sided difference quotient
argument, applied to $\Pi_i$ and including $j=i$, proves $F_\Sigma\ge0$ when
production is nonnegative on $\Omega$ and vanishes on the face.
Every principal submatrix $V_\Sigma$ is a nonsingular $M$-matrix;
the nonnegative-production splitting $F_\Sigma-V_\Sigma$ gives the
spectral-radius sign test.
The stability assertion is Proposition~\ref{p:Tri}(2).

\item On the face,
$\Pi_i(0,y)=b_i(0)\Phi_i(h_i+P_i(0,y))$.
For $b_i(0)>0$ its vanishing is precisely the response zero-set
condition; for $b_i(0)=0$ it is automatic. Together with
\eqref{eq:KB-face}, this proves all zero-set assertions.

\item At $y=0$ the argument of $\Phi_i$ equals $h_i$, so $h_i\le0$ is necessary.
Any positive-coefficient monomial containing only resident
variables grows without bound along $y=t\mathbf1$ and eventually
makes this argument positive. Hence such monomials are excluded on
the full orthant. Conversely, every monomial containing an absent
variable vanishes on the face.
On a box the nonnegative polynomial attains its largest resident
value at $\bar x_{\Sigma^c}$, proving the exact box test; for $h_i=0$ every
resident-only monomial is strictly positive there.
For linear $P_i$ and diagonal $V$, the criterion says that no
edge enters $\Sigma$ from outside. Such a set is a union of strongly
connected components and contains a source component of the full
graph. Each source component is itself predecessor closed, so
these and only these are the minimal nonempty siphons.
With nondiagonal $V$, the same proof includes the transfer edges.

\item A differentiable $O(\|z\|^2)$ remainder has zero derivative
at $z=0$. Also $H(0,y)=0$ implies $\partial_yH(0,y)=0$.
These facts identify the preserved entries of the face Jacobian.
\end{enumerate}
\end{proof}

\ssec{Re-splitting the family of the introduction.}
For the family \eqref{eq:intro-family}, write the
bracket as $B_i(x)$ and work on an invariant box $\Omega$ with
$\rho_ix_i\le1$. The production $\Pi_i=(1-\rho_ix_i)B_i\ge0$, but the
nonlinear removal is outside the linear form $-Vx$. Choose
\[
 \mu_i\ge\rho_i\sup_\Omega B_i+\sup_{[0,\bar x_i]}L_i',\qquad
 V'=V+\operatorname{Diag}(\mu_i),\quad \Pi' =\Pi-L+\operatorname{Diag}(\mu_i)x.
\]
Then $\partial_i\Pi_i'\ge0$ and $\Pi_i'(0,x_{-i})=B_i(0,x_{-i})\ge0$,
so $\Pi'\ge0$. The ODE and its upper-face condition are unchanged;
$V'$ remains a nonsingular $M$-matrix. Off-diagonal derivatives are
unchanged, so the Jacobian is Metzler on the box when the responses increase
in the other variables. Capacity factors can make $\partial_i\Pi_i<0$.

For SIS, with $\Pi_i=(1-x_i)H_i$, the
choice $\mu_i\ge\sup H_i$ restores monotonicity on the cube without changing the siphons or the
transverse Jacobians; the value of $R_\Sigma$ depends on the splitting, while its position relative to
one does not. What the siphon calculus uses holds in any case: on a siphon face,
$\Pi_\Sigma(0,y)=0$ and $\Pi\ge0$ give $\partial_z\Pi_\Sigma(0,y)\ge0$ entrywise, the diagonal included.

\beP[Positive networks have invariant boxes and cooperative networks have extremal equilibria]
[existence, bounds, and extremal equilibria][KBpositive]\label{p:KBwell}
Consider an autonomous production model of Definition~\ref{d:KB} with locally Lipschitz production. Let $\Omega=[0,\bar x]$, $\bar x\gg0$, and suppose $\Pi\ge0$ on $\Omega$ and
$f_i\le0$ on $x_i=\bar x_i$. Then each initial state in $\Omega$ has a unique
global solution remaining in $\Omega$, and an equilibrium exists there.
More generally, if $\Pi$ is defined on the full orthant and $0\le\Pi\le b$,
then without monotonicity all nonnegative solutions are global and
\[
 x(t)\le\bar x+e^{-Vt}(x(0)-\bar x),\qquad
 \bar x=V^{-1}b,\qquad\limsup_{t\to\infty}x(t)\le\bar x.
\]
\begin{enumerate}
\item[(a)] \emph{Cooperative box.} If in addition $J=-V+\partial_x\Pi$ is
Metzler on $\Omega$, the trajectories from $0$ and $\bar x$
converge to the least and greatest equilibria $x^-\le x^+$; every
equilibrium lies between them and every trajectory satisfies
\[
 x^-\le\liminf_{t\to\infty}x(t)\le\limsup_{t\to\infty}x(t)\le x^+.
\]
Consequently a unique equilibrium attracts the whole box, without
requiring contraction.
\item[(b)] \emph{Globally bounded monotone production.} If $0\le\Pi\le b$,
$b\gg0$, and $\partial_x\Pi\ge0$ entrywise (diagonal included) on the
orthant, then each $[0,a\bar x]$, $a\ge1$, is invariant and, for $a>1$,
absorbs every nonnegative solution in finite time. All equilibria lie in
$[0,\bar x]$; the least and greatest are also the limits of
\[
 x^-_0=0,\quad x^+_0=\bar x,\qquad x^\pm_{n+1}=V^{-1}\Pi(x^\pm_n).
\]
\end{enumerate}
\eeP
\begin{proof}
At $x_i=0$,
$f_i=-\sum_{j\ne i}v_{ij}x_j+\Pi_i\ge0$; at $x_i=\bar x_i$ the
assumed inequality is reversed. These inward signs give invariance.
Local uniqueness and compactness then give global solutions.
For existence of an equilibrium, apply Brouwer to the continuous map
$x\mapsto\operatorname{proj}_\Omega(x+\varepsilon f(x))$.
At a fixed point, an interior coordinate forces $f_i=0$;
a zero coordinate forces $f_i\le0$, and an upper coordinate forces
$f_i\ge0$. The inward signs force equality in both boundary cases.

On the orthant, positivity of $e^{-Vt}$ and variation of constants give
\[
 x(t)=e^{-Vt}x(0)+\int_0^t e^{-V(t-s)}\Pi(x(s))\,ds
 \le e^{-Vt}x(0)+V^{-1}(I-e^{-Vt})b.
\]
This bound excludes finite-time blow-up; $e^{-Vt}\to0$ gives
the asserted componentwise upper limit.

In (a), the Metzler Jacobian makes the flow $\varphi_t$ order preserving
\cite{smith1995}, \cite{Hirsch}. Invariance and the semigroup identity give
\[
 \varphi_t(0)\le\varphi_{t+h}(0),\qquad
 \varphi_{t+h}(\bar x)\le\varphi_t(\bar x)\qquad(t,h\ge0).
\]
Both bounded trajectories therefore have limits $x^-,x^+$.
For each $h\ge0$, continuity gives
$\varphi_h(x^\pm)=\lim_{t\to\infty}\varphi_{t+h}(0\text{ or }\bar x)=x^\pm$,
so both limits are equilibria.
Comparison with $\varphi_t(0)$ and $\varphi_t(\bar x)$ bounds every trajectory
and every equilibrium, proving extremality and the limit inequalities.

In (b), $V^{-1}\ge0$ has no zero row, so $\bar x\gg0$.
Since $f(a\bar x)\le-(a-1)b$, cooperativity implies the upper-face
inequalities on $[0,a\bar x]$: if $x_i=a\bar x_i$ and $x\le a\bar x$,
then $f_i(x)\le f_i(a\bar x)\le0$.
For $a>1$, the preceding global bound and decay of $e^{-Vt}$ give
$x(t)\le a\bar x$ for all sufficiently large $t$.
Finally $T=V^{-1}\Pi$ is continuous and increasing, with
$T(0)\ge0$ and $T(\bar x)\le\bar x$. Thus
\[
 0=x_0^-\le x_1^-\le\cdots\le x_n^-\le x_n^+\le\cdots
 \le x_1^+\le x_0^+=\bar x .
\]
Both sequences converge to fixed points of $T$. Every equilibrium
$e$ satisfies $e=T(e)\le\bar x$ and, by induction,
$x_n^-\le e\le x_n^+$; hence the limits are precisely the least
and greatest equilibria.
\end{proof}

\beP[Every positive smooth box model admits a  production representation]
[representation of positive box models][KBrepresentation]
Let $f$ be $C^1$ on a neighbourhood of an invariant box $[0,\bar x]$,
$\bar x\gg0$. Choose
\[
 \mu_i>\max\{0,\sup_{x\in[0,\bar x]}(-\partial_i f_i(x))\},\qquad
 V=\operatorname{Diag}(\mu_i),\qquad \Pi_i(x)=f_i(x)+\mu_i x_i.
\]
Then $\Pi\ge0$ on the box, so the same ODE has a representation with
nonnegative production and inward upper faces. If $f$ is cooperative there, $\partial_x\Pi\ge0$ entrywise,
including the diagonal.
The chosen diagonal coefficients describe a mathematical splitting;
they need not equal the physical removal coefficients of a given model.
\eeP
\begin{proof}
Invariance gives $f_i(0,x_{-i})\ge0$. Integration in the $i$th
coordinate yields
\[
 \Pi_i(x)=f_i(0,x_{-i})+
 \int_0^{x_i}[\partial_i f_i(s,x_{-i})+\mu_i]\,ds\ge0.
\]
For $j\ne i$, $\partial_j\Pi_i=\partial_j f_i$, while
$\partial_i\Pi_i=\partial_i f_i+\mu_i>0$, proving the derivative assertion. The original inward boundary conditions are
unchanged because the vector field itself is unchanged.
\end{proof}
\beR[Smooth ecological and regulatory models]
[ecological models][KBecoG]
Thus the results that assume only nonnegative production on an invariant
box also apply to smooth ecological or regulatory models, even when their natural
production and loss terms have a different form. Global bounds on the
whole orthant and the special threshold bifurcation formulas remain
additional hypotheses.
\eeR

\beR[Smooth onset of a response]
[smooth onset][KBonset]
If $\Phi$ is $C^1$ on all of $\mathbb R$ and vanishes on
$(-\infty,0]$, necessarily $\Phi'(0)=0$. A positive linear onset
$\Phi'(0+)>0$ instead uses a response defined and $C^1$ from the right
on $[0,\infty)$, with a smooth extension for local calculations;
its zero extension to negative arguments is not $C^1$. The results
use smoothness of $\Pi$ on the admissible phase space, in the sense of
Definition~\ref{d:KB}, not an incompatible smoothness assumption on that
zero extension. For the odds responses of Definition~\ref{d:Odd} the
formula $\phi/(1+\phi)$ itself is a smooth extension to a neighbourhood
of $[0,\infty)$, since $1+\phi>0$ near zero.
\eeR

\paragraph{The results using each subclass.}
\begin{enumerate}

\item \emph{Positive Lur'e systems} are used in Theorems~\ref{t:KBsharp} and~\ref{t:SpdBd} (common autonomous $\Phi$, zero stimulus and a slope bound). Their linear production subclass also admits the contraction bound \eqref{eq:KB-shared-contraction} with $0\le\Phi_i'\le L_i$.

\item \emph{Polynomial $P_i$} are used in Theorem~\ref{t:InvF}
(criterion with $\Phi_i=0$ on $(-\infty,0]$, $\Phi_i>0$ on $(0,\infty)$),
Theorem~\ref{t:siphMulti} ($h_i=0$, common $\Phi$ in the class $\mathcal C_0$).

\item \emph{State-dependent factors} \eqref{eq:KB-general-response} are used in
Theorem~\ref{t:InvF} (zero sets of $\Phi_i$ and siphons, with
$b_i(0)>0$) and Proposition~\ref{p:KBshared} (shared responses; $b_i$ constant in the contraction
test).
With
$b_i(x_i)=1-x_i$, $\Phi_i(q)=q$, $h_i=0$ and $P_i=H_i$ polynomial it is the SIS
equation.
\end{enumerate}

\paragraph{Transverse dynamics and higher-order activation.}

\beR[The transverse matrix of the capacity-factor form and of SIS]
[transverse matrix][KBtransG]
The initial siphon decomposition gives the transverse Metzler property
without full cooperativity; Theorem~\ref{t:InvF} supplies its
production splitting and spectral test. For the
capacity-factor form \eqref{eq:KB-capacity} the factor $1-\rho_ix_i$ equals one
on the face and drops out of the transverse matrix; for SIS contagion,
\[
 J^\perp(y)=A_{\Sigma\Sigma}+\partial_{x_\Sigma}g_\Sigma(0,y)-V_\Sigma.
\]
Only at the origin, or when $\partial_{x_\Sigma}g_\Sigma(0,y)=0$, does this
reduce to $A_{\Sigma\Sigma}-V_\Sigma$. For instance, $g_i=x_jx_k$ with
$i,j\in \Sigma$ and $k\notin \Sigma$ contributes the resident value $y_k$
to the $(i,j)$ transverse entry. Total polynomial degree and degree
in the absent variables must therefore be distinguished.
\eeR

\section{Autonomous threshold networks admit arbitrarily many stable equilibria}
The following two theorems construct arbitrary multistability under
different restrictions on the response.
Theorem~\ref{t:siphMulti} fixes the response zero set and
allows its slope to vary, preserving every siphon and every transverse
Jacobian while creating arbitrarily many stable positive equilibria.
Indeed, for its production--loss splitting, on every siphon face,
\[
 K_\Sigma(y)
 =\partial_{x_\Sigma}(\Pi_\Phi)_\Sigma(0,y)
   (V_\Sigma)^{-1}=0,\qquad R_\Sigma(y)=0,
\]
because $\Phi'(0)=0$.
Thus unchanged, strictly subcritical boundary invasion numbers can
coexist with arbitrarily large interior multistability.
Theorem~\ref{t:KBsharp} instead fixes a slope bound and identifies a sharp
extension $\bar R_0$ of the next-generation number: below $1/L$, every response gives global exponential
convergence; above $1/L$, suitable responses give arbitrarily many stable
equilibria.
Both proofs prescribe response values at selected arguments and make
the response derivative vanish there, giving $J=-V$ at the
constructed equilibria; neither theorem requires the other.

\paragraph{Arbitrary multistability with unchanged siphons.}
Without a uniform derivative bound, polynomials $P_i$ give an exact
support criterion for constructing stable positive equilibria while
preserving every siphon and its transverse Jacobian. Here production reaches every coordinate if every index is reachable from $I=\{i:P_i\not\equiv0\}$ through the directed edges $j\to i$ with $i\ne j$ and $-v_{ij}>0$, allowing paths of length zero. For diagonal $V$ this means that every $P_i$ is nonzero.

\beT[Polynomial production threshold networks with a common flat-onset CDF
admit arbitrarily many stable positive equilibria exactly when production
reaches every coordinate]
[][siphMulti]
Fix a nonsingular $M$-matrix $V$, $b\gg0$, and production polynomials
$P_i(x)=\sum_{\nu\ne0}a_{i\nu}x^\nu$ with $a_{i\nu}\ge0$,
allowing $P_i\equiv0$. Consider the production threshold model \eqref{eq:KB}, \eqref{polR}
with $h_i=0$ and a common response $\Phi$, and write $\Pi_\Phi$ for its production,
where $\mathcal C_0$ consists of the $C^\infty$ CDFs satisfying
$\Phi(q)=0$ for $q\le0$, $\Phi(q)>0$ for $q>0$, and $\Phi'\ge0$.
Such a $\Phi$ vanishes identically on $(-\infty,0]$, so all its derivatives
vanish at $0$; this is used in items 2, 3 and 4 below.
Putting
\be{tB}
 I=\{i:P_i\not\equiv0\},\qquad
 \tilde b=V^{-1}\bigl(b_i\mathbf1_{\{i\in I\}}\bigr)_i,
\ee
the following hold:
\begin{enumerate}
\item A positive equilibrium for any $\Phi\in\mathcal C_0$ requires
 $\tilde b\gg0$.
\item If $\tilde b\gg0$, every $m\ge1$ admits a common
 $\Phi_0\in\mathcal C_0$ with strictly ordered positive equilibria
 \[
 x^{(k)}=t_k\tilde b,\quad 0<t_1<\cdots<t_m<1,\qquad J(x^{(k)})=-V.
 \]
 These and the origin are locally exponentially stable. A small
 common-response perturbation makes $\Phi'(q)>0$ for all $q>0$ while
 retaining $m$ strictly ordered stable positive equilibria.
\item\rev{At the equilibria of item 2 the response satisfies $\Phi'=0$, so
 $J(x^{(k)})=-V$. Since $\Phi'(0)=0$, the transverse production derivative
 vanishes at every point of a siphon face, $J_\Sigma^\perp=-V_\Sigma$, and
 the associated next-generation matrix is zero, so every invasion number of
 every siphon face vanishes there. Arbitrarily large interior multistability
 is thus compatible with a strictly subcritical boundary, and no link
 between the two is asserted. The perturbation with $\Phi'>0$ keeps the
 count, at Jacobians that are Metzler and Hurwitz without being $-V$. What
 the counts become for a sigmoid fixed in advance is open.}
\item Every allowed response has the same siphons, namely
 \[
 \Sigma\ne\varnothing:\quad V_{\Sigma\Sigma^c}=0,
 \qquad \operatorname{supp}\nu\cap\Sigma\ne\varnothing
 \quad(i\in\Sigma,\ a_{i\nu}>0).
 \]
 On each face $J_\Sigma^\perp(y)=-V_\Sigma$ for every
 admissible resident state; in particular $J(0)=-V$.
\end{enumerate}
Thus $\tilde b\gg0$ is necessary and sufficient for arbitrarily many stable positive equilibria as the common response varies
in this class; additional equilibria are not excluded.
\rev{In this construction $V$, the production polynomials $P_i$, and their coefficients are fixed; the common response is chosen for each prescribed count. This differs from results which fix an activation function or vary the connection matrix \cite{WangLuChen2010,Betteti2025}.}
\eeT
\begin{proof}
Throughout, we work with the polynomial production threshold class of the statement, with zero stimuli and a common response in $\mathcal C_0$. The common box $[0,V^{-1}b]$ is invariant, since
$0\le\Pi_\Phi\le V\tilde b\le b$ and $\partial_x\Pi_\Phi\ge0$:
if $\bar x=V^{-1}b$ and $x_i=\bar x_i$, then
$f_i(x)\le-(V\bar x)_i+b_i=0$, and the lower-face derivatives are
nonnegative.
\begin{enumerate}
\item Every equilibrium satisfies $x=V^{-1}\Pi_\Phi(x)\le\tilde b$, because
 $\Pi_\Phi$ is bounded by $b_i$ in the coordinates of $I$ and vanishes in
 the others. A positive equilibrium therefore forces $\tilde b\gg0$.

\item Let $\tilde b\gg0$. Each $q_i(t)=P_i(t\tilde b)$, $i\in I$, increases strictly from
 zero. Choose $t_m<1$ and recursively choose earlier $t_k$ so that
 \[
 a_k:=\min_{i\in I}q_i(t_k)\le \widehat b_k:=\max_{i\in I}q_i(t_k)<a_{k+1}.
 \]
 Construct a smooth density, zero on $(-\infty,0]$ and every
 $[a_k,\widehat b_k]$, with respective masses $t_1$, $t_{k+1}-t_k$, and
 $1-t_m$ on $(0,a_1)$, $(\widehat b_k,a_{k+1})$, and
 $(\widehat b_m,\infty)$. Positive pieces flat at their endpoints are
 normalized multiples of
 \[
 e^{-1/((q-a)(b-q))}\quad(a<q<b),\qquad
 e^{-1/(q-\widehat b_m)-(q-\widehat b_m)}\quad(q>\widehat b_m).
 \]
 Their cumulative integral $\Phi_0$ belongs to $\mathcal C_0$, equals
 $t_k$ on $[a_k,\widehat b_k]$, and has zero slope there. Hence
 $\Pi_{\Phi_0}(t_k\tilde b)=t_kV\tilde b$ and $J(t_k\tilde b)=-V$. All derivatives of a
 response in $\mathcal C_0$ vanish at zero, so the origin carries the same
 stable Jacobian. Perturb by
 $\Phi_\varepsilon=(1-\varepsilon)\Phi_0+\varepsilon G$ with
 $G(q)=e^{-1/q}$ for $q>0$ and zero otherwise. The implicit-function
 theorem preserves the selected positive equilibria, their strict ordering
 and their stability, while $\Phi_\varepsilon'>0$ on $(0,\infty)$.

\item\rev{The equilibria of item 2 lie in the flat intervals
 $[a_k,\widehat b_k]$, where $\Phi_0'=0$, whence $J(x^{(k)})=-V$. On a
 siphon face $P_i(0,y)=0$ for $i\in\Sigma$ by item 4 below, and every derivative
 of a response in $\mathcal C_0$ vanishes at $0$; the transverse production
 derivative therefore vanishes, $J_\Sigma^\perp=-V_\Sigma$, and the
 associated next-generation matrix is zero, so each invasion number is zero.
 The two halves are established separately, at one and the same response, so
 no implication between them is claimed. The implicit-function argument of
 item 2 carries the selected equilibria through the perturbation, at
 Jacobians that are Metzler and Hurwitz without equalling $-V$. The counts
 for a sigmoid fixed in advance are left open.}

\item On a zero face,
 $f_\Sigma(0,y)=-V_{\Sigma\Sigma^c}y+
 (b_i\Phi(P_i(0,y)))_{i\in\Sigma}$ is a sum of nonnegative terms. It
 vanishes identically exactly under the displayed transfer and monomial
 conditions, which is the asserted characterization. There
 $P_i(0,y)=0$ for $i\in\Sigma$, so $\Phi'(0)=0$ gives the stated
 transverse matrix for every response in $\mathcal C_0$, the perturbation
 included, and $J(0)=-V$ is the case $\Sigma=\{1,\dots,n\}$.
\end{enumerate}
\end{proof}

\paragraph{Support and scope.}
A common response meets any response-sharing condition within the groups of Section~\ref{s:KBred}
or overlapping siphons. Resident dynamics may change although the
siphons and transverse matrices do not.

\paragraph{Arbitrary multistability under a slope bound.}
\beT[Autonomous positive Lur'e systems with a common slope-bounded
response $L$ have a sharp threshold at $\bar R_0=1/L$, with one globally
attracting equilibrium below it and arbitrary multistability above it]
[sharp threshold $\bar R_0=1/L$][KBsharp]\label{p:KBinterpolate}\label{o:KBsharp}\label{p:KBmaxstable}\label{o:KBmaxstable}
Fix a nonsingular $M$-matrix $V$, nonnegative matrices $B$ and $W$, and
$L>0$.
Let $\mathcal C_L$ consist of the functions $\Phi$, smooth on $[0,\infty)$, with $\Phi\ge0$ and
$0\le \Phi'\le L$ there; the arguments $\beta Wx$ are nonnegative on the orthant. Distribution functions
with this slope bound belong to $\mathcal C_L$, and so does the identity when $L\ge1$.
Consider
\begin{equation}\label{eq:KB-sharp-family}
 \dot x=-Vx+B\Phi(\beta Wx),\qquad x\in\mathbb R_+^{n_S},
 \quad \Phi\in\mathcal C_L,\quad\beta\ge0,
\end{equation}
with the same scalar response applied to every coordinate of $\beta Wx$.
Here $B\in\mathbb R_+^{n_S\times n_R}$ and $W\in\mathbb R_+^{n_R\times n_S}$ for some $n_R\ge1$. Put
\[
 r_0=\rho(BWV^{-1})=\rho(V^{-1}BW),\qquad
 \bar R_0=\rho(\beta BWV^{-1})=\beta r_0.
\]

\begin{enumerate}
\item If $\bar R_0<1/L$, every $\Phi\in\mathcal C_L$ gives a unique
globally exponentially attracting equilibrium.
\item If $\bar R_0>1/L$, for every $m\ge2$ there exists
$\Phi\in\mathcal C_L$ giving at least $m$ distinct locally exponentially
stable equilibria. One can choose $\Phi(q)=0$ for $q\le0$ and obtain
\[
 x^{(k)}=t_kr,\quad 0=t_0<t_1<\cdots<t_{m-1},
 \qquad r\ge0,\quad r\ne0,\qquad J(x^{(k)})=-V.
\]
Alternatively, if $V^{-1}B\mathbf1\gg0$, one can choose $0<\Phi(q)<1$ and
$\Phi'(q)>0$ on $\mathbb R$ and obtain $m$ strictly positive stable
equilibria.
The response may depend on $m$, and in both cases it is a distribution function, so
the construction stays within bounded responses.
\end{enumerate}
\rev{Item 1 is the Lur'e contraction theorem of \cite{Davydov} for this class, and is quoted here; what is new is item 2, the construction of a response in $\mathcal C_L$ for each finite $m$. Accordingly $\bar R_0=1/L$ is the exact transition, within $\mathcal C_L$, between the two alternatives above: below it every response of the class has a unique globally attracting equilibrium, and above it some response of the class has at least $m$ locally exponentially stable equilibria. By \emph{coexistence of equilibria} is meant here the second alternative, namely the existence of two or more locally exponentially stable equilibria for one response of the class; the response is allowed to depend on $m$, so this is not a statement about a response fixed in advance. Multiple stable equilibria of neural networks with a fixed piecewise-linear or sigmoid activation, the synaptic matrix being the free datum, are counted in \cite{WangLuChen2010,Betteti2025}; here $V$ and $W$ are fixed and the response is constructed.}
No irreducibility is required; uniqueness at equality is not asserted.
A common additive stimulus can be absorbed into $\Phi$, provided the shifted response belongs to $\mathcal C_L$.

Thus the stable count is unbounded over this response class, even
for $n_S=1$. With $\Phi>0$ and $\Phi'>0$, the interaction graph remains
fixed; if $V^{-1}B\mathbf1\gg0$, there are no nonempty siphons. In the zero-onset construction, the selected equilibria are ordered
componentwise and lie on the ray $\{tr:t\ge0\}$.
This ray need not be invariant. The full-support perturbation
preserves their number and stability, but need not preserve
their location on that ray.
\eeT

\beR[$\bar R_0$ extends the next-generation number][]
[NGMg]
Write
\[
 J(x)=F(x)-V,\qquad
 F(x)=\beta B\operatorname{Diag}(\Phi'(\beta Wx))W.
\]
Since $0\le F(x)\le\beta LBW$,
\[
 R(x):=\rho(F(x)V^{-1})\le L\bar R_0.
\]
At an equilibrium $x_*$, the NGM matrix criterion gives
$J(x_*)$ Hurwitz if and only if $R(x_*)<1$. If $\Phi(0)=0$, the origin is an equilibrium and
$R(0)=\rho(\beta\Phi'(0)BWV^{-1})=\Phi'(0)\bar R_0$ is the classical $R_0$, so the classical test is $\Phi'(0)\bar R_0<1$ (equivalently $\bar R_0<1/\Phi'(0)$ if $\Phi'(0)>0$), which is $\bar R_0<1/L$ when $L=\Phi'(0)$.
In general the test $\bar R_0<1/L$ is uniform over all states and
responses in $\mathcal C_L$, and the bar marks that $\bar R_0$ is not the reproduction number
of a specified equilibrium. The condition $\bar R_0<1/L$ ensures global
exponential convergence; $\bar R_0>1/L$ permits arbitrarily many stable
equilibria through a suitable choice of response.
\eeR

\begin{proof}
\emph{Interpolation used in the construction.}
Let $0<c<L$ and $0<q_1<\cdots<q_d<1/c$. There is a $C^\infty$ CDF
$\Phi_0:\mathbb R\to[0,1]$ satisfying
\[
 \Phi_0(q)=0\ (q\le0),\qquad 0\le \Phi_0'(q)<L,
 \qquad \Phi_0(q_j)=cq_j,\quad \Phi_0'(q_j)=0.
\]
Its density has compact support in $(0,\infty)$ and vanishes near each
$q_j$, including $\Phi_0'(0)=0$.

Choose $\eta\in C_c^\infty(0,1)$, $\eta\ge0$, with $\int\eta=1$ and
$\|\eta\|_\infty<L/c$. Such densities have supremum arbitrarily close
to one. Put $q_0=0$ and choose $a>(1-cq_d)\|\eta\|_\infty/L$. Set
\[
 \varphi(q)=\begin{cases}
 c\eta((q-q_{j-1})/(q_j-q_{j-1})),&q_{j-1}<q<q_j,\\
 (1-cq_d)\eta((q-q_d)/a)/a,&q_d<q<q_d+a,\\
 0,&\text{otherwise}.
 \end{cases}
\]
The pieces have masses $c(q_j-q_{j-1})$ and $1-cq_d$, vanish near
their endpoints, and have supremum less than $L$. Their cumulative
integral is the required response.

\emph{Contraction below the bound.}

For $\bar R_0<1/L$, the matrix $V-\beta LBW$ is a nonsingular $M$-matrix, because
$\rho(\beta LBWV^{-1})<1$. Put
\[ p=(V-\beta LBW)^{-1}\mathbf1\gg0 .\]
The bound $\partial_x\Pi\le\beta LBW$ gives $J(x)\le-V+\beta LBW$ entrywise, so
$J(x)p\le-\mathbf1$ for every $x$. Since $J(x)$ is Metzler, this gives contraction
in the weighted maximum norm $\max_i|x_i|/p_i$, at rate $1/\max_ip_i$. The vector field is globally
Lipschitz on the orthant, whose faces $x_i=0$ point inward, so nonnegative solutions are global. The
time-$t$ maps are strict contractions of the complete orthant, so they have one common fixed point,
which is the unique equilibrium, and every solution converges to it exponentially.
\emph{Stable equilibria above the bound.}
For $\bar R_0>1/L$, choose a nonnegative Perron vector
\[
 V^{-1}BWr=r_0r,\qquad r_0=\rho(V^{-1}BW),\qquad r\ne0,
 \qquad BWr=r_0Vr\ne0,
\]
and put $c=1/\bar R_0<L$. Choose the positive $t_k$ so that
$\max_i\beta t_{m-1}(Wr)_i<1/c$. Apply the interpolation construction above to all distinct positive arguments
$\beta t_k(Wr)_i$. At these and at zero,
\[
 B\Phi_0(\beta Wt_kr)=c\beta t_kBWr=c\bar R_0t_kVr=t_kVr,
 \qquad
 J(t_kr)=-V+\beta B\operatorname{Diag}(\Phi_0'(\beta Wt_kr))W=-V.
\]
The $m$ equilibria are exponentially stable and lie in
$[0,V^{-1}B\mathbf1]$.

\emph{Strictly positive equilibria and full support.}
Set
\[
 G(q)=(1+e^{-Lq})^{-1},\qquad
 \Phi_\varepsilon=(1-\varepsilon)\Phi_0+\varepsilon G.
\]
Since $0<G'\le L/4$, this remains in $\mathcal C_L$ and is strictly
increasing. The implicit-function theorem continues the equilibria, with
\[
 \left.\frac{d x^{(k)}(\varepsilon)}{d\varepsilon}\right|_0
 =V^{-1}BG(\beta Wx^{(k)})-x^{(k)},
\]
the term in $\Phi_0'$ dropping out because $\Phi_0'=0$ at these arguments.
The first term is strictly positive when $V^{-1}B\mathbf1\gg0$, that is when
$V^{-1}B\ge0$ has no zero row, which for $B=I$ is automatic.
Every initially zero coordinate therefore becomes positive, while
positive coordinates remain positive. Disjoint neighbourhoods and
hyperbolicity preserve the count and stability for one small
$\varepsilon>0$. Replacing $\Phi(q)$ by $\Phi(q+h)$ absorbs a common stimulus
$h\mathbf1$. Full support excludes zero faces; for these constructions $\beta>0$,
and the off-diagonal Jacobian entries are positive exactly when $-v_{ij}>0$ or $(BW)_{ij}>0$.
\end{proof}

\beR[Fixed production coefficients and the scope of the stable-equilibrium count]
[fixed coefficients and counting scope][KBfixedamp]\label{r:KBcountscope}
The amplification $B$ is carried by Theorem~\ref{t:KBsharp} itself, so the
case $B=\operatorname{Diag}(b_i)$, $b_i>0$, needs no separate treatment: the
interpolation values are $\Phi(q)=q/(\beta r_0)$ with
$r_0=\rho(V^{-1}BW)$, and the selected Jacobians are $-V$.

The response in Theorem~\ref{t:KBsharp} varies with the requested count. Finite bounds require
additional restrictions, such as a fixed distribution family, density
oscillation bounds, or fixed polynomial degree. This construction
preserves neither unimodality nor fixed-degree SIS structure, and does
not assert infinitely many hyperbolic equilibria for one response.
\eeR

\beC[Subcritical threshold networks entrain to periodic stimuli][periodic entrainment][Ent]
Let $V$ be a nonsingular $M$-matrix, $W\ge0$, $\beta\ge0$, and let $\Phi:\mathbb R\to[0,1]$ be $C^1$
with $0\le \Phi'\le L$, where $L>0$. Apply $\Phi$ coordinatewise and let
$h:\mathbb R\to\mathbb R^{n_S}$ be continuous and $T$-periodic, with $T>0$. If
\[
 \bar R_0=\beta\rho(V^{-1}W)<\frac1L,
\]
then
\[
 \dot x=-Vx+\Phi(\beta Wx+h(t))
\]
has a unique $T$-periodic solution $x_{\rm per}(t)$ in the nonnegative orthant. Every nonnegative
solution satisfies
\[
 \|x(t)-x_{\rm per}(t)\|_{\infty,p}
 \le e^{-ct}\|x(0)-x_{\rm per}(0)\|_{\infty,p},
\]
where $p=(V-\beta LW)^{-1}\mathbf1\gg0$, $c=\min_i p_i^{-1}>0$ and
$\|z\|_{\infty,p}=\max_i|z_i|/p_i$. The forcing need not be small or slowly varying.
This is the standard entrainment consequence of contraction
\cite{RussoDiBernardoSontag2010,Bullo2026}.
\eeC
\begin{proof}
As in the proof of Theorem~\ref{t:KBsharp}, the Jacobian is Metzler and obeys
\[
 J(x,t)=-V+\beta\operatorname{Diag}(\Phi'(\beta Wx+h(t)))W\le -V+\beta LW,
 \qquad J(x,t)p\le-\mathbf1\le-cp.
\]
The weighted maximum-norm estimate therefore holds for any pair of solutions with the same stimulus.
Put $q=V^{-1}\mathbf1\gg0$. On $x_i=0$, $x\ge0$, the vector field points into the orthant. On the
upper $i$th face of $[0,q]$,
\[
 f_i(x,t)\le-(Vq)_i+1=0,
\]
because $-V$ is Metzler. Thus $[0,q]$ is invariant, and its time-$T$ map has Lipschitz constant at
most $e^{-cT}<1$. Banach's theorem gives a periodic solution. Bounded production ensures global
existence, and contraction proves attraction from every nonnegative initial state and uniqueness
among all nonnegative $T$-periodic solutions. The minimal period may divide $T$.
\end{proof}

\paragraph{Arbitrary multistability with convergence of every trajectory.}
Theorem~2 of \cite{OfirOvseevichMargaliot2024} complements Theorem~\ref{t:KBsharp} by providing
convergence without uniqueness. For
\[
 \dot x=-dx+\Phi(\beta Wx),\qquad d>0,\qquad \Phi\ \text{bounded},\quad 0\le \Phi'\le L,
\]
their sufficient condition for $2$-contraction is
\[
 \beta^2L^2\bigl(\sigma_1(W)^2+\sigma_2(W)^2\bigr)<2d^2,
\]
where $\sigma_1(W)\ge\sigma_2(W)$ are the two largest singular values. Since $\Phi$ is bounded,
every trajectory then converges to an equilibrium, which may depend on the initial state.

This condition is compatible with $\bar R_0>1/L$. Let $W=I+C$, with $C$ the cyclic permutation matrix
of three coordinates, so that $\rho(W)=\sigma_1(W)=2$ and $\sigma_2(W)=1$. If
\[
 \frac12<\frac{\beta L}d<\sqrt{\frac25},
\]
then, for every $m\ge2$, Theorem~\ref{t:KBsharp} provides a common smooth distribution function
with $0<\Phi<1$ and $0<\Phi'\le L$ giving at least $m$ strictly positive locally exponentially stable
equilibria, and every nonnegative trajectory converges to an equilibrium. For rank-one
$W=w\,ee^\top$ with $e\gg0$, convergence needs no such criterion: $z=e^\top x$ satisfies a scalar
autonomous equation, and each coordinate is then driven by a convergent production term.

\section{Invariant partitions preserve stability and lift bifurcations}\lbl{s:KBred}
The reduction below replaces species with equal state values by one
variable per group. We first state the invariance condition, then show
when reduced equations determine full equilibrium stability and
bifurcations, and finally apply the condition to shared threshold responses.

\beD[Invariant partitions give exact equations for group-constant states][][LinRed]
Consider $\dot x=f(x;\theta)$ for $n_S$ species, with
$x\in\mathbb R^{n_S}$ and $f$ locally Lipschitz in $x$.
Partition $\{1,\ldots,n_S\}$ into $n_G$ nonempty groups
$C_1,\ldots,C_{n_G}$. Define the characteristic matrix $H$ and
averaging matrix \rev{$H^{+}$, written $L$ in the submitted version, where $L$ is
also the slope bound,} by
\be{HL}
\begin{gathered}
 H\in\mathbb R^{n_S\times n_G},\qquad
 H_{ia}=\begin{cases}1,&i\in C_a,\\0,&i\notin C_a,\end{cases}
 \qquad H^{+}=(H^\top H)^{-1}H^\top,\\
 (H^{+}x)_a=\frac{1}{|C_a|}\sum_{i\in C_a}x_i,\qquad
 H^{+}H=I_{n_G},\qquad H\mathbf1_{n_G}=\mathbf1_{n_S}.
\end{gathered}
\ee
Here $|C_a|$ is the size of group $a$, and $I_d$ and $\mathbf1_d$
are the identity matrix and all-ones vector of dimension $d$.
For $y\in\mathbb R^{n_G}$,
$(Hy)_i=\sum_{b=1}^{n_G}H_{ib}y_b=y_a$ when $i\in C_a$.
Thus $x=Hy$ assigns the common value $y_a$ to every species in $C_a$.
The associated \emph{synchrony subspace} is
\be{syncSpace}
 \Delta_C=\{x:x_i=x_j\text{ whenever }i,j\in C_a,
                    \ a=1,\ldots,n_G\}=\operatorname{im}H.
\ee
It is invariant if solutions starting there remain there.
This holds precisely when
\be{groupInv}
 f(Hy;\theta)=H\widehat f(y;\theta),
 \qquad\text{where}\qquad \widehat f(y;\theta):=H^{+}f(Hy;\theta),
\ee
for all admissible $y,\theta$. Equivalently, the derivatives
$f_i(Hy;\theta)$ are equal for all $i$ in each group.
The common values then satisfy $\dot y=\widehat f(y;\theta)$.
The groups need not be siphons, and invariance does not imply
convergence to $\Delta_C$ from outside it.
\eeD

The matrix form of the balance condition used below is explicit.
For $A\in\mathbb R^{n_S\times n_S}$,
\[
 AH=H\widehat A
 \quad\Longleftrightarrow\quad
 \sum_{j\in C_b}A_{ij}=\widehat A_{ab}
 \quad(i\in C_a,\quad a,b=1,\ldots,n_G).
\]
Thus the sum over each column group is independent of the row chosen
within a row group. This is the equitable (or weighted balanced)
partition condition for a matrix. It is used for $A$ in the next theorem
and for the coefficient matrices $V,W$ in the threshold application below.
General balanced partitions and synchrony subspaces are studied in
\cite{SGP03,golubitsky2005patterns}. For row-stochastic threshold networks
the exact quotient reduction is in \cite{Boud26}.

\beT[Invariant partitions preserve the spectral abscissa of Metzler matrices and lift saddle-node and cusp bifurcations under spectral conditions]
[stability and bifurcation lifting][Spe]
The result applies to invariant reductions of general vector fields, with the matrix and cooperativity hypotheses specified in each item. Let $H,H^{+}$ be the characteristic and averaging matrices of
Definition~\ref{d:LinRed}. Write $J(x;\theta):=\partial_x f(x;\theta)$ and, assuming
\eqref{groupInv}, let $\widehat J(y;\theta):=\partial_y\widehat f(y;\theta)$
denote the Jacobian of the reduced vector field. Use the spectral abscissa $s$ of Theorem~\ref{t:InvF}.
\begin{enumerate}
\item Assume that $A\in\mathbb R^{n_S\times n_S}$ is Metzler and
that there exists a matrix $\widehat A\in\mathbb R^{n_G\times n_G}$
satisfying $AH=H\widehat A$.
Then $\widehat A=H^{+}AH$ is Metzler, and
\begin{equation}\label{eq:spectral-equality}
 \|e^{At}\|_\infty=\|e^{\widehat A t}\|_\infty
 \quad(t\ge0),\qquad s(A)=s(\widehat A)
\end{equation}
(no irreducibility is required).
\item For a $C^1$ family satisfying \eqref{groupInv},
$J(Hy;\theta)H=H\widehat J(y;\theta)$.
At a lifted equilibrium with Metzler $J$, reduced and full local
exponential stability are equivalent. If the full system is cooperative
on a forward-invariant box $[0,H\bar y]$, $\bar y\gg0$, its least
and greatest equilibria are the lifts of the reduced ones. Full and
reduced uniqueness are equivalent, and uniqueness implies convergence
of every trajectory in the box.
\item Suppose the family is $C^k$, $k\ge3$, and the reduced system
has a nondegenerate saddle-node (a nonzero quadratic coefficient and a transversal parameter on its one-dimensional center manifold) or a generically unfolded cusp (zero quadratic coefficient, nonzero cubic coefficient and a rank-two parameter unfolding on its one-dimensional center manifold).
If the full Jacobian is irreducible Metzler and
$s(\widehat J)=0$, the same bifurcation occurs in the full system
at the same parameter values, with strictly stable complementary
 directions. More generally this conclusion holds for any fixed
full-column-rank $H$ defining an invariant reduction, provided the full
Jacobian has a simple zero eigenvalue and all other eigenvalues have
negative real parts. Stability of the center branches then determines
full local stability.
\end{enumerate}
\eeT
\begin{proof}
Off-diagonal entries of $\widehat A=H^{+}AH$ are averages of sums of
nonnegative entries between distinct groups, hence are nonnegative.
\rev{From $AH=H\widehat A$, induction gives $A^kH=H\widehat A^k$ for every
$k\ge0$, so the exponential series gives the intertwining identity
\begin{equation}\label{eq:intertwine}
 e^{At}H=He^{\widehat At}\qquad(t\ge0).
\end{equation}}
The intertwining identity and positivity give
\[
 \|e^{At}\|_\infty=\|e^{At}\mathbf1_{n_S}\|_\infty
 =\|He^{\widehat A t}\mathbf1_{n_G}\|_\infty
 =\|e^{\widehat A t}\|_\infty,
\]
since $\|Hz\|_\infty=\|z\|_\infty$. Taking exponential growth rates
proves item 1. Differentiating \eqref{groupInv} proves the
Jacobian identity and hence the stability assertion in item 2.
For its remaining claims, Proposition~\ref{p:KBpositive}(a) applies
to a cooperative invariant box (its order argument uses only these
hypotheses). The trajectories from $0$ and $H\bar y$ stay in
$\operatorname{im}H$ and converge to the full extremal equilibria;
these are therefore the lifts of the reduced extrema. Comparison
proves the uniqueness and convergence claims.

For item 3, item 1 and irreducible Metzler theory give a simple zero
and a strictly stable complement; the more general alternative
assumes these properties directly. In coordinates $x=Hy+Zz$,
where $Z\in\mathbb R^{n_S\times(n_S-n_G)}$ has columns
complementing $\operatorname{im}H$ and
$z\in\mathbb R^{n_S-n_G}$ measures departures from that subspace, include
$\dot\theta=0$. Invariance for every nearby parameter gives the
linearization
\[
 \begin{pmatrix}
 \widehat J&\partial_\theta\widehat f&*\\
 0&0&0\\
 0&0&A_\perp
 \end{pmatrix}
 \quad\text{in the order }(y,\theta,z),
\]
with $A_\perp$ Hurwitz. The full generalized center subspace lies
in $z=0$. A center manifold of the extended reduced system therefore
lifts to a full center manifold with unchanged dynamics. The stable
complement persists nearby, proving the bifurcation and stability
claims.
\end{proof}

\paragraph{Shared responses and invariant groups.}
Consider the production model \eqref{eq:KB}. For a specified family $\mathcal S$ of siphons (for example, the minimal
siphons of a zero-onset model), impose
\begin{equation}\label{eq:KB-share-siphons}
 \Phi_i=\Phi_j\quad\text{as functions whenever }i,j\in \Sigma,\quad \Sigma\in\mathcal S.
\end{equation}
Overlaps force agreement. Let $C_1,\ldots,C_{n_G}$ be the resulting $n_G$ equivalence
classes of the $n_S$ species indices, with uncovered indices initially
singleton classes. Here $\kappa:\{1,\ldots,n_S\}\to\{1,\ldots,n_G\}$
assigns each species to its class. Equivalently,
\begin{equation}\label{eq:KB-shared-response}
 \Phi_i=\widehat\Phi_{\kappa(i)},\quad i\in C_{\kappa(i)},\qquad
 \Pi_i(x)=b_i(x_i)\widehat\Phi_{\kappa(i)}(h_i+P_i(x)).
\end{equation}
Covered classes are unions of overlapping siphons and hence siphons:
their zero faces are intersections of invariant faces. Uncovered classes
need not be siphons; further response sharing is optional. Sharing on
\emph{every} siphon would force one global response whenever the full index
set is a siphon. Preserving the response zero sets preserves the chosen
family (Theorem~\ref{t:InvF}).

\beP[Sharing responses preserves the general theory but does not imply reduction]
[shared responses and theorem scope][KBshared]
Consider \eqref{eq:KB-shared-response} with the positivity and
invariance assumptions of Proposition~\ref{p:KBpositive}.
\begin{enumerate}
\item For constant coefficients $b_i>0$, linear polynomials
$P_i=\beta(Wx)_i$, $\beta\ge0$, $W\ge0$, and
$0\le\widehat\Phi_a'\le L_a$, a sufficient condition for exponential contraction is
\begin{equation}\label{eq:KB-shared-contraction}
 \beta\rho\!\left(
 \operatorname{Diag}(b_iL_{\kappa(i)})WV^{-1}\right)<1.
\end{equation}
\item On a siphon face $x_\Sigma=0$, write $y=x_{\Sigma^c}$ for the resident
coordinates. If $P_i(0,y)=0$ for $i\in \Sigma$ and
$\widehat\Phi_{\kappa(i)}(h_i)=0$, the transverse matrix is
\begin{equation}\label{eq:KB-shared-transverse}
 J^\perp(y)=-V_\Sigma+
 \operatorname{Diag}\!\left(
 b_i(0)\widehat\Phi_{\kappa(i)}'(h_i):i\in \Sigma\right)
 \partial_{x_\Sigma}P_\Sigma(0,y).
\end{equation}
Derivatives at a boundary argument are understood from the admissible
side when needed. If the response, stimulus and $b_i(0)$ are common
on $\Sigma$, their scalar factor is common; the remaining matrix need
not have constant rows or equal eigenvalues.
\item Use the characteristic matrix $H$ of these response classes.
Sufficient conditions for \eqref{groupInv} are
\[
 VH=H\widehat V,\qquad
 b_i(r)=\widehat b_a(r),\quad h_i=\widehat h_a,\quad
 P_i(Hy)=\widehat P_a(y)\quad(i\in C_a).
\]
For linear $P(x)=\beta Wx$, the two matrix conditions are
\[
 \sum_{j\in C_b}V_{ij}=\widehat V_{ab},\qquad
 \sum_{j\in C_b}W_{ij}=\widehat W_{ab}
 \quad(i\in C_a),
\]
with sums independent of $i$ within each group; then
$\widehat P(y)=\beta\widehat Wy$.
Thus variables initially equal within each $C_a$ remain equal,
and their common values satisfy
\[
 \dot y_a=-(\widehat V y)_a+
 \widehat b_a(y_a)\widehat\Phi_a(\widehat h_a+\widehat P_a(y)).
\]
When the system is cooperative, Theorem~\ref{t:Spe}(2) transfers local
exponential stability exactly between a reduced equilibrium and its lift.
Without cooperativity, exact reduction alone gives no such conclusion.
\end{enumerate}
\eeP
\begin{proof}
The contraction assertion is a case of the multivariable Lur'e theorem of
\cite{Davydov}. Directly, put
$B=\beta\operatorname{Diag}(b_iL_{\kappa(i)})W$ and $M=B-V$.
Here $\beta\ge0$. The spectral inequality gives $s(M)<0$.
Since $M$ is Metzler and Hurwitz,
$q=(-M)^{-1}\mathbf1_{n_S}\gg0$ satisfies
$Mq=-\mathbf1_{n_S}\le-cq$,
where $c=1/\max_iq_i>0$. Moreover $J(x)\le M$ entrywise and $J(x)$ is Metzler.
The mean-value formula on the convex box therefore gives, for two solutions,
\[
 D^+\|x(t)-\widetilde x(t)\|_{\infty,q}
 \le-c\|x(t)-\widetilde x(t)\|_{\infty,q},
 \qquad \|w\|_{\infty,q}=\max_i |w_i|/q_i.
\]
Integration proves exponential contraction.

Differentiating the production gives
\[
 \partial_j\Pi_i
 =\delta_{ij}b_i'(x_i)\widehat\Phi_{\kappa(i)}(h_i+P_i(x))
 +b_i(x_i)\widehat\Phi_{\kappa(i)}'(h_i+P_i(x))\partial_jP_i(x).
\]
The first term vanishes on the specified face, proving
\eqref{eq:KB-shared-transverse}.
For the reduction, direct substitution gives \eqref{groupInv}.
Uniqueness then preserves equality within each group, and
Theorem~\ref{t:Spe}(2) gives the stability transfer.
\end{proof}
\beR[Shared responses and correlated thresholds]
[shared thresholds][KBsharenote]
Sharing a response means sharing a threshold \emph{distribution};
it does not mean that all variables use the same random threshold
realization or the same noise. Such correlated updates require a
separate microscopic model. For disjoint minimal siphons of a linear
network the sharing requirement is particularly simple, but even
then unequal removal coefficients or unequal row sums can prevent
reduction. For example, $V=\operatorname{Diag}(1,2)$,
$P_1=x_2$, $P_2=x_1$, $h=0$, and a common zero-onset response
have the minimal siphon $\{1,2\}$, whereas at $x_1=x_2=r>0$,
$f_1=-r+\Phi(r)$ and $f_2=-2r+\Phi(r)$ are unequal.
\eeR

\beR[Response shape, zero sets and sharpness]
[shape and zeros][KBshapenote]
Changing response shapes while preserving their zero sets preserves
siphons, but can change slopes in \eqref{eq:KB-shared-transverse},
contraction constants, equilibria, saddle-nodes and cusps. For a globally $C^1$ zero-onset response,
Remark~\ref{r:KBonset} gives $\widehat\Phi_a'(0)=0$; hence zero stimulus
and $P_i$ vanishing on the face give $J^\perp=-V_\Sigma$.
The common-response sharpness theorem
\ref{t:KBsharp} already satisfies any prescribed equality among
responses, but its full-support perturbation removes all nonempty
siphons. It must not be used to assert arbitrary multistability
while simultaneously preserving a prescribed nontrivial siphon
family. For SIS the common identity response is already shared
by every index; restrictions on its infection coefficients and
susceptible factors remain essential.
\eeR

\section{Equilibria, folds, and multistability for threshold networks}\label{s:KBmulti}

\subsection{Autonomous response networks undergo folds and cusps}
The derivative tests below are the standard scalar reduction of an
equilibrium equation with a simple zero eigenvalue
\cite{kuznetsov2004}. Theorem~\ref{t:Spe}(3) transfers these bifurcations from invariant
group reductions when the complementary directions are stable.
On an invariant siphon face,
full stability additionally requires stability of the absent directions.

\beP[Equilibria admit saddle-node and cusp tests, and nondegenerate cusps persist]
[][sadCus]\label{p:cusprobust}
Let $f$ be $C^2$ near an equilibrium $(x_*,\theta_*)$ of
$\dot x=f(x;\theta)$, with $x\in\mathbb R^{n_S}$, $\theta\in\mathbb R^p$.
Assume $J=\partial_xf$ has a simple zero eigenvalue and all its other
eigenvalues have negative real parts.
A saddle-node bifurcation is also called a fold bifurcation.
Evaluate all derivatives below
at this equilibrium and choose
\[
 Jr=0,\quad \ell^\top J=0,\quad \ell^\top r=1,\qquad
 Q=I-r\ell^\top,\qquad \mathcal B=\partial_x^2f.
\]
\BEN \im
Varying $\theta_j$ gives a nondegenerate saddle-node if
\[
 a_j:=\ell^\top f_{\theta_j}\ne0,\qquad
 b_2:=\ell^\top\mathcal B[r,r]\ne0.
\]
\im
For the cusp test assume additionally that $f$ is $C^3$ and put
$\mathcal C=\partial_x^3f$.
If $b_2=0$, define uniquely
\[
 Jh_{20}=-\mathcal B[r,r],\quad\ell^\top h_{20}=0,
 \qquad Jh_j=-Qf_{\theta_j},\quad\ell^\top h_j=0.
\]
In this case, a nondegenerate cusp is unfolded by $(\theta_1,\theta_2)$ if
\[
 b_3:=\ell^\top\bigl(\mathcal C[r,r,r]
                  +3\mathcal B[r,h_{20}]\bigr)\ne0,
 \qquad
 \det\begin{pmatrix}a_1&a_2\\d_1&d_2\end{pmatrix}\ne0,
\]
where $d_j=\ell^\top\bigl((\partial_xf_{\theta_j})r+
\mathcal B[r,h_j]\bigr)$.
\im
Under the cusp conditions above, every sufficiently small $C^3$
perturbation of the vector-field family has a unique nearby cusp,
with shifted state and unfolding parameters $(\theta_1,\theta_2)$.
The cubic and unfolding conditions and the strictly stable complementary
spectrum persist. If the family depends smoothly on additional
perturbation parameters, the continued cusp depends smoothly on them
to the order supplied by the implicit-function theorem.
The continuation remains admissible when the original state and
unfolding parameters are interior to their admissible domains;
at a parameter boundary, admissibility must be checked separately.
\EEN
In all three assertions, parameters other than the two designated
unfolding parameters (or the single saddle-node parameter) are held fixed.
\eeP
\begin{proof}
The projection $Q$ has range $\ker\ell^\top$, which $J$ preserves,
and $J$ is invertible on that range. The implicit-function theorem
therefore solves
\[
 Qf(x_*+\xi r+w;\theta)=0,\qquad w\in\ker\ell^\top,
\]
uniquely as $w=w(\xi,\theta)$ near $(0,\theta_*)$. Put
\[
 g(\xi,\theta)=\ell^\top
 f(x_*+\xi r+w(\xi,\theta);\theta).
\]
Nearby equilibria correspond exactly to $g=0$. At the base point,
$w=0$, $w_\xi=0$, and differentiation of the projected equation gives
\[
 Jw_{\xi\xi}=-Q\mathcal B[r,r],\qquad
 Jw_{\theta_j}=-Qf_{\theta_j}.
\]
If $b_2=0$, then $Q\mathcal B[r,r]=\mathcal B[r,r]$, proving both
solvability and uniqueness of $h_{20}$; also $w_{\theta_j}=h_j$.
Since $\ell^\top J=0$, the scalar derivatives are
\[
 g_{\theta_j}=a_j,\qquad g_{\xi\xi}=b_2,
 \qquad g_{\xi\theta_j}=d_j,
 \qquad g_{\xi\xi\xi}=b_3\quad(b_2=0).
\]
These are the usual saddle-node and cusp tests for the equilibrium equation
\cite{kuznetsov2004}. The stable complementary spectrum gives the
corresponding differential-equation bifurcations. In particular, $g$
is obtained by solving the complementary equilibrium equations; the
center-manifold differential equation has the same nondegeneracy tests.

For persistence, apply the implicit-function theorem to
$\mathcal H=(g,g_\xi,g_{\xi\xi})$ in the variables
$(\theta_1,\theta_2,\xi)$. At the cusp,
\[
 D\mathcal H=
 \begin{pmatrix}a_1&a_2&0\\d_1&d_2&0\\ *&*&b_3\end{pmatrix},
 \qquad
 \det D\mathcal H=b_3
 \det\begin{pmatrix}a_1&a_2\\d_1&d_2\end{pmatrix}\ne0.
\]
The same complementary-equation reduction is available for nearby
families. Small $C^3$ perturbations give small $C^1$ perturbations
of $\mathcal H$, so its locally unique zero persists. Continuity
preserves the strict nondegeneracy and spectral conditions, and
interior admissibility is open.

For linear $P_i$ in \eqref{polR}, with $\varphi_i=\Phi_i'$,
the required state derivatives are
\[
 \mathcal B_i[r,r]=b_i\beta^2\varphi_i'(h_i+\beta(Wx)_i)(Wr)_i^2,
 \qquad
 \mathcal C_i[r,r,r]=b_i\beta^3\varphi_i''(h_i+\beta(Wx)_i)(Wr)_i^3.
\]
The correction terms use the same bilinear form:
\[
 \mathcal B_i[r,z]
 =b_i\beta^2\varphi_i'(h_i+\beta(Wx)_i)(Wr)_i(Wz)_i.
\]
\end{proof}

\subsection{Odds-polynomial responses: a siphon-preserving cusp}
\label{ss:oddsCusp}

The constructions above give two distinct mechanisms: flat zero-onset
responses preserve siphons but have vanishing onset slope, whereas the
full-support responses in the explicit two-population cusp calculations
exclude nontrivial zero faces. Polynomial odds provide a common class in
which the response zero set is preserved while its onset slope and
nonlinear shape can be varied independently.

We use the cusp continuation criterion of Proposition~\ref{p:sadCus}
and persistence of hyperbolic equilibria under weak coupling
\cite{golubitsky2023dynamics,kuznetsov2004,Boud26}.
The constructions below additionally preserve specified siphon faces
and control their transverse stability.

\beD[Odds-polynomial threshold response]
[][Odd]
Fix $m\ge1$ and $a_1,\ldots,a_m\ge0$ with $a_m>0$, and put
\begin{equation}\label{eq:oddsFamily}
 \varphi(q)=\sum_{k=1}^m a_kq^k,\qquad
 \Phi(q)=\begin{cases}\dfrac{\varphi(q)}{1+\varphi(q)},&q\ge0,\\[4pt] 0,&q\le0.\end{cases}
\end{equation}
Then $\Phi$ is a nondecreasing distribution function, is strictly positive
for $q>0$, tends to one as $q\to\infty$, and has zero set $(-\infty,0]$.
It is smooth on $[0,\infty)$ from the right, with
\[
 \Phi'(0+)=a_1,\qquad \Phi''(0+)=2(a_2-a_1^2).
\]
If $a_1>0$, its extension by zero to negative arguments is not $C^1$ at
zero. Derivatives at zero below are taken in the admissible nonnegative
directions; local cusp calculations are performed at strictly
positive arguments, where the vector field is smooth. Since $1+\varphi>0$ on a
neighbourhood of $[0,\infty)$, the formula $\varphi/(1+\varphi)$ is a smooth
extension of $\Phi$ there, so threshold models built from $\Phi$ with nonnegative
polynomials $P_i$ satisfy the smoothness convention of Definition~\ref{d:KB}
(Remark~\ref{r:KBonset}).
\eeD

\beP[Odds responses keep the siphon lattice and scale the invasion matrix][odds invasion matrix][Odd]
Consider a production threshold model \eqref{eq:KB}, \eqref{polR} with $h_i=0$, a common
response $\Phi$, and polynomials with $P_i(0)=0$. Keep $V$, $b_i$, and $P_i$ fixed while varying the
admissible coefficients of $\Phi$ in \eqref{eq:oddsFamily}.
Write $z=x_\Sigma$ and $y=x_{\Sigma^c}$.
A nonempty set $\Sigma$ is a siphon if and only if
\[
 P_i(0,y)=0\quad(i\in\Sigma)
 \quad\text{for every admissible resident state }y,
\]
and
\begin{equation}\label{eq:oddsNoLinCoupling}
 V_{\Sigma\Sigma^c}=0.
\end{equation}
Consequently the entire siphon lattice is independent of the odds
coefficients. On such a face the transverse Jacobian is
\begin{equation}\label{eq:oddsTransverse}
 J_\Sigma^\perp(y)
 =-V_\Sigma
 +a_1\operatorname{Diag}(b_i:i\in\Sigma)D_zP_\Sigma(0,y).
\end{equation}
Thus only $a_1$ enters the transverse matrix at a \emph{fixed} resident
state $y$, and it changes that matrix only if
$D_zP_\Sigma(0,y)\ne0$. Higher odds coefficients may change invasion
thresholds indirectly, by moving the resident equilibrium at which
\eqref{eq:oddsTransverse} is evaluated.
\eeP
\begin{proof}
The polynomials $P_i$ are nonnegative on the state space, so there $\Phi$ may be replaced by
$\varphi/(1+\varphi)$, which is smooth near $q=0$; this gives a $C^1$ extension of the vector field,
with $\Pi\ge0$ on the admissible state space. Theorem~\ref{t:InvF}(1) then applies: $\Sigma$ is a siphon exactly when
$V_{\Sigma\Sigma^c}=0$ and $b_i\Phi(P_i(0,y))=0$ for $i\in\Sigma$. Since $b_i>0$ and $\Phi(q)>0$ for
$q>0$, the second condition is $P_i(0,y)=0$, which proves the equivalence and independence from
the response coefficients. Since $\Phi'(0+)=a_1$, the chain rule gives
$\partial_z\Pi_\Sigma(0,y)=a_1\operatorname{Diag}(b_i)D_zP_\Sigma(0,y)$, and the formula for
$J^\perp$ in Theorem~\ref{t:InvF}(1) proves \eqref{eq:oddsTransverse}.
\end{proof}

We next construct a quartic odds response having a positive cusp
together with a stable boundary equilibrium.

\beT[Scalar threshold models with polynomial odds require degree at least four for a stable origin and a nondegenerate positive cusp][minimal degree and quartic odds cusp][Qua]\label{p:Min}
For $\dot x=-x+\Phi(x)$ with $\Phi$ as in Definition~\ref{d:Odd},
degree four of the odds polynomial $\varphi$ is the minimum admitting both a hyperbolically stable
origin and a positive nondegenerate cusp, in the sense of Proposition~\ref{p:sadCus}(2). The following quartic
family attains it and also has a second hyperbolically stable equilibrium.

Choose $0<r<s<1$ with $3r+s\le1$, and define
\[
 c=\frac{1}{(1-r)^3(1-s)},\qquad
 \varphi_{r,s}(q)=\frac{q-cq(q-r)^3(q-s)}{1-q}.
\]
The apparent singularity at $q=1$ is removable: the expression is the
quartic polynomial
\[
 \begin{aligned}
 \varphi_{r,s}(q)&=a_1q+a_2q^2+a_3q^3+a_4q^4,\\
 a_1&=1-cr^3s,\\
 a_2&=c(1-3r-s+3r^2+3rs),\\
 a_3&=c(1-3r-s),\qquad a_4=c.
 \end{aligned}
\]
The coefficients $a_1,a_2,a_4$ are positive, and $a_3\ge0$, with equality
exactly when $3r+s=1$; thus $3r+s\le1$ is precisely the admissibility
condition $a_k\ge0$. With the response $\Phi_{r,s}$
constructed from $\varphi_{r,s}$, the scalar equation
$\dot x=f_{r,s}(x)=-x+\Phi_{r,s}(x)$ on $[0,1]$ satisfies
\[
 f_{r,s}(x)=-\frac{cx(x-r)^3(x-s)}{1+\varphi_{r,s}(x)}.
\]
It has a hyperbolically stable origin, a nondegenerate positive cusp
at $x=r$, and a hyperbolically stable positive equilibrium at $x=s$.
More precisely,
\[
 f_{r,s}'(0)=-cr^3s<0,\qquad
 f_{r,s}'''(r)=6cr(s-r)(1-r)>0,
\]
\[
 f_{r,s}'(s)=-cs(s-r)^3(1-s)<0.
\]
The cusp is generically unfolded by $(a_1,a_2)$, keeping $a_3,a_4$
fixed, since for any quartic odds polynomial
\begin{equation}\label{eq:oddsUnfold}
 \det\frac{\partial(f,f_x)}{\partial(a_1,a_2)}
 =\frac{x^2}{(1+\varphi(x))^4}>0\qquad(x>0).
\end{equation}
In particular, at $x=r$ the determinant equals $r^2(1-r)^4$.

The root parameters control the boundary stability margin and the cusp
strength,
\begin{equation}\label{eq:rsMargin}
 1-a_1=c\,r^3s,\qquad f_{r,s}'''(r)=6cr(s-r)(1-r),
\end{equation}
and over the admissible region
\begin{equation}\label{eq:rsSup}
 \sup_{0<r<s,\ 3r+s\le1}\,(1-a_1)=\frac1{81}.
\end{equation}
The supremum is not attained. It is approached only as $r,s\to1/4$,
where $f_{r,s}'''(r)\to0$ and the cusp degenerates. For instance, $(r,s)=(1/5,2/5)$ has $a_3=0$ and margin $1/96$.

For $r=1/8$ and $s=1/2$ this is exactly the reference polynomial
\[
 \varphi(q)=\tfrac1{343}(342q+368q^2+128q^3+1024q^4),
\]
for which
\begin{equation}\label{eq:oddsFactor}
 f(x)=-\frac{x(2x-1)(8x-1)^3}{343(1+\varphi(x))}.
\end{equation}
Here $f'(0)=-1/343$, $f'(1/2)=-27/686$, and
\[
 f=f_x=f_{xx}=0,\qquad f_{xxx}=\frac{36}{49}
 \quad\text{at }x=\frac18.
\]
At this cusp the unfolding determinant in \eqref{eq:oddsUnfold}
is $2401/262144$, and
\[
 \det\frac{\partial(f,f_x,f_{xx})}{\partial(x,a_1,a_2)}
 =\frac{441}{65536}\ne0
 \quad\text{at }(x,a_1,a_2)
 =\left(\tfrac18,\tfrac{342}{343},\tfrac{368}{343}\right).
\]
\eeT
\begin{proof}
\emph{Minimum degree.}
The numerator $(1-x)\varphi(x)-x$ has degree at most $m+1$,
and contains the factor $x$. A positive triple equilibrium thus
requires $m\ge3$. If $m=3$, its leading coefficient is $-a_3$
and necessarily
\[
 (1-x)\varphi(x)-x=-a_3x(x-r)^3,\qquad r>0.
\]
Comparison of the coefficients of $x$ yields
$a_1-1=a_3r^3>0$. Thus the origin is unstable.
For $m\le2$ the numerator cannot contain both $x$ and a
positive triple root. This proves the lower bound; the construction
below proves attainment.

\emph{Quartic construction.}
The choice of $c$ makes
$q-cq(q-r)^3(q-s)$ vanish at $q=1$, so division by $1-q$
gives a polynomial of degree four with the displayed coefficients.
Since $c>0$, $a_3\ge0$ is $3r+s\le1$, and
$a_2=a_3+3cr(r+s)>0$. Also $r<s\le1-3r$ gives $r<1/4$ and
$1-s\ge3r$, which imply $(1-r)^3(1-s)>r^3s$, hence $a_1>0$;
finally $a_4=c>0$.
The identity
$(1-x)\varphi_{r,s}(x)-x=-cx(x-r)^3(x-s)$ proves the
general factorization and hence its specialization
\eqref{eq:oddsFactor}. The three derivatives follow by factoring at
the respective equilibria and using
$1+\varphi_{r,s}(r)=(1-r)^{-1}$ and
$1+\varphi_{r,s}(s)=(1-s)^{-1}$.
For the unfolding determinant, differentiating
$f=(\varphi-x-x\varphi)/(1+\varphi)$ with respect to $a_1,a_2$
and then $x$ gives \eqref{eq:oddsUnfold}; at the cusp,
$f=f_x=f_{xx}=0$ and $f_{xxx}\ne0$. Hence the Jacobian
of $(f,f_x,f_{xx})$ in $(x,a_1,a_2)$ equals
$f_{xxx}$ times the unfolding determinant, and is nonzero. Proposition~\ref{p:sadCus} in one dimension,
therefore gives the claimed generically unfolded cusp.
Since $f_{r,s}'''(r)>0$, the cusp equilibrium itself is unstable;
the two stable equilibria are $0$ and $s$.
\emph{Stability margin.}
The margin identity in \eqref{eq:rsMargin} is $f_{r,s}'(0)=a_1-1$.
For \eqref{eq:rsSup}, the margin
$\Delta(r,s)=1-a_1=r^3s/((1-r)^3(1-s))$ increases with $s$, so
$\Delta(r,s)\le\Delta(r,1-3r)=r^2(1-3r)/(3(1-r)^3)$, and
\[
 \frac{d}{dr}\,\frac{r^2(1-3r)}{3(1-r)^3}=\frac{2r(1-4r)}{3(1-r)^4}>0
 \qquad(0<r<\tfrac14).
\]
The bound therefore increases to its value $1/81$ at $r=1/4$, which is not
admissible since there $s=1-3r=r$. Along any sequence with
$\Delta\to1/81$ both $r$ and $s$ tend to $1/4$, so $s-r\to0$ and
$f_{r,s}'''(r)\to0$.
Substitution gives the stated exact values.
\end{proof}

\paragraph{The quartic odds construction has an explicit scalar energy.}
For $\dot x=-x+\Phi_a(x)$ with an odds response $\Phi_a$ of Theorem~\ref{t:Qua}, $\Phi_a$ is strictly
increasing on $[0,\infty)$ when $a_1>0$, maps that interval onto $[0,1)$, and has $\Phi_a^{-1}(0)=0$. Define
\[
 \mathcal E_a(x)=\int_0^x\Phi_a^{-1}(s)\,ds-\frac{x^2}{2},\qquad 0\le x<1.
\]
Then $\dot{\mathcal E}_a(x)=-(\Phi_a^{-1}(x)-x)(x-\Phi_a(x))\le0$, with equality exactly at scalar equilibria.
At a positive equilibrium $e$, $\mathcal E_a''(e)=1/\Phi_a'(e)-1$, so a hyperbolically stable equilibrium
is a strict local minimum of this energy. At a scalar cusp $r$, where $\Phi_a(r)=r$, $\Phi_a'(r)=1$ and
$\Phi_a''(r)=0$, differentiating the inverse gives
\[
 \mathcal E_a'(r)=\mathcal E_a''(r)=\mathcal E_a'''(r)=0,\qquad \mathcal E_a''''(r)=-\Phi_a'''(r).
\]
For the reference quartic,
\[
 \mathcal E_a''(0+)=\frac1{342}>0,\qquad
 \mathcal E_a''(1/2)=\frac{27}{659}>0,\qquad
 \mathcal E_a''''(1/8)=-\frac{36}{49}<0.
\]
The stable origin and the stable positive equilibrium are therefore local energy minima, while the
cusp is a degenerate local maximum at the reference parameter. This is the generator of the Fenchel--Bregman divergence
associated with the quartic construction, not an additional minimum-degree
claim. The coupled odds networks have nonlinear $P_i$, so this energy does not simply sum over them.

\paragraph{Comparison with rational feedback models.}
The genetic feedback model of \cite{Hasty2001} has production
\begin{equation}\label{eq:Hasty}
 m\frac{1+x^2+\alpha\sigma_1x^4}
 {1+x^2+\sigma_1x^4+\sigma_1\sigma_2x^6},
\end{equation}
which is positive at $x=0$ and does not satisfy the stable-origin
assumption of Theorem~\ref{t:Qua}; see also \cite{Isaacs2003}.
Earlier polynomial-equilibrium cusp calculations appear in
\cite[Eq.~(21)]{vanHoek}. The minimum-degree statement here concerns
only the specified odds-polynomial class.

\beT[A two-dimensional polynomial production threshold network with a quartic odds response has an interior cusp and a stable boundary equilibrium for small couplings]
[][Net]
Let $\Phi_a$ denote the odds-polynomial threshold response of Definition~\ref{d:Odd}, parametrized by $a=(a_1,a_2)$, with $a_3=128/343$ and
$a_4=1024/343$ fixed, and initially
$(a_1,a_2)=a^0=(342/343,368/343)$ as in
Theorem~\ref{t:Qua}. Consider the two-population family on
the forward-invariant square $[0,1]^2$
\begin{equation}\label{eq:oddsNet}
 \dot x=-x+\Phi_a\bigl(x(1+\gamma y)\bigr),\qquad
 \dot y=-y+\Phi_a\bigl(y(10+\varepsilon x)\bigr),
 \qquad \gamma,\varepsilon\ge0.
\end{equation}
\BEN
\item For every $a\in\mathbb R^2_+$ and every $\gamma,\varepsilon\ge0$,
 the faces $x=0$ and $y=0$ are invariant and the minimal siphons
 are exactly $\{x\}$ and $\{y\}$.
\item At $(\gamma,\varepsilon,a)=(0,0,a^0)$, the resident equation
 $\dot y=-y+\Phi_{a^0}(10y)$ has exactly one positive equilibrium
 $y^*=q^*/10$, where $q^*$ is the unique positive root of
\begin{equation}\label{eq:oddsResident}
 p(q)=1024q^4-10112q^3-912q^2-3338q-3077=0.
\end{equation}
 Its eigenvalue is strictly negative, and, writing $\varphi$ for the odds polynomial of $\Phi_{a^0}$, is given exactly by
 \[
 \lambda_y=-\frac{q^*p'(q^*)}{343(1+\varphi(q^*))}<0.
 \]
 In particular, at $(1/8,y^*)$ the full Jacobian has a simple
 zero eigenvalue and a strictly negative complementary eigenvalue.
\item There exists $\delta>0$ and smooth functions
 $a_1(\gamma,\varepsilon)$, $a_2(\gamma,\varepsilon)$,
 and $x_c(\gamma,\varepsilon)$, for
 $0\le\gamma,\varepsilon<\delta$, reducing to
 $(a_1^0,a_2^0,1/8)$ at zero coupling, where $x_c$ is the first coordinate of the continued cusp, such that
 \eqref{eq:oddsNet} has a nondegenerate interior cusp near
 $(1/8,y^*)$. Its complementary eigenvalue remains strictly
 negative, and $(a_1,a_2)$ still unfold the cusp generically.
\item Let $y_*(a_1,a_2)$ be the unique resident equilibrium
 near $y^*$ for the current coefficients. On the continued cusp
 family, the boundary equilibrium $(0,y_*)$ has transverse eigenvalue
\begin{equation}\label{eq:oddsTransEig}
 \lambda_{\mathrm{tra}}(\gamma,\varepsilon)
 =-1+a_1(\gamma,\varepsilon)
 \Bigl(1+\gamma y_*\bigl(a_1(\gamma,\varepsilon),
                             a_2(\gamma,\varepsilon)\bigr)\Bigr).
\end{equation}
 This is $-1/343$ at zero coupling and remains negative for all
 sufficiently small nonnegative couplings.
\item Along the continued cusp family, a stable boundary
 equilibrium near $(0,y^*)$ and a stable interior equilibrium
 near $(1/2,y^*)$ coexist. For $\gamma,\varepsilon>0$, both off-diagonal Jacobian entries
 are strictly positive at every positive interior state.
\EEN
Consequently, sufficiently small positive coupling preserves
both minimal siphons, a nondegenerate interior cusp, and bistability
between a stable boundary equilibrium and a stable positive equilibrium,
with a strictly stable transverse invasion eigenvalue on the boundary.
\eeT
\begin{proof}
For part 1, both arguments of $\Phi_a$ contain the corresponding
variable as a factor. The zero-onset property therefore makes each
coordinate zero face invariant. These singleton faces are the
minimal siphon faces.

For part 2, with $q=10y$ one computes
\[
 -y+\Phi_{a^0}(10y)
 =-\frac{q\,p(q)}{10\,[343+343\varphi(q)]}.
\]
The coefficients of $p$ have exactly one sign change, so by
Descartes' rule $p$ has exactly one positive root, counted with
multiplicity. Since $p(0)<0$ and $p(q)\to\infty$, that root is
simple and $p'(q^*)>0$. Differentiating the displayed identity at
$q=q^*$ gives $\lambda_y$ as stated. The equations decouple at
zero coupling, and Theorem~\ref{t:Qua} gives the zero eigenvalue
in the $x$ direction.

For part 3, at zero coupling the first equation is the scalar
cusp of Theorem~\ref{t:Qua}, unfolded by $(a_1,a_2)$; the
second equation is independent of $x$ and has eigenvalue
$\lambda_y<0$. Eliminating its equilibrium equation therefore
leaves exactly $\Phi_a(x)-x$, so the full system satisfies the cusp
tests of Proposition~\ref{p:sadCus}. Its persistence assertion
continues the cusp for small $(\gamma,\varepsilon)$, preserving
the nonzero cubic, the unfolding determinant, and the stable
complementary eigenvalue. Smoothness of this family gives the
claimed smooth continuation. The base state and odds coefficients
are strictly positive; restricting the continued family to
nonnegative couplings preserves admissibility.

For part 4, restriction to $x=0$ gives
$\dot y=-y+\Phi_a(10y)$, independent of $\gamma,\varepsilon$.
Its hyperbolic positive equilibrium continues smoothly with $a$.
The face Jacobian is triangular, with transverse entry
$-1+\Phi_a'(0+)(1+\gamma y_*)$, which is
\eqref{eq:oddsTransEig}. This is strictly negative at the base
point, so it remains so for sufficiently small couplings.

For part 5, both $(0,y^*)$ and $(1/2,y^*)$ are hyperbolically
stable at zero coupling: their tangential resident eigenvalue is
$\lambda_y<0$, the former has transverse eigenvalue $-1/343$,
and the latter has $x$ eigenvalue $-27/686$. The implicit-function
theorem and spectral continuity preserve both equilibria and their
stability under the continued parameter adjustment. The first stays
on the invariant face $x=0$; the second remains strictly positive.
For positive couplings, direct differentiation gives
$\partial_y\dot x=\gamma x\Phi_a'(x(1+\gamma y))>0$
and
$\partial_x\dot y=\varepsilon y\Phi_a'(y(10+\varepsilon x))>0$
when $x,y>0$.
\end{proof}

\beR[What the construction does and does not settle]
For fixed odds coefficients from Theorem~\ref{t:Qua}, let $y^*$
be a positive equilibrium of the corresponding resident equation
$\dot y=-y+\Phi_{r,s}(10y)$. The transverse stability condition at
$(0,y^*)$ is, by \eqref{eq:rsMargin},
\begin{equation}\label{eq:rsGamma}
 \gamma<\frac{1-a_1}{a_1y^*}=\frac{c\,r^3s}{(1-c\,r^3s)\,y^*},
\end{equation}
which is $-1+a_1(1+\gamma y^*)<0$ solved for $\gamma$. For the reference
coefficients it is $\gamma<1/(342y^*)\approx2.9\times10^{-3}$\rev{, where the
resident value is $y^*\approx0.99997$, which is what makes the bound this small}.
This is \emph{not} an explicit bound along the continued cusp family: there $a_1$ changes with
the couplings, and within the quartic family changing $(r,s)$ also changes the resident state,
the cusp location and its cubic coefficient. The supremum $1/81$ in \eqref{eq:rsSup} therefore
bounds the uncoupled origin's stability margin, not the admissible coupling. The two-parameter quartic
family of Theorem~\ref{t:Qua} permits the cusp location, the
second stable equilibrium, and the boundary stability margin to be
varied while the response zero set is fixed; its coupled continuation
requires only the same hyperbolicity and unfolding hypotheses.
The theorem gives local continuation, not a classification of all cusp
branches by their invasion thresholds.
\eeR

\beR[Why the coupling must be multiplicative for singleton siphons][multiplicative coupling][Mul]
The arguments of \eqref{eq:oddsNet} are $x(1+\gamma y)$ and $y(10+\varepsilon x)$,
not $x+\gamma y$ and $10y+\varepsilon x$. For a zero-onset response that is strictly positive on $(0,\infty)$,
diagonal removal and nonnegative polynomials $P_i$, the singleton $\{i\}$ is a
siphon exactly when
\[
 P_i(x)\big|_{x_i=0}\equiv0,
\]
that is, when every monomial of $P_i$ contains $x_i$; this is the criterion of
Theorem~\ref{t:InvF} specialised to a singleton. The additive argument
$x+\gamma y$ has the monomial $\gamma y$, which survives at $x=0$, so
$\Phi(\gamma y)>0$ and the face $x=0$ is not invariant: the siphon is destroyed for
every $\gamma>0$. The multiplicative argument $x(1+\gamma y)$ carries $x$ in every
monomial and so preserves it for all $\gamma$. A nonzero off-diagonal linear term therefore destroys the
corresponding singleton siphon.

The necessity is specific to singleton siphons. For a general $\Sigma$, with
$V_{\Sigma\Sigma^c}=0$, the criterion of Proposition~\ref{p:Odd} requires only
\[
 \operatorname{supp}\nu\cap\Sigma\ne\varnothing
 \quad\text{for every monomial }x^\nu\text{ of }P_i,\ i\in\Sigma,
\]
so a coordinate of $\Sigma$ may have an additive term from another coordinate of
$\Sigma$. Overlapping minimal siphons therefore allow a richer class of couplings.
\eeR

\subsection{Overlapping minimal siphons support distinct invasion thresholds and interior cusps}

\paragraph{Reading siphons from the production.}
\beR[Gaussian and logistic responses have no siphons]
[no siphons][KBnosiph]
Gaussian and logistic CDFs are strictly positive at every finite argument.
They therefore exclude nonempty siphons whenever $b_i(0)>0$ for all $i$.
Zero-onset responses can admit siphons, depending on $h$, $P$ and $V$;
irreducibility of the interaction graph alone does not settle the question.
\eeR
\beR[Graph rule for siphons]
[graph rule][KBsiphgraph]
For zero stimulus, linear nonnegative $P_i$ and such zero-onset
responses, Theorem~\ref{t:InvF} gives a graph description;
with nondiagonal $V$, include transfer edges $j\to i$ whenever
$-v_{ij}>0$, as well as the positive edges of $W$. Siphons are the
sets containing every predecessor of each of their indices, and
minimal siphons are the source strongly connected components of
this combined graph. This graph rule requires a response positive
at every positive argument; a response with an additional flat zero
interval on a finite box requires \eqref{eq:KB-response-zero-test}.
\eeR
\beR[The response has a definite order along transverse rays]
[transverse order][KBtransOrder]
For the threshold model on a siphon face, fix $h_i=0$, a positive resident state $y$,
and a nonzero polynomial $P_i$ with $P_i(0,y)=0$. Put
\[
 d_i=\min_{a_{i\nu}>0}\sum_{j\in \Sigma}\nu_j\ge1,\qquad
 \Phi_i(t)=c_it^{k_i}+o(t^{k_i})\quad(t\downarrow0),\quad c_i,k_i>0.
\]
For each fixed $\zeta\gg0$, the leading coefficient
$p_i(\zeta,y)>0$ satisfies
\[
 P_i(\varepsilon\zeta,y)=\varepsilon^{d_i}p_i(\zeta,y)
 +o(\varepsilon^{d_i}),\qquad
 \Pi_i(\varepsilon\zeta,y)
 =b_ic_ip_i(\zeta,y)^{k_i}\varepsilon^{d_ik_i}
 +o(\varepsilon^{d_ik_i}).
\]
Thus the response has order $d_ik_i$ along positive transverse rays.
For the model with factors $b_i$ this order also holds if $b_i(0)>0$;
a vanishing $b_i(0)$ can increase it.
The density has order $k_i-1$ only if additionally
$\Phi_i'(t)=k_ic_it^{k_i-1}+o(t^{k_i-1})$; this does not follow
by differentiating the response expansion.
Grouping monomials by their degree in $z$ gives the leading
term; its coefficient is positive for $\zeta,y\gg0$.
Substitution into the response expansion gives the displayed
production expansion. Multiplication by $b_i(0)+o(1)$ proves the
assertion for the factors $b_i$ when $b_i(0)>0$. No differentiation of an
$o(\cdot)$ remainder is used.
\eeR

\beXa[Minimal siphons can overlap for higher-order $P_i$]
[overlapping minimal siphons][KBoverlap]
On $[0,1]^3$, let $V=I$, $h_i=0$, and
\[
 P_1(x)=x_2,\qquad P_2(x)=x_1x_3,\qquad P_3(x)=x_2,
 \qquad \Pi_i(x)=\Phi_i(P_i(x)),
\]
where $0\le \Phi_i\le1$, $\Phi_i(0)=0$ and $\Phi_i(q)>0$ for $q>0$.
The minimal siphons are $\{1,2\}$ and $\{2,3\}$.
Indeed, absence of $1$ or $3$ forces inclusion of $2$, and absence
of $2$ must be accompanied by absence of $1$ or $3$ to suppress
$x_1x_3$. Neither singleton is a siphon, and $\{1,3\}$ is not one.
Sharing a response on both minimal siphons forces
$\Phi_1=\Phi_2=\Phi_3$. Thus minimal siphons do not in general form a
partition. Smooth admissible responses exist, for example
$\Phi_i(q)=\exp(-1/q)$ for $q>0$ and zero for $q\le0$.
\eeXa

\paragraph{Why Example~\ref{e:KBoverlap} has no nontrivial boundary residents.}
Consider first the unmodified model of Example~\ref{e:KBoverlap},
\[
 \dot x_1=-x_1+\Phi(x_2),\qquad
 \dot x_2=-x_2+\Phi(x_1x_3),\qquad
 \dot x_3=-x_3+\Phi(x_2),
\]
where $\Phi(0)=0$ and $\Phi(q)>0$ for $q>0$. Its minimal siphons
are $\Sigma_{12}=\{1,2\}$ and $\Sigma_{23}=\{2,3\}$.
On $x_1=x_2=0$ the resident equation is $\dot x_3=-x_3$;
on $x_2=x_3=0$ it is $\dot x_1=-x_1$.
The origin is its only boundary equilibrium: at an
equilibrium, $x_1=0$ or $x_3=0$ forces $x_2=0$ and then all
coordinates to vanish, while $x_2=0$ forces $x_1x_3=0$.
For an odds response with $\Phi'(0+)=a_1$, its Jacobian at the
origin is
\[
 J(0)=\begin{pmatrix}
 -1&a_1&0\\ 0&-1&0\\ 0&a_1&-1
 \end{pmatrix},\qquad \sigma(J(0))=\{-1,-1,-1\}.
\]
The transverse blocks at resident levels $x_3=u$ and $x_1=u$
are, respectively,
\[
 J_{12}^\perp(u)=\begin{pmatrix}-1&a_1\\a_1u&-1\end{pmatrix},
 \qquad
 J_{23}^\perp(u)=\begin{pmatrix}-1&a_1u\\a_1&-1\end{pmatrix}.
\]
Both have spectral abscissa $-1+a_1\sqrt{u}$, but the
resident satisfies $\dot u=-u$, so $u>0$ is not a resident
equilibrium. Thus the overlap is genuine, while the original
model has no nontrivial boundary equilibrium at which to
compare invasion thresholds. In particular, varying the onset
slope does not change its boundary-equilibrium invasion
eigenvalues.

We now add self-feedback without changing the two minimal
siphons. Use the quartic odds response $\Phi_a$ of
equation~\eqref{eq:oddsFamily}, with
\[
 \phi_a(q)=a_1q+a_2q^2+\frac{128}{343}q^3
                    +\frac{1024}{343}q^4,\qquad
 a^0=\left(\frac{342}{343},\frac{368}{343}\right),
\]
and let $r=1/8$, $s=1/2$. By Theorem~\ref{t:Qua},
$-u+\Phi_{a^0}(u)$ has stable roots $0,s$ and a nondegenerate
cusp root $r$. For $a=(a_1,a_2)$ near $a^0$, denote by
$s_a$ the unique root near $s$; it depends smoothly on $a$,
and
\[
 s_a>0,\qquad \lambda_s(a):=\Phi_a'(s_a)-1<0,
 \qquad a_1<1.
\]

\beT[A three-dimensional polynomial production threshold network with a quartic odds response has overlapping minimal siphons, independently variable invasion thresholds, and an interior cusp coexisting with at least eight stable equilibria][][KBoverlapCusp]\label{t:Thr}
Consider, on $[0,1]^3$, the cooperative threshold network
\begin{equation}\label{eq:KBoverlapCusp}
\begin{aligned}
 \dot x_1&=-x_1+\Phi_a(x_1+\alpha x_2),\\
 \dot x_2&=-x_2+\Phi_a(x_2+\beta x_1x_3),\\
 \dot x_3&=-x_3+\Phi_a(x_3+\gamma x_2),
\end{aligned}
\qquad \alpha,\beta,\gamma>0.
\end{equation}
Use the odds-polynomial response of Definition~\ref{d:Odd}, with $a_3=128/343$, $a_4=1024/343$ and $a=(a_1,a_2)$ near $a^0=(342/343,368/343)$ of Theorem~\ref{t:Qua}. Let $r=1/8$, $s=1/2$ be the positive cusp and stable scalar equilibrium of Theorem~\ref{t:Qua} at $a^0$, respectively, and let $s_a$ be the scalar equilibrium near $s$ satisfying $\Phi_a(s_a)=s_a$, with $\lambda_s(a):=\Phi_a'(s_a)-1<0$. Then the following conclusions hold.
\begin{enumerate}
\item The minimal siphons are exactly
$\Sigma_{12}=\{1,2\}$ and $\Sigma_{23}=\{2,3\}$;
the only other nonempty siphon is $\{1,2,3\}$.
Every equilibrium has support either $\varnothing$,
$\{1\}$, $\{3\}$, or $\{1,2,3\}$.
For $x\gg0$ the Jacobian is irreducible and Metzler.

\item Besides $E_0=0$, the equilibria
\[
 E_{12}=(0,0,s_a),\qquad E_{23}=(s_a,0,0)
\]
lie on the respective proper siphon faces. The tangential
eigenvalue at each is $\lambda_s(a)<0$. Their transverse
matrices, with absent variables ordered as $(x_1,x_2)$ and
$(x_2,x_3)$, respectively, are
\begin{equation}\label{eq:KBoverlapJtra}
 J_{12}^\perp=\begin{pmatrix}
 a_1-1&a_1\alpha\\ a_1\beta s_a&a_1-1
 \end{pmatrix},\qquad
 J_{23}^\perp=\begin{pmatrix}
 a_1-1&a_1\beta s_a\\ a_1\gamma&a_1-1
 \end{pmatrix}.
\end{equation}
Consequently the two invasion numbers and spectral abscissae, in the notation of Theorem~\ref{t:InvF}, are exactly
\begin{equation}\label{eq:KBoverlapR}
\begin{aligned}
 R_{\Sigma_{12}}&=a_1\bigl(1+\sqrt{\alpha\beta s_a}\bigr),&
 s(J_{12}^\perp)&=R_{\Sigma_{12}}-1,\\
 R_{\Sigma_{23}}&=a_1\bigl(1+\sqrt{\gamma\beta s_a}\bigr),&
 s(J_{23}^\perp)&=R_{\Sigma_{23}}-1.
\end{aligned}
\end{equation}
In particular $E_{12}$ and $E_{23}$ are locally
exponentially stable precisely when their respective
reproduction numbers are less than one. With $a$ and
$\beta>0$ fixed, the two thresholds can be varied
independently by $\alpha$ and $\gamma$; the reproduction numbers differ whenever
$\alpha\ne\gamma$. For example, the critical values are
\[
 \alpha_c=\gamma_c=
 \frac{(1-a_1)^2}{a_1^2\beta s_a}.
\]
These exact boundary-threshold formulas do not presume the
existence of an interior cusp at the critical couplings.

\item There is $\delta>0$ and a smooth map
\[
 (\alpha,\beta,\gamma)\longmapsto
 \bigl(a_1(\alpha,\beta,\gamma),
       a_2(\alpha,\beta,\gamma),
       x_c(\alpha,\beta,\gamma)\bigr)
\]
defined for $0\le\alpha,\beta,\gamma<\delta$, where $x_c$ is the first coordinate of the continued cusp,
with $(a_1,a_2,x_c)=(342/343,368/343,r)$ at zero
coupling, such that \eqref{eq:KBoverlapCusp} has a
nondegenerate interior cusp near $(r,s,s)$. Its zero
eigenvalue is simple, its other two eigenvalues have
negative real parts, and $(a_1,a_2)$ generically unfold
the cusp.

\item Along this cusp continuation, for every sufficiently
small $\alpha,\beta,\gamma>0$, the system has
\emph{at least eight distinct locally exponentially stable
equilibria}: the origin, $E_{12}$, $E_{23}$, and at least
five strictly positive equilibria. These eight stable
equilibria are distinct from the interior cusp equilibrium.
Both proper siphon faces are transversely stable at
$E_{12}$ and $E_{23}$, respectively. No assertion of
exactly eight stable equilibria is made.
\rev{These eight are the product states of three weakly coupled bistable units: at the coefficients used, each coordinate equation has two stable roots, and weak coupling makes each of the $2^3$ sign patterns persist. The overlap of the two siphons gives the pair of boundary residents with invasion numbers distinct when $\alpha\ne\gamma$, while the count comes from the weak coupling; deeper in the cusp wedge the same argument gives $3^n$ (Remark~\ref{r:count}).}
\end{enumerate}
\eeT

\begin{proof}
For the siphon classification, on $x_1=0$ the first
production argument equals $\alpha x_2$, so a siphon
containing $1$ must contain $2$. Similarly, one containing
$3$ must contain $2$. On $x_2=0$ the second argument equals
$\beta x_1x_3$, so a siphon containing $2$ must contain
at least one of $1$ and $3$. This gives exactly the three
listed siphons. At an equilibrium, $x_2>0$ implies
$x_1,x_3>0$ by the first and third equations; whereas
$x_2=0$ implies $x_1x_3=0$ by the second equation.
Thus the asserted supports are exhaustive. In the strictly
positive interior, the four cross derivatives associated
with the edges $2\to1$, $1\to2$, $3\to2$, $2\to3$ are
strictly positive. Hence the Jacobian is irreducible Metzler.
The box is invariant because $\Phi_a(q)\ge0$ for $q\ge0$ and
$\Phi_a(q)<1$ for every finite $q\ge0$.

On $\Sigma_{12}$ the sole resident coordinate satisfies
$\dot x_3=-x_3+\Phi_a(x_3)$, independently of the couplings;
the same scalar equation holds for $x_1$ on $\Sigma_{23}$.
Differentiating the absent equations gives
\eqref{eq:KBoverlapJtra}. Since the linear removal matrix
is the identity, the nonnegative production blocks have
Perron eigenvalues
$a_1(1+\sqrt{\alpha\beta s_a})$ and
$a_1(1+\sqrt{\gamma\beta s_a})$, respectively.
The tangential eigenvalue is $\lambda_s(a)<0$.
The triangular face Jacobians then give the stability
equivalences. At the origin, the $x_1x_3$ term has no
linear part, and the Jacobian is triangular after reordering;
all its eigenvalues equal $a_1-1<0$.

For cusp continuation, regard $(x_2,x_3)$ as the
hyperbolic variables and $x_1$ as the scalar critical
variable. At $(a,\alpha,\beta,\gamma)=(a^0,0,0,0)$
and $(x_1,x_2,x_3)=(r,s,s)$, the derivatives of the last
two equilibrium equations with respect to $(x_2,x_3)$
form $\lambda_s(a^0)I_2$, which is invertible.
The implicit-function theorem therefore solves those
equilibrium equations locally as
\[
 (x_2,x_3)=Y(x_1;a,\alpha,\beta,\gamma).
\]
Substitution in the first equilibrium equation produces
the scalar Lyapunov--Schmidt equation
\[
 g(x_1;a,\alpha,\beta,\gamma)
 =-x_1+\Phi_a\bigl(x_1+\alpha Y_1(x_1;a,\alpha,\beta,\gamma)\bigr).
\]
At zero coupling, $g(x_1)=\Phi_a(x_1)-x_1$. By
Theorem~\ref{t:Qua}, at $(x_1,a)=(r,a^0)$,
\[
 g=g_{x_1}=g_{x_1x_1}=0,\qquad
 g_{x_1x_1x_1}=\frac{36}{49}\ne0,
\]
and
\[
 \det\frac{\partial(g,g_{x_1})}{\partial(a_1,a_2)}
 =r^2(1-r)^4>0.
\]
At zero coupling $J=\operatorname{Diag}(0,\lambda_s(a^0),\lambda_s(a^0))$ and $f_1$ is independent of
$(x_2,x_3)$. Both null vectors are $e_1$ and $h_{20}=0$, so these are the cusp conditions of
Proposition~\ref{p:sadCus}(2). Proposition~\ref{p:sadCus}(3) gives the stated smooth cusp
continuation in $(\alpha,\beta,\gamma)$, with the cubic, unfolding, and stable
complementary-eigenvalue conditions. Strict positivity of the cusp coordinates and of
all four odds coefficients persists by continuity.

Finally we prove the stable-equilibrium count without
assuming that the nonlinear equations can be solved
explicitly. Let $S\subset\{1,2,3\}$ and write $z^S(a)$ for
the vector having coordinate $s_a$ at each index in $S$
and zero elsewhere. At zero coupling the eight $z^S(a)$
are distinct, hyperbolically stable equilibria of the
three decoupled scalar equations. Define the increasing
map $T_{a,\alpha,\beta,\gamma}:[0,1]^3\to[0,1]^3$ by
\[
 T(x)=\bigl(\Phi_a(x_1+\alpha x_2),
             \Phi_a(x_2+\beta x_1x_3),
             \Phi_a(x_3+\gamma x_2)\bigr).
\]
For every $S$, nonnegativity of the couplings gives
$T(z^S(a))\ge z^S(a)$ componentwise.
Since $\Phi_a'(0+)<1$ and $\Phi_a'(s_a)<1$, one can choose a
single small $h>0$, uniformly for $a$ near $a^0$, such
that the eight boxes
\[
 I_S(a)=[z^S(a),z^S(a)+h\mathbf1]
\]
are pairwise disjoint, lie in $[0,1]^3$, and satisfy
$\Phi_a(t+h)<t+h$ for $t\in\{0,s_a\}$.
For sufficiently small nonnegative couplings,
$T(z^S(a)+h\mathbf1)\le z^S(a)+h\mathbf1$.
Monotonicity then makes every $I_S(a)$ invariant under
$T$. By reducing $h$ and the coupling bound if necessary,
$\Phi_a'$ is uniformly less than one at all arguments attained
in these boxes, and the row sums of $D T$ are uniformly
less than one, since the row sums of the Jacobian of the arguments tend to one as
the couplings vanish. Therefore $T$ is a contraction on
each box. Each box contains a unique equilibrium
$e^S$, and $J(e^S)=DT(e^S)-I$ is Hurwitz.

The three boxes indexed by $S=\varnothing,\{1\},\{3\}$
contain the exact stable equilibria $E_0,E_{23},E_{12}$.
For the remaining five subsets, the equilibrium is
strictly positive: if $2\in S$ then $x_2>0$ forces
$x_1,x_3>0$; if $S=\{1,3\}$ then $x_1,x_3>0$
forces $x_2>0$. Thus at least five stable equilibria
are interior. Their boxes are separated from the cusp,
whose first coordinate is close to $r$ instead of
to $0$ or $s$. Along the cusp continuation all parameters
remain in the neighborhoods just used. This completes
the proof.
\end{proof}

\paragraph{Method: certify cusp continuation and simultaneous stability.}
Choose a square augmented system $\mathcal G(X;\eta)=0$ for an equilibrium, its simple zero
eigenvalue, the vanishing quadratic coefficient, and two adjustable unfolding parameters. Here $X$
collects the unknown state, normalized nullvectors and adjustable parameters; $\eta$ denotes the
remaining couplings. The standard cusp formulas already used in the paper supply these equations.
State the normalization explicitly and check that the augmented Jacobian is invertible at the seed.
For a box $\mathcal X$, centre $X_0$, coupling box $E$, and a fixed numerical inverse approximation
$C$, enclose
\[
 \mathcal K(\mathcal X,E)
 =X_0-C\mathcal G(X_0;E)
  +\bigl(I-C D_X\mathcal G(\mathcal X;E)\bigr)(\mathcal X-X_0).
\]
Certify strict inclusion in $\operatorname{int}\mathcal X$ and a uniform norm bound
\[
 \sup_{X\in\mathcal X,\,\eta\in E}\|I-C D_X\mathcal G(X;\eta)\|<1.
\]
With $C$ nonsingular, these conditions give a unique augmented solution in $\mathcal X$ for each
fixed $\eta\in E$. Separately certify positivity of the required coefficients and states,
nonvanishing of the cubic coefficient and unfolding determinant, and negativity of the real parts
of all complementary eigenvalues. Enclose the continued stable equilibria in disjoint boxes with
validated existence and Hurwitz tests, and enclose the resident equilibria and their transverse
invasion inequalities. A successful cusp solve alone certifies none of these additional assertions.

\beC[Three-population odds networks admit nine stable equilibria near the interior cusp][nine stable equilibria][Nin]
Fix any sufficiently small positive coupling triple
$(\alpha,\beta,\gamma)$ covered by Theorem~\ref{t:KBoverlapCusp}, and let
$a^c=(a_1^c,a_2^c)$ be its continued cusp coefficients, with
$a_3=128/343$ and $a_4=1024/343$ fixed.
There exists a nonempty open set of admissible coefficient pairs
$(a_1,a_2)$ whose closure contains $a^c$, on which the network has
at least nine distinct locally exponentially stable equilibria:
the origin, two equilibria on the proper siphon faces, and at least
six strictly positive equilibria.
Two additional strictly positive equilibria near the former cusp
each have exactly one unstable eigenvalue.
The minimal siphons remain $\Sigma_{12}=\{1,2\}$ and
$\Sigma_{23}=\{2,3\}$, and the two continued stable boundary equilibria
remain transversely stable.
\eeC
\begin{proof}
At $a=a^c$, Theorem~\ref{t:KBoverlapCusp} supplies eight distinct hyperbolically
stable equilibria, separated from the interior cusp equilibrium.
The implicit-function theorem continues these equilibria to a common
neighbourhood of $a^c$, preserving their distinctness and stability.
The origin remains an equilibrium. The two stable boundary equilibria
remain on their respective invariant siphon faces, and the other five
remain strictly positive.

Write the vector field as
$f=(f_1,f_2,f_3)^\top$ and put $H=(f_2,f_3)^\top$.
As in the proof of Theorem~\ref{t:KBoverlapCusp}, the last two equilibrium
equations can be solved near the cusp as
\[
 (x_2,x_3)=Y(x_1;a),
 \qquad
 g(x_1;a)=f_1(x_1,Y(x_1;a);a).
\]
At the cusp,
\[
 g=g_{x_1}=g_{x_1x_1}=0,\qquad
 g_{x_1x_1x_1}>0,\qquad
 \det\frac{\partial(g,g_{x_1})}{\partial(a_1,a_2)}\ne0.
\]
The positivity of the cubic follows by continuity from its reference
value $36/49$. The generic scalar cusp unfolding therefore has an open
parameter region, accumulating at $a^c$, with three simple nearby
equilibria. In their increasing $x_1$ order, the derivatives
$g_{x_1}$ are positive, negative, and positive, respectively.

To transfer these signs to the full system, let
\[
 J_{\mathrm{res}}
   =D_{(x_2,x_3)}H
\]
at a nearby equilibrium. This matrix remains Hurwitz, since it is close
to the strictly stable resident matrix used in the cusp construction.
The Schur-complement identity gives
\[
 \det D_xf=\det J_{\mathrm{res}}\,g_{x_1}.
\]
The two eigenvalues complementary to the critical eigenvalue remain
in the open left half-plane. Their product is positive, as is
$\det J_{\mathrm{res}}$. Hence the remaining real eigenvalue has the
same sign as $g_{x_1}$. The middle equilibrium is locally exponentially
stable, whereas each outer equilibrium has exactly one unstable
eigenvalue.

All three equilibria are strictly positive and lie in a neighbourhood
disjoint from the eight continued stable equilibria. The middle one
therefore raises the stable count to at least nine, including at least
six strictly positive equilibria. Shrinking the parameter region
preserves positivity of the odds coefficients and the strict
inequalities $R_{\Sigma_{12}}<1$ and $R_{\Sigma_{23}}<1$ at the two
continued boundary equilibria. The fixed positive couplings and
unchanged production supports preserve the minimal siphons.
\end{proof}

\beR[The count nine, and the count $3^n$ deeper in the cusp wedge][what is counted][count]
{\color{revcolor}
Two different counts must be kept apart. The eight equilibria of
Theorem~\ref{t:KBoverlapCusp} are the product states of three weakly coupled
bistable units: at the couplings of Table~\ref{tab:KBoverlapNum} each
coordinate equation has two stable roots, and weak-coupling persistence makes
every one of the $2^3$ sign patterns an equilibrium of the network. The ninth
equilibrium of Corollary~\ref{c:Nin} is the one produced by the interior cusp,
on the lower outer branch.

Nine is therefore not a maximum. All three coordinates carry the same response
$\Phi_a$, so moving $(a_1,a_2)$ further into the cusp wedge gives each
coordinate equation three stable roots, and the same persistence argument then
gives $3^3=27$ stable equilibria, and $3^n$ for $n$ coordinates coupled in this
way; the siphon faces carry further stable residents. Nine is the count
located by the continuation reported here, at the coefficients of
Table~\ref{tab:KBnineNum}, not the largest count the construction allows.}
\eeR

The weak-coupling continuation is an application of the implicit-function theorem; compare \cite{MacKaySepulchre1995,ChengLinShih2006}.

\beL[Weak coupling preserves products of stable scalar equilibria]
[product continuation][ProdCount]
Let $\Phi$ be $C^1$ on a neighbourhood of $[0,1]$, with values in
$[0,1]$, and suppose the scalar equation $\dot u=-u+\Phi(u)$
has $k$ distinct hyperbolic locally attracting equilibria
$u^{(1)},\ldots,u^{(k)}$. For nonnegative $C^1$ functions
$Q_i:[0,1]^n\to\mathbb R_+$, consider
\[
\dot x_i=-x_i+\Phi(x_i+\varepsilon Q_i(x)),
\qquad i=1,\ldots,n.
\]
There exists $\varepsilon_0>0$ such that, for every
$0\le\varepsilon<\varepsilon_0$, this system has at least
$k^n$ distinct hyperbolic locally attracting equilibria.
In particular, within an open cusp-unfolding region of the quartic
response in Theorem~\ref{t:Qua}, where the scalar equation has
three hyperbolic attracting roots, the three-coordinate system
\eqref{eq:KBoverlapCusp} has at least $27$ stable equilibria for
all sufficiently small positive $\alpha,\beta,\gamma$.
No particular numerical coupling window is asserted.
\eeL
\begin{proof}
At zero coupling the $k^n$ products of the scalar stable roots
are distinct equilibria. Their Jacobians are diagonal, with strictly
negative diagonal entries. The implicit-function theorem continues
each product uniquely under small couplings, and continuity preserves
its Hurwitz Jacobian. Finiteness gives a common positive coupling
neighbourhood. For the quartic cusp unfolding with positive cubic
coefficient, the middle of the three nearby cusp roots is stable;
together with the stable origin and the continuation of the stable
root near $1/2$, this gives three scalar stable roots.
\end{proof}

\paragraph{A larger numerical equilibrium count.}
For the same three-population model at
$(a_1,a_2)=(0.997213369,1.071847276)$ and
$(\alpha,\beta,\gamma)=(2\cdot10^{-7},10^{-6},5\cdot10^{-8})$,
Newton solves initialized at all $5^3$ products of the five scalar
fixed points recover $125$ distinct equilibria, of which $27$
have Hurwitz Jacobians. The largest reported residual is
$5.6\cdot10^{-17}$ and the minimum separation between recovered
states is $0.02639$. The calculation is reproduced by
\texttt{verify\_125.py}. These are numerical counts; they are not
interval-certified or a proof that no additional equilibria exist.

\beP[A composite energy certifies an explicit subset of the basin of a boundary equilibrium][][ComBas]
For \eqref{eq:KBoverlapCusp} with $\alpha,\beta,\gamma>0$ and $a_1>0$, let $\mathcal E_a$ be the
scalar energy following Theorem~\ref{t:Qua}. Take $\eta>0$, and put
$\ell=(\sqrt\beta,\sqrt\alpha)$, $\bar z=\eta/\min\{\sqrt\alpha,\sqrt\beta\}$,
$\bar q=\bar z(1+\max\{\alpha,\beta\})$ and $\rho=\sup_{0<q\le\bar q}\Phi_a(q)/q$. Assume
\[
 c:=1-\rho\bigl(1+\sqrt{\alpha\beta}\bigr)>0,
\]
and that $0<y_-<s_a<y_+<1$ satisfy $\Phi_a(q)>q$ on $[y_-,s_a)$, $\Phi_a(q)<q$ on $(s_a,y_+]$,
$\Phi_a(y_++\gamma\bar z)\le y_+$ and $\Phi_a'(s_a)<1$. Put
\[
 M_U=\max_{y\in[y_-,y_+]}|\Phi_a^{-1}(y)-y|,\qquad
 M_\Phi=\max_{q\in[y_-,\,y_++\gamma\bar z]}\Phi_a'(q).
\]
Then
\[
 K=\{x:\ x_1,x_2\ge0,\ \sqrt\beta\,x_1+\sqrt\alpha\,x_2\le\eta,\ y_-\le x_3\le y_+\}
\]
is forward invariant, and for $\kappa>\gamma M_UM_\Phi/(c\sqrt\alpha)$,
\[
 H(x)=\mathcal E_a(x_3)-\mathcal E_a(s_a)+\kappa\bigl(\sqrt\beta\,x_1+\sqrt\alpha\,x_2\bigr)
\]
is a strict Lyapunov function on $K$. Thus $E_{12}=(0,0,s_a)$ is globally asymptotically stable
relative to $K$, so $K$ is an explicit subset of its basin of attraction. Exchanging $x_1$ with $x_3$ and $\alpha$ with $\gamma$ gives the same statement for
$E_{23}$.
\eeP
\begin{proof}
Put $z=(x_1,x_2)$. On $K$, $x_1,x_2\le\bar z$ and $x_3\le1$, so the arguments of the first two equations lie
in $[0,\bar q]$. Hence $\dot x_1\le-x_1+\rho(x_1+\alpha x_2)$ and $\dot x_2\le-x_2+\rho(x_2+\beta x_1)$, and
therefore $\ell^\top\dot z\le-c\,\ell^\top z$. The faces $x_1=0$, $x_2=0$ point inward, and
$\ell^\top z$ decreases at the level $\eta$. On $x_3=y_-$, $\dot x_3\ge-y_-+\Phi_a(y_-)>0$; on $x_3=y_+$,
$\dot x_3\le-y_++\Phi_a(y_++\gamma\bar z)\le0$. So $K$ is forward invariant.

Let $U(y)=\mathcal E_a(y)-\mathcal E_a(s_a)$, so $U'(y)=\Phi_a^{-1}(y)-y$. The sign conditions and
monotonicity of $\Phi_a$ give $U'<0$ on $[y_-,s_a)$ and $U'>0$ on $(s_a,y_+]$; thus $U\ge0$ vanishes only at
$s_a$. Next, $U'(y)(-y+\Phi_a(y))=-(\Phi_a^{-1}(y)-y)(y-\Phi_a(y))$, and the quotient of
$(\Phi_a^{-1}(y)-y)(y-\Phi_a(y))$ by $(y-s_a)^2$ is continuous and positive on $[y_-,y_+]$, with limit
$(1-\Phi_a'(s_a))^2/\Phi_a'(s_a)$ at $s_a$. This gives the resident dissipation estimate used below. Finally,
$U'(x_3)(\Phi_a(x_3+\gamma x_2)-\Phi_a(x_3))\le M_UM_\Phi\gamma x_2\le b\,\ell^\top z$, with
$b=\gamma M_UM_\Phi/\sqrt\alpha$. The positive resident quotient
established above gives $U'(y)(-y+\Phi_a(y))\le-a(y-s_a)^2$
for some $a>0$. Therefore
$\dot H\le-a(x_3-s_a)^2-(\kappa c-b)\ell^\top z<0$
away from $E_{12}$.
\end{proof}

The hypothesis $c>0$ is stronger than $R_{\Sigma_{12}}<1$, because it holds on all of $K$: it uses
$\rho\ge a_1$ and the bound $x_3\le1$ instead of the value $s_a$. At the couplings of
Table~\ref{tab:KBoverlapNum}, with the reference coefficients $(a_1,a_2)=(342/343,368/343)$,
$s_a=1/2$, $\bar z=1/10$, $y_-=1/5$ and $y_+=7/10$, one has $\rho\approx0.99998$ and
$c\approx1.6\cdot10^{-5}$; this numerical check, with the derivative bound tested at sampled points,
is illustrative.

\beP[A Lyapunov--Volterra function certifies attraction to a boundary resident]
[][LVboundary]
Under the hypotheses of Proposition~\ref{p:ComBas}, retain its invariant set $K$ and constants $c>0$, $M_\Phi$, $y_-$, $y_+$ and $s_a$. Define
\[
 U_V(y)=y-s_a-s_a\log\frac{y}{s_a},\qquad
 M_V=\max_{y\in[y_-,y_+]}\left|1-\frac{s_a}{y}\right|,\qquad
 b_V=\frac{\gamma M_\Phi M_V}{\sqrt{\alpha}}.
\]
For every $\kappa>b_V/c$, the function
\[
 H_V(x)=x_3-s_a-s_a\log\frac{x_3}{s_a}+\kappa\bigl(\sqrt{\beta}\,x_1+\sqrt{\alpha}\,x_2\bigr)
\]
is a strict Lyapunov function on $K$ for $E_{12}=(0,0,s_a)$. More precisely, put
\[
 a_*=\min_{y\in[y_-,y_+]}\frac{y-\Phi_a(y)}{y-s_a}>0,
\]
where the quotient at $y=s_a$ is defined as $1-\Phi_a'(s_a)$. Then
\[
 \dot H_V\le-\frac{a_*}{y_+}(x_3-s_a)^2-(\kappa c-b_V)\bigl(\sqrt{\beta}\,x_1+\sqrt{\alpha}\,x_2\bigr),
\]
and
\[
 H_V(x(t))\le e^{-\omega_Vt}H_V(x(0)),\qquad
 \omega_V=\min\left\{\frac{2a_*y_-^2}{s_ay_+},\ c-\frac{b_V}{\kappa}\right\}>0.
\]
Consequently every trajectory starting in $K$ converges to $E_{12}$.
\eeP
\begin{proof}
Write $y=x_3$ and $Z=\sqrt{\beta}\,x_1+\sqrt{\alpha}\,x_2$. The proof of Proposition~\ref{p:ComBas} gives $\dot Z\le-cZ$. For the resident equation,
\[
 U_V'(y)\bigl(-y+\Phi_a(y)\bigr)=-\frac{y-s_a}{y}\bigl(y-\Phi_a(y)\bigr)\le-\frac{a_*}{y_+}(y-s_a)^2.
\]
The additional production caused by $x_2$ satisfies
\[
 U_V'(y)\bigl[\Phi_a(y+\gamma x_2)-\Phi_a(y)\bigr]\le M_VM_\Phi\gamma x_2\le b_VZ.
\]
Adding $\kappa\dot Z$ proves the derivative inequality. The sign assumptions of Proposition~\ref{p:ComBas} and $\Phi_a'(s_a)<1$ make the quotient defining $a_*$ continuous and strictly positive.

Moreover, $U_V(s_a)=U_V'(s_a)=0$ and $U_V''(y)=s_a/y^2$, so
\[
 \frac{s_a}{2y_+^2}(y-s_a)^2\le U_V(y)\le\frac{s_a}{2y_-^2}(y-s_a)^2.
\]
These bounds give $\dot H_V\le-\omega_VH_V$ and the asserted convergence.
\end{proof}

This is the classical weighted Lyapunov--Volterra form about $E_{12}$. The positive resident coordinate contributes $x_3-s_a-s_a\log(x_3/s_a)$, and the coordinates whose equilibrium values are zero contribute linear terms, using $\lim_{p\downarrow0}(x-p-p\log(x/p))=x$. This gives a boundary-equilibrium Lyapunov function for the threshold network. At the values above, $M_V=3/2$, $b_V\approx1.69\cdot10^{-4}$ and $b_V/c\approx10.4$ (illustrative).

For Theorem~\ref{t:Net} at $a=a^0$, suppose $0<\bar x<1/8$ and
$c_2:=1-(1+\gamma)\sup_{0<q\le\bar x(1+\gamma)}\Phi_{a^0}(q)/q>0$. Then for every $\varepsilon\ge0$ every
solution starting in $[0,\bar x]\times(0,1]$ converges to $(0,y^*)$. Indeed, $y\le1$ gives
$\dot x\le-c_2x$ on $[0,\bar x]$, so $x$ decays exponentially. Monotonicity of the response gives
$\dot y\ge-y+\Phi_{a^0}(10y)$, so comparison with the scalar resident equation, whose only positive
equilibrium is $y^*$, gives $\liminf y\ge y^*$. Since $y\le1$ and $x\to0$, the difference
$\Phi_{a^0}(y(10+\varepsilon x))-\Phi_{a^0}(10y)$ tends to zero uniformly, and the scalar field is negative
above $y^*$, so $\limsup y\le y^*$. On compact subsets of $[0,\bar x]\times(0,1)$ that are forward
invariant, the resident storage $\int\Psi^{-1}-y^2/2$ about $y^*$, with $\Psi(y)=\Phi_{a^0}(10y)$, plus a
linear term $\kappa x$ gives a Lyapunov function by the same composite derivative estimate used above.

\beR[Two boundary invasion crossings]
[double boundary transcritical relay][KBrelay]
Fix $a$ with $a_1<1$ and $\alpha,\gamma>0$, and vary only $\beta$. Species $2$ belongs to both $\Sigma_{12}$ and
$\Sigma_{23}$, so its production, which contains $\beta x_1x_3$, gives a row of both transverse matrices. Since
$x_2=0$ at $E_{12}$ and $E_{23}$, the residents do not depend on $\beta$. Hence both invasion numbers increase
with $\beta$ at fixed residents, and
\[
 R_{\Sigma_{12}}=1\iff\beta=\beta_{12}:=\frac{(1-a_1)^2}{a_1^2\alpha s_a},\qquad
 R_{\Sigma_{23}}=1\iff\beta=\beta_{23}:=\frac{(1-a_1)^2}{a_1^2\gamma s_a}.
\]
Since the tangential eigenvalue $\lambda_s(a)$ is negative, each resident is stable below its value and unstable
above it. The two stability crossings concern different residents on overlapping siphon faces. Their transcritical bifurcation direction is proved in Theorem~\ref{t:OvlRelay}.

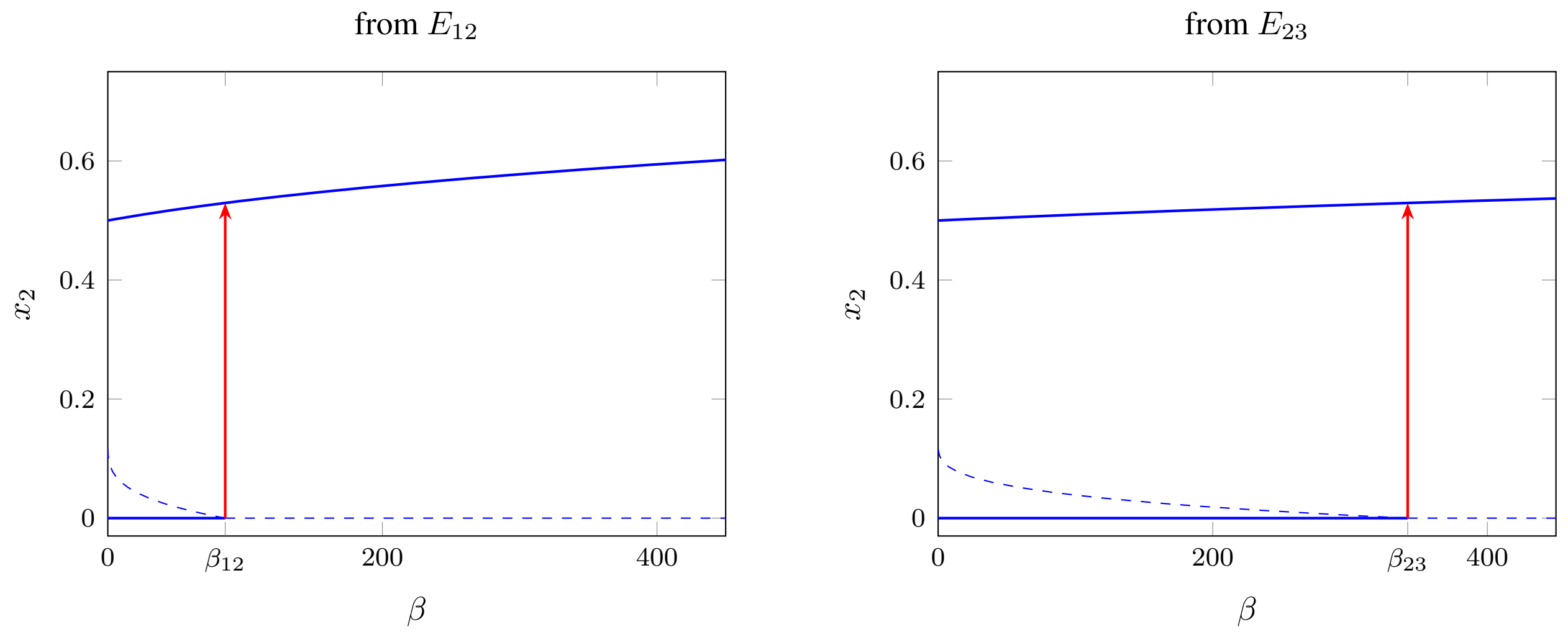
\begin{figure}[H]
\centering
\begin{minipage}{.48\linewidth}
\centering
\begin{tikzpicture}
\begin{axis}[width=\linewidth,height=.8\linewidth,xlabel={$\beta$},ylabel={$x_2$},xmin=0,xmax=450,ymin=-.03,ymax=.75,
 xtick={0,200,400},tick label style={font=\scriptsize},label style={font=\small},title style={font=\small},extra x ticks={85.5},extra x tick labels={$\beta_{12}$},title={from $E_{12}$}]
\addplot[blue,thick] coordinates {(0,0) (85.50,0)};
\addplot[blue,dashed] coordinates {(85.50,0) (450,0)};
\addplot[blue,dashed] coordinates {(85.50,0.0000) (85.03,0.0002) (84.79,0.0003) (84.43,0.0005) (83.88,0.0007) (83.04,0.0011) (81.79,0.0016) (79.91,0.0025) (77.12,0.0037) (73.01,0.0057) (67.07,0.0087) (58.74,0.0132) (47.56,0.0200) (40.44,0.0250) (29.80,0.0338) (21.34,0.0426) (14.76,0.0515) (9.77,0.0603) (6.12,0.0691) (3.55,0.0779) (1.85,0.0867) (0.82,0.0955) (0.27,0.1044) (0.05,0.1132) (0.00,0.1220) (0.00,0.1250)};
\addplot[blue,thick] coordinates {(0.00,0.5000) (4.85,0.5020) (24.61,0.5096) (46.37,0.5173) (70.22,0.5249) (96.25,0.5326) (124.58,0.5402) (155.28,0.5478) (188.48,0.5555) (224.28,0.5631) (262.80,0.5708) (304.15,0.5784) (348.45,0.5861) (395.85,0.5937) (446.47,0.6013) (455.00,0.6026)};
\draw[red,thick,-{Stealth[length=5pt]}] (axis cs:85.50,0) -- (axis cs:85.50,0.5295);
\end{axis}
\end{tikzpicture}
\end{minipage}\hfill
\begin{minipage}{.48\linewidth}
\centering
\begin{tikzpicture}
\begin{axis}[width=\linewidth,height=.8\linewidth,xlabel={$\beta$},ylabel={$x_2$},xmin=0,xmax=450,ymin=-.03,ymax=.75,
 xtick={0,200,400},tick label style={font=\scriptsize},label style={font=\small},title style={font=\small},extra x ticks={342.0},extra x tick labels={$\beta_{23}$},title={from $E_{23}$}]
\addplot[blue,thick] coordinates {(0,0) (341.99,0)};
\addplot[blue,dashed] coordinates {(341.99,0) (450,0)};
\addplot[blue,dashed] coordinates {(341.99,0.0000) (340.13,0.0002) (339.17,0.0003) (337.71,0.0005) (335.51,0.0007) (332.18,0.0011) (327.16,0.0016) (319.65,0.0025) (308.47,0.0037) (292.03,0.0057) (268.30,0.0087) (234.97,0.0132) (190.26,0.0200) (161.75,0.0250) (119.20,0.0338) (85.38,0.0426) (59.06,0.0515) (39.10,0.0603) (24.46,0.0691) (14.18,0.0779) (7.39,0.0867) (3.26,0.0955) (1.09,0.1044) (0.20,0.1132) (0.00,0.1220) (0.00,0.1250)};
\addplot[blue,thick] coordinates {(0.00,0.5000) (19.42,0.5020) (98.51,0.5096) (185.60,0.5173) (281.08,0.5249) (385.30,0.5326) (455.00,0.5373)};
\draw[red,thick,-{Stealth[length=5pt]}] (axis cs:341.99,0) -- (axis cs:341.99,0.5295);
\end{axis}
\end{tikzpicture}
\end{minipage}
\caption{Numerical continuation of the two boundary invasion crossings of Remark~\ref{r:KBrelay} for $\dot x_1=-x_1+\Phi_a(x_1+\alpha x_2)$, $\dot x_2=-x_2+\Phi_a(x_2+\beta x_1x_3)$, $\dot x_3=-x_3+\Phi_a(x_3+\gamma x_2)$: equilibria, with the backward branches given by $\beta=(\Phi_a^{-1}(x_2)-x_2)/(x_1x_3)$, for the reference coefficients and the couplings
$\alpha,\gamma$ of Table~\ref{tab:KBoverlapNum}: the shared species $x_2$ against $\beta$. Solid curves are stable and dashed curves
unstable. Each resident ($x_2=0$) loses stability at its own value, $\beta_{12}\approx85.5$ or $\beta_{23}\approx342$; the unstable
positive branch leaves it backward and ends near $x_2=1/8$ as $\beta\to0$. Past its value the state jumps to the upper
branch, where it remains when $\beta$ is decreased.}
\label{fig:relay}
\end{figure}

For $a_2>a_1^2$, the expansion
$\Phi_a(q)=a_1q+(a_2-a_1^2)q^2+O(q^3)$
provides the quadratic coefficient for the local bifurcation calculation.
At each crossing the transverse block has a simple zero eigenvalue with
positive Perron vectors $w,v$. The numerical calculation below evaluates
$v^\top(\partial_\beta J^\perp)w$ and
$\frac12v^\top D^2f[w,w]$. The positive quadratic coefficients are proved in Theorem~\ref{t:OvlRelay}. Along the numerically continued branch
from $E_{12}$, the second equation gives
\[
 \beta=\frac{\Phi_a^{-1}(x_2)-x_2}{x_1x_3},
\]
with $x_1,x_3$ the roots of the first and third equations near $0$ and $s_a$. Figure~\ref{fig:relay} is computed
from this formula.

\emph{Numerical evidence.} The following was computed for the reference coefficients $a=(342,368,128,1024)/343$, so
$s_a=1/2$ and $a_2-a_1^2\approx0.079$, with $\alpha=2\cdot10^{-7}$, $\gamma=5\cdot10^{-8}$, which give $\beta_{12}\approx85.5$ and
$\beta_{23}\approx342$. Trajectories were integrated to $t=5\cdot10^7$ from $x_2(0)=10^{-4}$.
\begin{enumerate}
\item At $E_{12}$, with the critical Perron vectors normalized to have second coordinate one,
$v^\top(\partial_\beta J)w\approx3.41\cdot10^{-5}$ and $\frac12v^\top D^2f[w,w]\approx0.079$. The predicted branch agrees with Newton
solves: at $\beta=0.999\,\beta_{12}$ it predicts $x_2\approx3.682\cdot10^{-5}$, and the solve gives $3.684\cdot10^{-5}$.
\item The largest real part of the Jacobian spectrum lies in $[3.3\cdot10^{-6},1.04\cdot10^{-3}]$ along the sampled backward branches,
which are therefore unstable, and in $[-2.9\cdot10^{-3},-2.6\cdot10^{-3}]$ along the upper branches, which are stable.
\item Below its value each resident attracts the perturbed state: the trajectories return to $E_{12}$ at $0.8\,\beta_{12}$ and to $E_{23}$
at $0.8\,\beta_{23}$. Past its value the state jumps to an upper branch: the end states are $(3.66\cdot10^{-5},0.5343,0.5000)$ at
$1.2\,\beta_{12}$ from $E_{12}$, and $(0.5000,0.5343,9.1\cdot10^{-6})$ at $1.2\,\beta_{23}$ from $E_{23}$.
\item Between the two values one resident has been lost while the other survives: at $\beta=1.2\,\beta_{12}$ the trajectory from $E_{23}$ returns to
$E_{23}$.
\item The jump is not undone by decreasing $\beta$, down to the smallest value tested: from the upper state reached at
$1.2\,\beta_{12}$, the trajectory stays on the upper branch at $0.5\,\beta_{12}$, with $x_2\approx0.516$, and at $0.05\,\beta_{12}$, with $x_2\approx0.502$.
This is consistent with the backward branch reaching $\beta=0$ near $x_2=1/8$, so that no fold closes a hysteresis loop for $\beta>0$.
\end{enumerate}

The relay requires non-disjoint siphons. If $\Sigma$ and $\Sigma'$ are disjoint, the production of a species of $\Sigma$
gives no row of $F_{\Sigma'}$, so changing it can affect $J_{\Sigma'}^\perp$ only by moving the resident.
Folds and cusps do not distinguish the two cases: Theorem~\ref{t:Net} has them with the disjoint singleton siphons
$\{x\}$ and $\{y\}$. The values above concern the boundary residents, and do not extend the small-coupling
range of Theorem~\ref{t:KBoverlapCusp}: they lie far outside it, so the relay is not shown to occur next to the interior cusp.
\eeR

\beT[Two nonzero boundary equilibria of a three-dimensional polynomial production threshold network undergo backward transcritical bifurcations]
[backward boundary relay][OvlRelay]
{\color{revcolor}
For the three-dimensional system~\eqref{eq:KBoverlapCusp} of Theorem~\ref{t:KBoverlapCusp}, take an odds-polynomial response $\Phi_a$ as in Definition~\ref{d:Odd} with $0<a_1<1$, $a_2>a_1^2$, and fix $\alpha,\gamma>0$. Suppose $s_a>0$ is a scalar equilibrium with $\Phi_a(s_a)=s_a$ and $0<\Phi_a'(s_a)<1$. Keep the response and $\alpha,\gamma$ fixed and vary $\beta$. At $E_{12}=(0,0,s_a)$ and $E_{23}=(s_a,0,0)$ the two invasion
thresholds $\beta_{12},\beta_{23}$ of Remark~\ref{r:KBrelay} are
nondegenerate transcritical bifurcations. At either threshold a local
branch of strictly positive equilibria exists for $\beta$ below its
threshold, and each equilibrium on this branch has exactly one
unstable eigenvalue. The corresponding boundary resident is locally
exponentially stable below its threshold and has one unstable
eigenvalue above it.}
\eeT
\begin{proof}
{\color{revcolor}
Put $\delta=1-a_1>0$ and $\kappa=a_2-a_1^2>0$. The scalar response has
the expansion
\[
 \Phi_a(q)=a_1q+\kappa q^2+O(q^3)
\]
about zero, using its smooth rational extension. The residents do not
move with $\beta$. At either threshold the transverse matrix has a
simple zero eigenvalue, its other eigenvalue is $-2\delta<0$, and the
resident eigenvalue is $\lambda_s=\Phi_a'(s_a)-1<0$. Choose a right
zero eigenvector $r\gg0$ and a left zero eigenvector $\pi\ge0$, with
$\pi^\top r=1$; the resident component of $\pi$ is zero.

At $E_{12}$, write $r=(r_1,r_2,r_3)$ and
$\pi=(\pi_1,\pi_2,0)$, with $\pi_1,\pi_2>0$. The resident component
of the right eigenvector is
\[
 r_3=\frac{\gamma\Phi_a'(s_a)}{1-\Phi_a'(s_a)}r_2>0.
\]
The crossing and quadratic coefficients in the one-dimensional
center-manifold equation are
\[
\begin{aligned}
 A_{12}&=\pi^\top(\partial_\beta J)r
       =\pi_2a_1s_ar_1>0,\\
 B_{12}&=\frac12\pi^\top D_x^2f(E_{12})[r,r]\\
 &=\kappa\left[\pi_1(r_1+\alpha r_2)^2
     +\pi_2(r_2+\beta_{12}s_ar_1)^2\right]
     +\pi_2a_1\beta_{12}r_1r_3>0.
\end{aligned}
\]
At $E_{23}$ the resident is the first component, so
$\pi=(0,\pi_2,\pi_3)$ with $\pi_2,\pi_3>0$ and
\[
 r_1=\frac{\alpha\Phi_a'(s_a)}{1-\Phi_a'(s_a)}r_2>0.
\]
This gives
\[
\begin{aligned}
 A_{23}&=\pi^\top(\partial_\beta J)r
       =\pi_2a_1s_ar_3>0,\\
 B_{23}&=\frac12\pi^\top D_x^2f(E_{23})[r,r]\\
 &=\kappa\left[\pi_2(r_2+\beta_{23}s_ar_3)^2
     +\pi_3(r_3+\gamma r_2)^2\right]
     +\pi_2a_1\beta_{23}r_1r_3>0.
\end{aligned}
\]
Because the resident branch is invariant for every $\beta$, the
center-manifold equation at either threshold has the form
\[
 \dot u=A_j(\beta-\beta_j)u+B_ju^2
  +O\bigl(u^3+|\beta-\beta_j|u^2
                 +(\beta-\beta_j)^2|u|\bigr),
 \qquad j\in\{12,23\}.
\]
The nonresident branch satisfies
$u=-A_j(\beta-\beta_j)/B_j+O((\beta-\beta_j)^2)$.
Since $r\gg0$, this branch is strictly positive when
$\beta<\beta_j$ is sufficiently close to $\beta_j$.
Its center eigenvalue is $-A_j(\beta-\beta_j)+o(|\beta-\beta_j|)>0$,
while the other two eigenvalues remain in the open left half-plane.
The resident branch has the opposite center eigenvalue.}
\end{proof}

\beR[The siphon lattice changes at zero coupling]
The proof starts from the decoupled system $\alpha=\beta=\gamma=0$, at which
every singleton $\{i\}$ is a siphon and every coordinate face is invariant. For
$\alpha,\beta,\gamma>0$ the minimal siphons are instead $\Sigma_{12}$ and
$\Sigma_{23}$. The implicit-function arguments only need the decoupled point as
a base for continuation, so this causes no gap; but the siphon lattice is not
constant across that point. Theorem~\ref{t:KBoverlapCusp} is a statement about
the positive-coupling family, which retains the overlapping-siphon lattice of
Example~\ref{e:KBoverlap}; it is not a continuation through zero coupling with an
unchanged lattice. Likewise, $R_{\Sigma_{12}}$ and $R_{\Sigma_{23}}$ depend on the
higher coefficients $a_2,a_3,a_4$ through $s_a$, although these coefficients do
not appear explicitly in \eqref{eq:oddsTransverse}.
\eeR

\paragraph{A numerical instance.}
Table~\ref{tab:KBoverlapNum} reports one coupling triple,
$(\alpha,\beta,\gamma)=(2\cdot10^{-7},10^{-6},5\cdot10^{-8})$. The cusp was
located by Newton's method on the seven equations $f=0$, $Jv=0$ ($v_1=1$),
$w^\top D^2f(v,v)=0$ with $w$ a left null vector, in 60-digit arithmetic.
The cubic coefficient and the unfolding determinant are those of the
Lyapunov--Schmidt function $g$ of the proof. They were computed by central
differences and agree, to the digits shown, with implicit differentiation of
the resident equations. The eight equilibria were found by Newton's method from
the points $z^S(a)$. The computation is illustrative, with no interval validation;
the supplementary script \texttt{verify\_cusp\_tables.py} independently computes the cusp coordinates, continued equilibria, and their spectra.
\begin{table}[H]
\centering\small
\begin{tabular}{ll}
\hline
Quantity & Value\\\hline
continued coefficients $(a_1,a_2)$ & $(0.997082757,\ 1.072892269)$\\
cusp $(x_1,x_2,x_3)$ & $(0.124986385,\ 0.500005319,\ 0.500004404)$\\
complementary eigenvalues & $-0.0393616,\ -0.0393609$\\
$g_{x_1x_1x_1}$ (reference $36/49\approx0.734694$) & $0.734463$\\
$\det\partial(g,g_{x_1})/\partial(a_1,a_2)$ (reference $0.00915909$) & $0.00915761$\\
resident level $s_a$ & $0.500003793$\\
$R_{\Sigma_{12}}$, $R_{\Sigma_{23}}$ & $0.997083072,\ 0.997082914$ \rev{(difference $1.58\cdot10^{-7}$)}\\\hline
stable equilibrium, box $S$ & spectral abscissa\\\hline
$\varnothing$: $(0,0,0)$ & $-0.002917$\\
$\{3\}$: $(0,0,0.500003793)=E_{12}$ & $-0.002917$\\
$\{1\}$: $(0.500003793,0,0)=E_{23}$ & $-0.002917$\\
$\{2\}$: $(3.42\cdot10^{-5},\ 0.500003793,\ 8.55\cdot10^{-6})$ & $-0.002912$\\
$\{2,3\}$: $(3.42\cdot10^{-5},\ 0.500003794,\ 0.500004404)$ & $-0.002912$\\
$\{1,3\}$: $(0.500003794,\ 8.56\cdot10^{-5},\ 0.500003793)$ & $-0.002904$\\
$\{1,2\}$: $(0.500006234,\ 0.500003793,\ 8.55\cdot10^{-6})$ & $-0.002916$\\
$\{1,2,3\}$: $(0.500006234,\ 0.500009895,\ 0.500004404)$ & $-0.039361$\\\hline
\end{tabular}
\caption{Theorem~\ref{t:KBoverlapCusp} at
$(\alpha,\beta,\gamma)=(2\cdot10^{-7},10^{-6},5\cdot10^{-8})$, with $a_3,a_4$
fixed. Both invasion numbers are below one and differ, since $\alpha\ne\gamma$.}
\label{tab:KBoverlapNum}
\end{table}

\paragraph{A nearby instance with nine stable equilibria.}
Corollary~\ref{c:Nin} concerns coefficient pairs near, but not at, the cusp.
Keeping the couplings of Table~\ref{tab:KBoverlapNum}, the same script moves
$(a_1,a_2)$ from the cusp values by the solution $\delta a$ of
$\partial_a(g,g_{x_1})\,\delta a=(0,-10^{-6})$, which makes $g_{x_1}$ negative at
the cusp point. Newton's method from the eight points $z^S(b)$ and from three
points near the former cusp gives eleven distinct equilibria \rev{among those found,
the search not being exhaustive (Remark~\ref{r:count})}, with residuals
below $10^{-52}$ and minimal pairwise distance $2.85\cdot10^{-3}$
(Table~\ref{tab:KBnineNum}). Nine of them have negative spectral abscissa. The
middle equilibrium near the former cusp is stable, with spectral abscissa
$-1.0\cdot10^{-6}$, of the order of the shift in $g_{x_1}$; the two outer ones
have one positive eigenvalue each. This is again illustrative, not
interval-validated.
\begin{table}[H]
\centering\small
\begin{tabular}{ll}
\hline
Quantity & Value\\\hline
coefficients $(a_1,a_2)$ & $(0.997084063,\ 1.072881819)$\\
resident level $s_b$ & $0.499991347$\\
$R_{\Sigma_{12}}$, $R_{\Sigma_{23}}$ & $0.997084378,\ 0.997084221$\\\hline
equilibrium & spectral abscissa\\\hline
the eight continued equilibria of Table~\ref{tab:KBoverlapNum} & between $-0.03935$ and $-0.002902$\\
$(0.122134360,\ 0.499992838,\ 0.499991957)$ & $+1.97\cdot10^{-6}$\\
$(0.124986385,\ 0.499992872,\ 0.499991957)$ & $-1.00\cdot10^{-6}$\\
$(0.127850842,\ 0.499992907,\ 0.499991957)$ & $+2.03\cdot10^{-6}$\\\hline
\end{tabular}
\caption{Corollary~\ref{c:Nin} at the couplings of Table~\ref{tab:KBoverlapNum}:
nine locally exponentially stable equilibria, the ninth being the middle
equilibrium of the cusp unfolding.}
\label{tab:KBnineNum}
\end{table}

\paragraph{Interpretation.}
The bare model of Example~\ref{e:KBoverlap} has overlapping minimal
siphons but no positive boundary resident; its self-feedback
augmentation \eqref{eq:KBoverlapCusp} retains the same siphon lattice.
Weak-coupling persistence of combinatorial families of stable
equilibria is already in \cite{Boud26}, and the ninth stable
equilibrium of Corollary~\ref{c:Nin} exists off the cusp parameter and relies
on the positive cubic coefficient. What the theorem adds is the combination:
the stable equilibria coexist with a nondegenerate interior cusp,
the coupling preserves two overlapping minimal siphons, the proper
siphon faces carry nonzero resident equilibria, and distinct,
explicit invasion numbers decide their transverse stability.
This is local coexistence; the theorem does not assert that the cusp
itself crosses either invasion-threshold hypersurface as a coupling is varied.
\section{Threshold responses with spatial diffusion}\label{s:spaTh}
\beT[Autonomous positive Lur'e systems with diagonal diffusion retain the transition at $\bar R_0=1/L$ and admit a uniform upper spreading-speed bound]
[speed bound][SpdBd]
Let $V$, $B$, $W$, $L$, $\beta$ and $\mathcal C_L$ be as in Theorem~\ref{t:KBsharp}, with $\bar R_0=\rho(\beta BWV^{-1})$ and $s$ the spectral abscissa of Theorem~\ref{t:InvF}. Let $\mathcal D=\operatorname{Diag}(d_i)$ with $d_i>0$ and $n\ge1$. For a common response $\Phi\in\mathcal C_L$ applied coordinatewise, consider
\begin{equation}\label{eq:spaTh}
 x_t=\mathcal D\Delta_\zeta x-Vx+B\Phi(\beta Wx),\qquad x(t,\zeta)\in\mathbb R_+^{n_S},\quad\zeta\in\mathbb R^n.
\end{equation}
\begin{enumerate}
\item If $\bar R_0<1/L$, then for every $\Phi\in\mathcal C_L$ and every $\mathcal D$, with $\bar x$ the unique equilibrium of
Theorem~\ref{t:KBsharp}, $p=(V-\beta LBW)^{-1}\mathbf1$ and $c=\min_ip_i^{-1}$, every bounded nonnegative solution of \eqref{eq:spaTh} satisfies
\[
 \sup_\zeta\|x(t,\zeta)-\bar x\|_{\infty,p}\le e^{-ct}\sup_\zeta\|x(0,\zeta)-\bar x\|_{\infty,p},\qquad \|z\|_{\infty,p}=\max_i|z_i|/p_i.
\]
In particular \eqref{eq:spaTh} has no nonconstant entire solution bounded uniformly in $(t,\zeta)$, and hence no traveling front.
\item If $\bar R_0>1/L$, then for every $m\ge2$ the response of Theorem~\ref{t:KBsharp} gives, for every $\mathcal D$, at least $m$ spatially
constant, linearly exponentially stable equilibria of \eqref{eq:spaTh}.
\item Suppose that $\bar R_0>1/L$ and that $\beta LBW-V$ is irreducible. Put
\[
 \bar c(L)=\inf_{\lambda>0}\frac1\lambda\,s\bigl(\lambda^2\mathcal D+\beta LBW-V\bigr)>0.
\]
For every $\Phi\in\mathcal C_L$ with $\Phi(0)=0$, every nonnegative, bounded, compactly supported initial datum, and every $\varepsilon>0$,
\[
 \sup_{|\zeta|\ge(\bar c(L)+\varepsilon)t}\|x(t,\zeta)\|\to0.
\]
\end{enumerate}
\rev{Items 1 and 2 are statements about the whole class $\mathcal C_L$: the value $\bar R_0=1/L$ separates uniform attraction for every response of $\mathcal C_L$ from the existence of responses with arbitrarily many stable equilibria, and it is not a threshold for propagation under one given response. Item 3 is an upper bound on the speed: it neither produces a front nor proves nonlinear linear determinacy.}
\eeT
\begin{proof}
(1) Let $x,y$ be two bounded nonnegative solutions and $z=x-y$. Then
$z_t=\mathcal D\Delta_\zeta z-Vz+\beta BE(t,\zeta)Wz$ with $E$ diagonal and entries in $[0,L]$, a cooperative linear system. Put
$Z(t)=Ce^{-ct}p$, with $C=\sup_\zeta\|z(0,\zeta)\|_{\infty,p}$. Since $(-V+\beta BEW)p\le(-V+\beta LBW)p=-\mathbf1\le-cp$, both $Z-z$ and $Z+z$ are supersolutions
of that system with nonnegative initial data, so the comparison principle gives $|z|\le Z$. Taking $y\equiv\bar x$ gives the estimate. Applying it
to an entire solution from time $-T$, and letting $T\to\infty$, gives $x\equiv\bar x$.

(2) The zero-slope construction in Theorem~\ref{t:KBsharp} gives equilibria with $J=-V$. Its full-support perturbation instead gives nearby Metzler--Hurwitz Jacobians. In either case, Corollary~\ref{c:DifSta} preserves their exponential stability under every positive diagonal diffusion matrix.

(3) By the next-generation theorem, $\bar R_0>1/L$ is equivalent to $s(\beta LBW-V)>0$. Since $\mathcal D$ is positive diagonal, the continuous function
$g(\lambda)=s(\lambda^2\mathcal D+\beta LBW-V)$ satisfies
$g(0)>0$ and $g(\lambda)/\lambda\to+\infty$ as
$\lambda\downarrow0$ or $\lambda\to+\infty$.
It therefore has a minimizer $\lambda_*>0$; irreducibility supplies
its Perron vector $q_*\gg0$ and $\bar c(L)>0$. Since $\Phi(0)=0$ and $\Phi'\le L$, for $x\ge0$ one has $\Phi(\beta Wx)\le\beta LWx$, and $B\ge0$ gives $B\Phi(\beta Wx)\le\beta LBWx$, so $x$ is a subsolution of the cooperative system
$u_t=\mathcal D\Delta_\zeta u+(\beta LBW-V)u$. For every unit vector $e$, the function $Ce^{-\lambda_*(e\cdot\zeta-\bar c(L)t)}q_*$ solves this system and dominates the compactly supported initial data for large $C$
independent of $e$. Choosing $e=\zeta/|\zeta|$, comparison gives $\|x(t,\zeta)\|\le C'e^{-\lambda_*\varepsilon t}$ for $|\zeta|\ge(\bar c(L)+\varepsilon)t$.
\end{proof}

\beC[Diffusion preserves every constructed count of stable equilibria]
[stable equilibria under diffusion][DifSta]
Let $f$ be $C^1$ and cooperative on a convex set $K\subset\mathbb R^{n_S}$ invariant under $x_t=\mathcal D\Delta_\zeta x+f(x)$, with $\mathcal D$ a positive diagonal matrix, and let $e\in K$ be an equilibrium of $f$ with $J=Df(e)$ Hurwitz. Put $p=(-J)^{-1}\mathbf1\gg0$ and $c=\min_ip_i^{-1}$. Then
\[
 u_t=\mathcal D\Delta_\zeta u+Ju\implies\sup_\zeta\|u(t,\zeta)\|_{\infty,p}\le e^{-ct}\sup_\zeta\|u(0,\zeta)\|_{\infty,p},
\]
for bounded solutions on $\mathbb R^n$ and on bounded domains with Neumann or periodic boundary conditions, and $e$ is locally exponentially stable for the nonlinear equation in the supremum norm, for every $\mathcal D$. Consequently, at the same parameter values, the stable equilibria of Theorems~\ref{t:siphMulti}, \ref{t:KBsharp}, \ref{t:Net}, \ref{t:KBoverlapCusp} remain stable spatially constant equilibria, and the nine stable equilibria of Corollary~\ref{c:Nin} remain nine. This does not assert fronts connecting them.
\eeC
\begin{proof}
$J$ is Metzler and Hurwitz, so $-J$ is a nonsingular $M$-matrix and $Jp=-\mathbf1\le-cp$. Put $Z(t)=Ce^{-ct}p$, with $C=\sup_\zeta\|u(0,\zeta)\|_{\infty,p}$. Then $Z-u$ and $Z+u$ are supersolutions of the cooperative system $w_t=\mathcal D\Delta_\zeta w+Jw$ with nonnegative initial data, so comparison gives $|u|\le Z$, the first estimate. For a solution $x$ in $K$ near $e$, $u=x-e$ solves $u_t=\mathcal D\Delta_\zeta u+A(t,\zeta)u$, with $A=\int_0^1Df(e+\theta u)\,d\theta$ Metzler, and by continuity $Ap\le-\frac c2p$ when $\|u\|$ is small. The same comparison, with rate $c/2$, keeps $u$ small and gives the decay.
\end{proof}

\section{Conclusions and remaining problems}\label{s:concShort}

Siphon sets cut out invariant siphon faces. First-order terms on such a face
assemble its transverse invasion matrix, while terms at least quadratic in the
absent variables leave that matrix unchanged and can still create bistability
(Theorem~\ref{t:InvF}). The transverse block is Metzler, so its spectral
abscissa determines linear invasion or decay when nonzero; the critical case
requires nonlinear analysis.

For threshold networks, Theorem~\ref{t:KBsharp} shows that the threshold $\bar R_0=1/L$ is
sharp over the response class with slope bound $L$.
Theorem~\ref{t:siphMulti} instead fixes the siphon lattice and every transverse
matrix, using responses flat at onset. In both, the response depends on the
prescribed count; neither gives arbitrary multistability for a fixed response or
a fixed polynomial degree.

Odds-polynomial responses join these directions. Their zero set fixes the siphon
lattice, their onset slope enters the transverse matrix directly, and their higher
coefficients enter it only through the resident state (Proposition~\ref{p:Odd}).
Theorem~\ref{t:Net} carries the quartic cusp of Theorem~\ref{t:Qua} into a
network over a transversely stable siphon face; singleton siphons force
multiplicative coupling. Overlapping minimal siphons allow additive
cross-couplings, and Theorem~\ref{t:KBoverlapCusp} combines them with nonzero
residents, distinct explicit invasion numbers and an interior cusp. This cusp has
a positive cubic coefficient: it is unstable, the middle equilibrium of its
unfolding is stable, the stable equilibria at the cusp parameter lie elsewhere,
and Corollary~\ref{c:Nin} adds one stable equilibrium to the eight.

Three problems remain. Local bifurcation analysis does not determine which
attractor is reached after a fold. Theorem~\ref{t:Net} gives no explicit coupling
window along its continued family, where $a_1$ moves with the couplings. Finally,
it remains open how the cusp branches intersect the invasion-threshold
hypersurfaces $R_{\Sigma}=1$ associated with overlapping minimal siphons.
Theorem~\ref{t:KBoverlapCusp} proves coexistence and supplies the thresholds, but
does not establish such intersections; the same question arises for asymmetric
directed coupling, with $\gamma$ and $\varepsilon$ independent instead of tied.

\ifjournalversion\else
\clearpage
\part{Additional developments: arXiv version only}
The following sections offer general invasion-wave calculations,
fixed-response examples, front computations, additional Lyapunov
certificates, and background theory.

\paragraph{Related fixed-point results.}
For $T=V^{-1}\Pi$, the nonnegative neural-network map interpretation
is considered in \cite{piotrowski2024fixed}; its weak-scalability
hypotheses cover different response classes from the present flat-step
constructions. The clipped affine family of \cite{Retnaraj2026}
has a timescale-independent equilibrium set, but is only piecewise
smooth at clipping thresholds.
Converse contraction and state-dependent metrics are developed in
\cite{LohmillerSlotine1998,ForniSepulchre2014,Giesl2015}; the
positive-weight criterion is related to \cite{Coogan2019}.
\subsection{Fixed-response activation and explicit cusp calculations}

\paragraph{Onset with fixed responses.}
The next result varies the coefficient of the polynomials $P_i$ instead of
the shape of the responses.

\beT[Increasing higher-order terms creates positive equilibria
while zero remains stable][activation by higher-order terms][ActF]\label{t:KBactivation}
Consider
\[
 \dot x=f(x;s)=-Vx+\Pi_s(x),\qquad
 (\Pi_s)_i=b_i\Phi_i(sP_i(x)),\qquad s\ge0,
\]
where $V$ is a nonsingular $M$-matrix, $b\gg0$, and the invariant
state space is $[0,\bar x]$, $\bar x=V^{-1}b$.
Each $P_i$ is a nonzero polynomial with nonnegative coefficients and
monomials of degree at least two. Assume
$\Phi_i\in C^1([0,\infty))$, with the derivative at zero understood from
the right, and
\[
 \Phi_i(0)=0,\qquad \Phi_i(t)\uparrow1\quad(t\to\infty),
 \qquad \Phi_i'(t)>0\quad(t>0).
\]
The graph containing transfer edges $j\to i$ with $-v_{ij}>0$ and
edges $j\to i$ supplied by pure-power terms $a x_j^k$ in $P_i$, $a>0$,
is assumed strongly connected. Mixed monomials are allowed without
supplying these edges. Then $s_*\in(0,\infty)$, defined by
\begin{equation}\label{eq:KB-activation}
 s_*:=\inf\{s\ge0:\exists x\gg0,\ x\le\bar x,\ \Pi_s(x)\ge Vx\},
\end{equation}
is attained by a positive equilibrium. Positive equilibria exist
exactly for $s\ge s_*$, while zero is locally exponentially stable
for every finite $s$.

At the greatest equilibrium $x_*$ for $s=s_*$, the Jacobian has a
simple zero eigenvalue and a strictly stable complement. If $f$ is
$C^2$ near $(x_*,s_*)$ and
\[
 \ell^\top\partial_x^2f(x_*;s_*)[r,r]\ne0,
 \qquad Jr=0,\quad\ell^\top J=0,\quad r,\ell\gg0,
 \quad\ell^\top r=1,
\]
then, for $s>s_*$ sufficiently close to $s_*$, a saddle-node creates
a locally exponentially stable positive equilibrium and a positive
equilibrium with exactly one unstable eigenvalue.
\eeT
\begin{proof}
A positive subsolution $\Pi_s(x)\ge Vx$ generates a nondecreasing
trajectory bounded by $\bar x$, whose limit is a positive equilibrium.
Increasing $s$ preserves the subsolution inequality. For $0<t<1$,
$V(t\bar x)=tb$ and every $P_i(t\bar x)>0$, so $t\bar x$ is a
subsolution for sufficiently large $s$. For sufficiently small $s$,
\[
 \sup_{x\in[0,\bar x]}\|V^{-1}\partial_x\Pi_s(x)\|_\infty<1,
\]
since $\partial_x\Pi_s\to0$ uniformly on the box. Thus
$x\mapsto V^{-1}\Pi_s(x)$ is a contraction of the box into itself;
its fixed point is zero. Consequently $0<s_*<\infty$.
On bounded $s$-intervals, $\Pi_s(x)=O(\|x\|^2)$ uniformly near zero;
$x=V^{-1}\Pi_s(x)$ therefore excludes small nonzero equilibria.
Limits as $s\downarrow s_*$ give a nonzero equilibrium at $s_*$.
At an equilibrium, a zero coordinate forces its graph predecessors to
be zero; strong connectivity makes every nonzero equilibrium positive.
Also $J(0)=-V$, giving stability of zero.

The decreasing trajectory from $\bar x$ gives the greatest equilibrium
$x_*$. It lies strictly below $\bar x$, since
\[
 \bar x-x_*=V^{-1}(b-\Pi_{s_*}(x_*))\gg0.
\]
Its Jacobian is irreducible Metzler. If $s(J)>0$, the positive Perron
vector gives a strict subsolution $x_*+\varepsilon r$ inside the same
box, contradicting maximality. If $s(J)<0$, the implicit-function
theorem continues the positive equilibrium below $s_*$. Hence $s(J)=0$
is simple and all other eigenvalues have negative real parts.
Finally,
\[
 \ell^\top\partial_s f(x_*;s_*)
 =\sum_i\ell_i b_i\Phi_i'(s_*P_i(x_*))P_i(x_*)>0.
\]
The saddle-node test of Proposition~\ref{p:sadCus} applies when
its quadratic coefficient is nonzero.
That coefficient is negative because positive equilibria do not exist
below $s_*$. The two branches above it have zero and one unstable
eigenvalues, respectively; the latter is a saddle in dimension at least two.
\end{proof}

\paragraph{Homogeneous branches, saddle-nodes, and heterogeneous cusps.}

\beT[Balanced response arguments preserve scalar branches and their first instability]
[homogeneous branches][KBhom]
For a production model \eqref{eq:KB}, suppose
$V\mathbf1=v\mathbf1$, $v>0$, and
$\Pi(s\mathbf1)=r(s)\mathbf1$ throughout the relevant interval.
The states $x=s\mathbf1$, in which all variables are equal, then obey
\begin{equation}\label{eq:KBscalar}
 \dot s=-vs+r(s).
\end{equation}
If the Jacobian $J$ at such an equilibrium is Metzler, its spectral
abscissa is $s(J)=-v+r'(s)$.
When this number is zero the eigenvalue zero is semisimple, even
without irreducibility. If $J$ is irreducible it is simple, and all
heterogeneous eigenvalues have negative real parts.
For $r(s)=b\Phi(h+\beta s)$ with $b,\beta>0$,
write $\varphi=\Phi'$ and assume the derivatives used below exist and
are continuous. A scalar saddle-node candidate
satisfies
\[
 vs=b\Phi(z),\quad z=h+\beta s,\quad b\beta\varphi(z)=v.
\]
It is nondegenerate for variation of $h$ if $\varphi'(z)\ne0$.
If $\varphi(z_c)>0$, $\varphi'(z_c)=0$, and
$\varphi''(z_c)\ne0$, the scalar equation has a nondegenerate two-parameter cusp at
\[
 s_c=\frac b v\Phi(z_c),\qquad
 \beta_c=\frac v{b\varphi(z_c)},\qquad
 h_c=z_c-\beta_cs_c.
\]
At a saddle-node or cusp with irreducible Metzler full Jacobian, the scalar
bifurcation has only stable complementary directions in the full system.
For higher-order balanced $P_i$, one may take
$r(s)=b\Phi(h+\sum_k a_ks^k)$ instead.
\eeT
\begin{proof}
Invariance follows by substitution. With $H=\mathbf1$,
Theorem~\ref{t:Spe}(1) gives spectral equality. More explicitly,
$J\mathbf1=(-v+r'(s))\mathbf1$ implies
$\|e^{Jt}\|_\infty=e^{(-v+r'(s))t}$. This proves spectral ordering
and, at zero, boundedness of the semigroup excludes Jordan blocks.
The scalar derivatives are
$f_s=-v+b\beta\varphi$, $f_{ss}=b\beta^2\varphi'$, and
$f_{sss}=b\beta^3\varphi''$. At the cusp the unfolding determinant
in $(h,\beta)$ is $b^2\varphi(z_c)^2>0$.
Proposition~\ref{p:sadCus} verifies these scalar bifurcations.
The cusp or saddle-node lifts under Theorem~\ref{t:Spe}(3); in a
reducible network several center directions may instead coexist.
\end{proof}

For $V=I$ and row-stochastic $W$, the spectral ordering and the semisimple or simple zero
eigenvalue are in \cite{Boud26}; Theorem~\ref{t:KBhom} extends that argument to
$V\mathbf1=v\mathbf1$ and nonlinear balanced production.

For a unimodal density with $\varphi''(z_c)<0$, the cubic coefficient $f_{sss}$ is negative, so the outer branches of the cusp
attract, and for $\beta>\beta_c$ the two folds bound a hysteresis loop in $h$ (Figure~\ref{fig:hyst}).

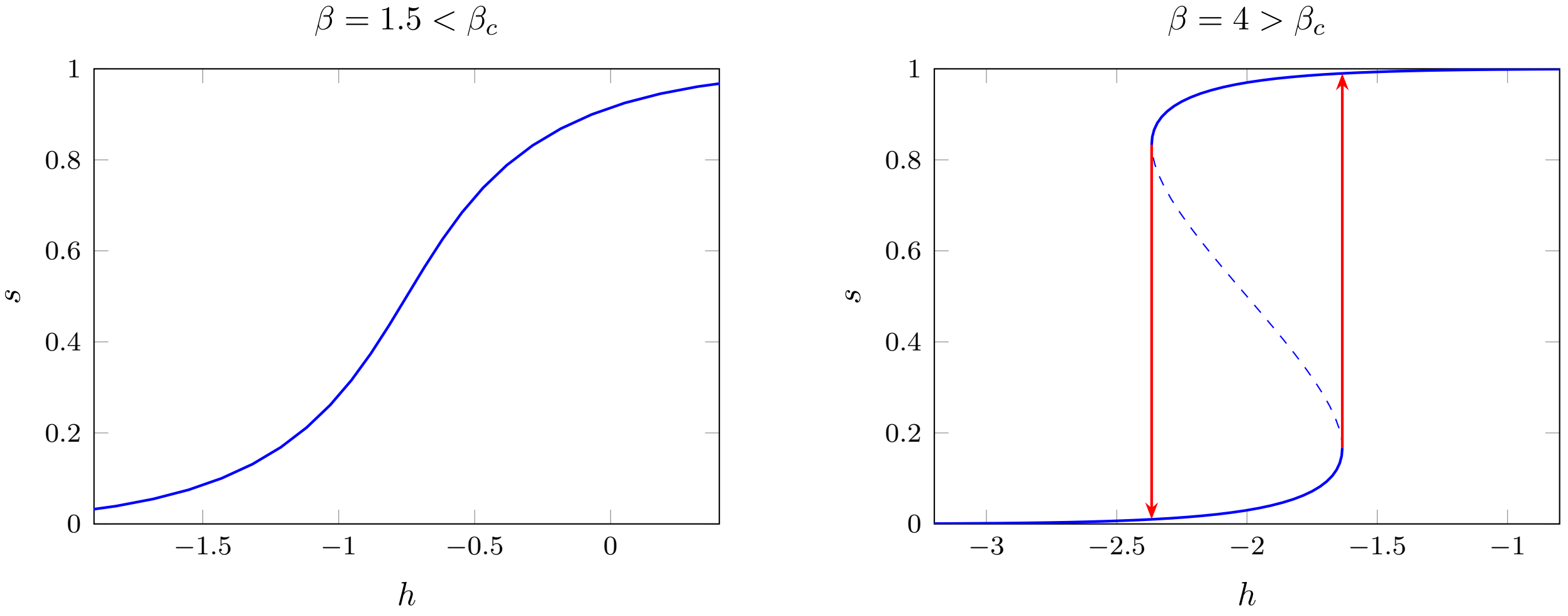
\begin{figure}[H]
\centering
\begin{minipage}{.48\linewidth}
\centering
\begin{tikzpicture}
\begin{axis}[width=\linewidth,height=.78\linewidth,xlabel={$h$},ylabel={$s$},ymin=0,ymax=1,
 tick label style={font=\scriptsize},label style={font=\small},title style={font=\small},xmin=-1.9,xmax=0.4,title={$\beta=1.5<\beta_c$}]
\addplot[blue,thick] coordinates {(-3.2010,0.0007) (-3.0418,0.0012) (-2.8830,0.0020) (-2.7249,0.0033) (-2.5679,0.0052) (-2.4123,0.0082) (-2.2588,0.0125) (-2.1081,0.0188) (-1.9611,0.0274) (-1.8188,0.0392) (-1.6822,0.0548) (-1.5524,0.0749) (-1.4304,0.1003) (-1.3170,0.1314) (-1.2128,0.1685) (-1.1178,0.2119) (-1.0316,0.2611) (-0.9534,0.3156) (-0.8817,0.3745) (-0.8147,0.4364) (-0.7500,0.5000) (-0.6853,0.5636) (-0.6183,0.6255) (-0.5466,0.6844) (-0.4684,0.7389) (-0.3822,0.7881) (-0.2872,0.8315) (-0.1830,0.8686) (-0.0696,0.8997) (0.0524,0.9251) (0.1822,0.9452) (0.3188,0.9608) (0.4611,0.9726) (0.6081,0.9812) (0.7588,0.9875) (0.9123,0.9918) (1.0679,0.9948) (1.2249,0.9967) (1.3830,0.9980) (1.5418,0.9988) (1.7010,0.9993)};
\end{axis}
\end{tikzpicture}
\end{minipage}\hfill
\begin{minipage}{.48\linewidth}
\centering
\begin{tikzpicture}
\begin{axis}[width=\linewidth,height=.78\linewidth,xlabel={$h$},ylabel={$s$},ymin=0,ymax=1,
 tick label style={font=\scriptsize},label style={font=\small},title style={font=\small},xmin=-3.2,xmax=-0.8,title={$\beta=4>\beta_c$}]
\addplot[blue,thick] coordinates {(-3.8003,0.0001) (-3.7296,0.0001) (-3.6588,0.0001) (-3.5882,0.0002) (-3.5176,0.0002) (-3.4470,0.0003) (-3.3765,0.0004) (-3.3061,0.0005) (-3.2358,0.0006) (-3.1657,0.0008) (-3.0957,0.0010) (-3.0259,0.0013) (-2.9564,0.0016) (-2.8872,0.0020) (-2.8183,0.0025) (-2.7499,0.0031) (-2.6820,0.0038) (-2.6148,0.0047) (-2.5482,0.0058) (-2.4825,0.0071) (-2.4177,0.0086) (-2.3541,0.0104) (-2.2917,0.0125) (-2.2308,0.0150) (-2.1715,0.0179) (-2.1141,0.0212) (-2.0588,0.0251) (-2.0058,0.0295) (-1.9553,0.0346) (-1.9076,0.0404) (-1.8629,0.0470) (-1.8216,0.0543) (-1.7838,0.0626) (-1.7498,0.0718) (-1.7198,0.0820) (-1.6940,0.0933) (-1.6726,0.1056) (-1.6559,0.1191) (-1.6438,0.1338) (-1.6365,0.1497) (-1.6341,0.1668)};
\addplot[blue,thick] coordinates {(-2.3659,0.8332) (-2.3635,0.8503) (-2.3562,0.8662) (-2.3441,0.8809) (-2.3274,0.8944) (-2.3060,0.9067) (-2.2802,0.9180) (-2.2502,0.9282) (-2.2162,0.9374) (-2.1784,0.9457) (-2.1371,0.9530) (-2.0924,0.9596) (-2.0447,0.9654) (-1.9942,0.9705) (-1.9412,0.9749) (-1.8859,0.9788) (-1.8285,0.9821) (-1.7692,0.9850) (-1.7083,0.9875) (-1.6459,0.9896) (-1.5823,0.9914) (-1.5175,0.9929) (-1.4518,0.9942) (-1.3852,0.9953) (-1.3180,0.9962) (-1.2501,0.9969) (-1.1817,0.9975) (-1.1128,0.9980) (-1.0436,0.9984) (-0.9741,0.9987) (-0.9043,0.9990) (-0.8343,0.9992) (-0.7642,0.9994) (-0.6939,0.9995) (-0.6235,0.9996) (-0.5530,0.9997) (-0.4824,0.9998) (-0.4118,0.9998) (-0.3412,0.9999) (-0.2704,0.9999) (-0.1997,0.9999)};
\addplot[blue,dashed] coordinates {(-1.6341,0.1668) (-1.6372,0.1877) (-1.6465,0.2102) (-1.6619,0.2342) (-1.6830,0.2596) (-1.7095,0.2864) (-1.7410,0.3144) (-1.7770,0.3435) (-1.8168,0.3736) (-1.8597,0.4045) (-1.9051,0.4360) (-1.9521,0.4679) (-2.0000,0.5000) (-2.0479,0.5321) (-2.0949,0.5640) (-2.1403,0.5955) (-2.1832,0.6264) (-2.2230,0.6565) (-2.2590,0.6856) (-2.2905,0.7136) (-2.3170,0.7404) (-2.3381,0.7658) (-2.3535,0.7898) (-2.3628,0.8123) (-2.3659,0.8332)};
\draw[red,thick,-{Stealth[length=5pt]}] (axis cs:-1.6341,0.1668) -- (axis cs:-1.6341,0.9900);
\draw[red,thick,-{Stealth[length=5pt]}] (axis cs:-2.3659,0.8332) -- (axis cs:-2.3659,0.0100);
\end{axis}
\end{tikzpicture}
\end{minipage}
\caption{Hysteresis on the homogeneous branch of Theorem~\ref{t:KBhom}, $\dot s=-s+\Phi_{\mathrm G}(h+\beta s)$, for the Gaussian response:
the equilibria are $s=\Phi_{\mathrm G}(z)$, $h=z-\beta\Phi_{\mathrm G}(z)$, and $\beta_c=\sqrt{2\pi}$. Left: the response to $h$ is
gradual and reversible. Right: solid branches are stable and the dashed branch is unstable. As $h$ increases the state jumps up at
$h\approx-1.634$; as $h$ decreases it jumps down at $h\approx-2.366$.}
\label{fig:hyst}
\end{figure}

\beT[Two identical response populations admit cusp formulas with linear exchange]
[cusps with linear exchange][KBexchange]
Consider $\dot x=-V_\kappa x+\Pi(x)$ for two populations with
the same strictly increasing response $\Phi:\mathbb R\to(0,1)$ and coefficients,
\[
 V_\kappa=\begin{pmatrix}v+\kappa&-\kappa\\-\kappa&v+\kappa\end{pmatrix},
 \quad W_\varepsilon=\begin{pmatrix}1-\varepsilon&\varepsilon\\
 \varepsilon&1-\varepsilon\end{pmatrix},
 \quad \Pi_i=b\Phi(h+\beta(W_\varepsilon x)_i),
\]
where $v,b,\beta>0$, $\kappa\ge0$, and
$0\le\varepsilon\le1/2$. The terms with coefficient $\kappa$
transfer activity between the two variables and cancel in their sum.
Write the two response arguments as $a+s$ and $a-s$, and put
\[
 M(a,s)=\tfrac12(\Phi(a+s)+\Phi(a-s)),\quad
 D(a,s)=\tfrac12(\Phi(a+s)-\Phi(a-s)),\quad
 \widehat\beta=\beta b/v.
\]
The heterogeneous equilibria are parameterized, whenever $D\ne0$, by
\begin{align*}
 m&=\frac b v M,&d&=\frac b{v+2\kappa}D,&x_{1,2}&=m\pm d,\\
 h&=a-\widehat\beta M,&
 \varepsilon&=\frac12\left(1-
 \frac{(v+2\kappa)s}{\beta bD}\right).
\end{align*}
At the same arguments $(a,s)$, define the parameter map without exchange,
with unit removal and production coefficient, by
\[
 U=a-\widehat\beta M,\qquad
 \varepsilon_0=\frac12\left(1-\frac{s}{\widehat\beta D}\right).
\]
The original parameters satisfy exactly
\begin{equation}\label{eq:KB-cusp-affine}
 (h,\varepsilon)=
 \left(U,(1+2\kappa/v)\varepsilon_0-\kappa/v\right).
\end{equation}
Consequently its generic parameter-map folds and cusps are preserved
by this affine transformation, retaining only those with admissible
$0\le\varepsilon\le1/2$ and positive activities.
Suppose $\Phi$ is $C^7$ near $z_c$ and, writing $\varphi=\Phi'$,
\[
 \varphi(z_c)>0,\qquad \varphi'(z_c)=0,\qquad
 \varphi''(z_c)<0.
\]
Put $a_0=b/v$, $\beta_c=1/(a_0\varphi(z_c))$,
$x_c=a_0\Phi(z_c)$, $h_c=z_c-\beta_cx_c$, and
$\delta=\beta-\beta_c\downarrow0$.
The two cusps related by interchanging the populations extend the
zero-exchange calculation of \cite{Boud26}. They satisfy,
for $\kappa=O(\delta)$,
\begin{align*}
 \varepsilon_c&=\frac{\delta}{3\beta_c}-\frac\kappa v+o(\delta),\\
 |x_{1,c}-x_{2,c}|&=
 2\left[-\frac{2}{a_0\varphi''(z_c)\beta_c^4}\right]^{1/2}
 \delta^{1/2}+o(\delta^{1/2}),\\
 h_c(\beta)-h_c&=-x_c\delta+o(\delta).
\end{align*}
Let $\varepsilon_{0,c}(\widehat\beta)$ be the value of $\varepsilon_0$ at the cusps. Their
couplings are admissible exactly when $\varepsilon_{0,c}(\widehat\beta)\ge\kappa/(v+2\kappa)$.
For fixed $\kappa>0$, these cusps leave the physical range
$\varepsilon\ge0$ sufficiently near the scalar cusp.
\eeT
\begin{proof}
The mean equilibrium equation is $vm=bM$, and the contrast equation
is $(v+2\kappa)d=bD$. The fields satisfy
$a=h+\beta m$ and $s=\beta(1-2\varepsilon)d$.
Eliminating $m,d$ gives the parameterization and
\eqref{eq:KB-cusp-affine}. Its target transformation is invertible,
so it preserves Whitney fold and cusp type. The source transformation
from activities to fields is locally invertible when
$\beta(1-2\varepsilon)\ne0$.

\emph{Unit model.} Take $v=b=1$, $\kappa=0$, and write
$\delta_u=\beta-\beta_c$. This is the model of \cite{Boud26}, whose theorem on
heterogeneous cusps born from the scalar cusp gives two exchange-related generic Whitney
cusps of the parameter map $(a,s)\mapsto(h,\varepsilon)$ and the three displayed laws.
We record an exact form of its nondegeneracy computation.
Let $\Xi=\partial_ah\,\partial_s\varepsilon-\partial_sh\,\partial_a\varepsilon$
and let $D_k\Xi=\partial_a\Xi\,\partial_s\varepsilon-\partial_s\Xi\,\partial_a\varepsilon$
be its derivative along the kernel direction $k=(\partial_s\varepsilon,-\partial_a\varepsilon)$. Put
\[
 R(a,s)=\frac{D(a,s)}s
       =\frac12\int_{-1}^1\varphi(a+\tau s)\,d\tau,
\]
\[
 H=(1-\beta M_a)R_s+\beta M_sR_a,
 \qquad \mathcal K=R_s\partial_a-R_a\partial_s.
\]
Since $\mathcal K R=0$, direct differentiation gives the exact identities
\[
 \Xi=\frac{H}{2\beta R^2},\qquad
 D_k\Xi=\frac{\mathcal K H}{4\beta^2R^4},\qquad
 D_k^2\Xi=\frac{\mathcal K^2H}{8\beta^3R^6}.
\]
With $\delta_u=t^2$, $s=tS$, $a-z_c=t^2A$, the cusp branches have
\[
 S(0)=\pm S_*,\quad S_*^2=-\frac{2}{\beta_c^2\varphi''(z_c)},\qquad
 A(0)=\frac{\varphi'''(z_c)}{4\beta_c^2\varphi''(z_c)^2},
\]
and, with $\varphi''=\varphi''(z_c)$, along them
\[
 \partial_s\Xi=\frac{\varphi''}{3}t^2+o(t^2),\qquad
 \partial_s\varepsilon=\frac{\beta_c\varphi''S}{6}t+o(t),\qquad
 D_k^2\Xi=-\frac{2\beta_c^2\varphi''^3S}{27}t^3+o(t^3).
\]
Thus $d\Xi\ne0$, the parameter map has rank one, and
$D_k^2\Xi\ne0$ for all sufficiently small $t>0$.

\emph{General model.} The activities are $a_0=b/v$ times the unit-model
ones. With $\widehat\beta-\widehat\beta_c=a_0\delta$ and
$\widehat\beta_c=a_0\beta_c$, the unit-model formulas give the displayed ones
for $\kappa=0$; the transformation \eqref{eq:KB-cusp-affine} gives
$\varepsilon_c=(1+2\kappa/v)\varepsilon_{0,c}-\kappa/v$, and the contrast carries the
factor $b/(v+2\kappa)$ in place of $b/v$. For $\kappa=O(\delta)$, the product
$(2\kappa/v)\varepsilon_{0,c}$ is $O(\delta^2)$ and the correction
to the contrast is $O(\delta^{3/2})$; the leading coupling shift
is $-\kappa/v$.
By \eqref{eq:KB-cusp-affine}, $0\le\varepsilon\le1/2$ is equivalent to
$\kappa/(v+2\kappa)\le\varepsilon_0\le1/2$, which is the admissibility test of the statement.
Near the cusp the other two-dimensional eigenvalue is negative:
at the same fields, its trace satisfies exactly
\[
 \operatorname{tr}J_\kappa
 =v[-2+(1-\varepsilon_{0,c})(q_1+q_2)]
  +\kappa[(1-2\varepsilon_{0,c})(q_1+q_2)-2],
\]
where
$q_i=\widehat\beta\varphi(a\pm s)$.
Here $q_1+q_2=2+o(\delta)$ and
$\varepsilon_{0,c}=\delta/(3\beta_c)+o(\delta)$.
For $\kappa=O(\delta)$, the trace is consequently
$-2v\delta/(3\beta_c)+o(\delta)<0$.
Thus physical cusps have one stable complementary direction.
Theorem~\ref{t:Spe}(3) lifts them to an exact cooperative
network realization whose Jacobian at the cusp is irreducible.
\end{proof}

\paragraph{Application: persistence of the two exchange cusps.}
For each fixed sufficiently small $\beta-\beta_c>0$ for which
Theorem~\ref{t:KBexchange} supplies two distinct nondegenerate
full-network cusps, Proposition~\ref{p:sadCus} continues both under
sufficiently small $C^3$ asymmetry or loss of equitability.
Apply it in disjoint neighbourhoods of the two cusp states.
Their stimulus and coupling values need no longer coincide,
and the permissible perturbation size may depend on $\beta-\beta_c$.
Interior admissibility persists; at a parameter boundary such as
$\varepsilon=0$, check it separately using the exchange-induced
admissibility inequality above.

\beT[Unequal stimuli preserve both threshold cusps precisely within an explicit asymmetry interval]
[][cuspAsym]
Consider
\[
 \dot x_{1,2}=-x_{1,2}+
 \Phi\bigl(h\pm\eta+\beta[(1-k)x_{1,2}+kx_{2,1}]\bigr),
 \qquad 0\le k\le1/2,\quad\beta>0.
\]
Let $\Phi:\mathbb R\to(0,1)$ be a strictly increasing CDF, $C^5$ near zero,
with $\Phi(-q)=1-\Phi(q)$. Put $\varphi=\Phi'$, assume $\varphi(0)>0$, $\varphi''(0)<0$,
and set $\beta_c=1/\varphi(0)$. For $\beta=\beta_c+\delta$, $\delta>0$
sufficiently small, there is a unique sufficiently small $s>0$ with
$\beta\varphi(s)=1$. Define
\[
 D=\Phi(s)-\tfrac12>0,\qquad \eta_* =\beta D-s>0,\qquad
 k_0=\frac{\eta_*}{2\beta D}.
\]
For $\zeta=\pm1$, the two continuations of the heterogeneous cusps at
$\eta=0$ are
\[
 x_\zeta^*=(\tfrac12+\zeta D,\tfrac12-\zeta D),\qquad
 h_\zeta^*=-\beta/2,\qquad k_\zeta^*=k_0+\frac{\zeta\eta}{2\beta D}.
\]
Both candidates lie in $0<k<1/2$ and are nondegenerate cusps of
the $(h,k)$ family precisely while
$|\eta|<\eta_*(\beta)$, with eigenvalues $0,-2k_\zeta^*<0$.
At $|\eta|=\eta_*$ one reaches $k=0$ with two zero eigenvalues,
ceasing to be a one-dimensional-centre cusp; beyond it that candidate
has inadmissible $k<0$. The other cusp remains nondegenerate for
asymmetries near this boundary. Moreover,
\[
 s=\sqrt{\frac{-2}{\beta_c^2\varphi''(0)}}\delta^{1/2}+O(\delta^{3/2}),\qquad
 \eta_* =\frac{2}{3\beta_c}\sqrt{\frac{-2}{\beta_c^2\varphi''(0)}}
 \delta^{3/2}+O(\delta^{5/2}).
\]
\eeT
\begin{proof}
In coordinates $m=(x_1+x_2)/2-1/2$, $d=(x_1-x_2)/2$, at
$h=-\beta/2$, $m=0$ the arguments are $q_1=\eta+\beta(1-2k)d$, $q_2=-q_1$.
The mean equation vanishes by symmetry; setting $d=\zeta D$, $q_1=\zeta s$
gives $\zeta s=\eta+\beta(1-2k)\zeta D$, hence $k=k_\zeta^*$.
Evenness of $\varphi$ and $\beta\varphi(s)=1$ give
\[
 J=-I+\begin{pmatrix}1-k&k\\k&1-k\end{pmatrix},
\]
with eigenvalues $0,-2k$ and mean-direction kernel for $k>0$.

Write the equations as $f_m,f_d$. Since $\partial_df_d=-2k$, solve
$f_d=0$ for $d=d(m,h,k)$ and put $g=f_m(m,d(m,h,k),h,k)$.
At either candidate $d_m=0$. Differentiating $f_d=0$ twice gives
\[
 d_{mm}=\frac{\zeta\beta^2\varphi'(s)}{2k},\qquad
 d_k=-\frac{\zeta D}{k}.
\]
Substitution into the derivatives of $g$ gives
\[
 g=g_m=g_{mm}=0,\qquad
 g_{mmm}=\beta^3\varphi''(s)+\frac{3\beta^4(1-2k)\varphi'(s)^2}{2k},
\]
\[
 g_h=1/\beta,\quad g_k=0,\quad g_{mk}=-\frac{\beta^2D\varphi'(s)}k,
 \qquad
 \det\begin{pmatrix}g_h&g_k\\g_{mh}&g_{mk}\end{pmatrix}
 =-\frac{\beta D\varphi'(s)}k\ne0.
\]
These formulas apply to both signs because $\varphi'$ is odd;
$\varphi'(s)<0$ for small positive $s$.
Evenness also gives
\[
 \varphi(s)=\varphi(0)+\tfrac12\varphi''(0)s^2+O(s^4),\quad
 s^2=-\frac{2\delta}{\beta_c^2\varphi''(0)}+O(\delta^2),
\]
\[
 D=\varphi(0)s+\tfrac16\varphi''(0)s^3+O(s^5).
\]
Thus
\[
 \eta_* =\frac\delta{\beta_c}s+\frac{\beta_c}6\varphi''(0)s^3
 +O(\delta^{5/2})=\frac{2\delta}{3\beta_c}s+O(\delta^{5/2}),
 \qquad k_0=\frac\delta{3\beta_c}+O(\delta^2),
\]
proving the asymptotics. Directly,
$\eta_* =\beta\int_0^s[\varphi(q)-\varphi(s)]\,dq>0$.
For $|\eta|<\eta_*$, $0<k_\zeta^*<2k_0<1/2$.
The cubic coefficient decreases with $k$ and at its upper endpoint is
\[
 g_{mmm}|_{k=2k_0}=-\tfrac72\beta_c^3\varphi''(0)+o(1)>0.
\]
Proposition~\ref{p:sadCus} therefore makes both candidates
nondegenerate cusps unfolded by $(h,k)$.
At the boundary the couplings are $0,2k_0$: the first has $J=0$,
whereas the second retains its nonzero cubic coefficient and stable
complementary eigenvalue, which persist nearby.
\end{proof}

\paragraph{Remarks.}
The exact inequality $|\eta|<\beta(\Phi(s)-1/2)-s$, $\beta \Phi'(s)=1$,
concerns the two inherited cusps, not possible distant ones.
Their states and common stimulus remain fixed, while
$k_+^*-k_-^*=\eta/(\beta D)$. Since $D=O(\delta^{1/2})$ and
$k_0=O(\delta)$, stimulus asymmetry of order $\delta^{3/2}$ changes
the coupling by its own order. For the standard Gaussian response $\Phi_{\mathrm G}$,
\[
 \beta_c=\sqrt{2\pi},\qquad s=\sqrt{2\log(\beta/\beta_c)},
\]
\[
 \begin{aligned}
 \eta_* &=\beta\left[\Phi_{\mathrm G}\!\left(\sqrt{2\log(\beta/\beta_c)}\right)-\tfrac12\right]
 -\sqrt{2\log(\beta/\beta_c)}\\
 &=\frac{2\sqrt2}{3\beta_c^{3/2}}(\beta-\beta_c)^{3/2}
 +O((\beta-\beta_c)^{5/2}).
 \end{aligned}
\]
Unequal removal coefficients or asymmetric response shapes require
separate unfolding calculations; these formulas assume a symmetric
response and unequal additive stimuli.

\section{Spatial spreading via transversal invasion waves}\label{s:KP}

\beAs[Hypotheses of the linear spreading theory]
\label{as:lindet}
{\color{revcolor}
Throughout this section the following are assumed, and nothing beyond them is
used for a linear speed.
\begin{enumerate}
\item The diffusion matrix $\mathcal D=\diag(d_i)$ is nonnegative diagonal and
nonzero, and the medium is homogeneous and isotropic, so that the Laplacian is
the same in every direction.
\item The resident state is spatially constant, and the transverse block $M$ at
it is irreducible Metzler with $s(M)>0$.
\item The linearization at the resident state is the one whose speed is computed.
No claim of linear determinacy is made: whether the nonlinear front travels at
the linear speed is a separate question, settled here in no case, and the
comparison with a nonlinear front is numerical wherever one is shown.
\item For a bound in the nonlinear system the response satisfies
$\Phi(0)=0$ and $\Phi'\le L$, so that $\Phi(\beta Wx)\le\beta LWx$ on the
orthant; this gives an upper bound on the speed, never a lower one.
\end{enumerate}}
\eeAs

\label{s:nonscalar-waves}

In simple monostable examples like rabies \cite{ruan2017spatiotemporal}, $R_0>1$ is equivalent to instability of the disease-free state and to invasion at positive speed. In general, for a spatially homogeneous resident with an irreducible transverse block, $R_\Sigma>1$ is equivalent to positive linear spreading speed; nonlinear invasion can also occur when $R_\Sigma<1$.

The transverse Jacobian block of a siphon face is Metzler, so the Freidlin--G\"artner--Weinberger--Lewis--Li formula gives the spreading speed of the linearized invasion equation. Identifying this speed with that of the nonlinear system requires additional linear-determinacy hypotheses.

\subsection{The transversal invasion wave optimization problem}
\label{s:opt}

\beL[Entrywise bounds for Metzler exponentials, and a uniform saddle-point expansion][Metzler bounds and the Laplace step][trotter]
{\color{revcolor}
Let $M\in\mathbb R^{m\times m}$ be Metzler, let $N$ be real diagonal and let
$\mathcal D=\diag(d_i)$ with $d_i\ge0$.
\begin{enumerate}
\item For every $t\ge0$, entrywise $\bigl|e^{t(M+iN)}\bigr|\le e^{tM}$, and
$\mu\mapsto e^{t(M+\mu\mathcal D)}$ is nondecreasing entrywise.
\item Suppose in addition that $M$ is irreducible, $s(M)>0$ and
$\mathcal D\ne0$. Put $G(\mu)=s(M+\mu\mathcal D)$ and
$g(\lambda)=G(\lambda^2)$. For $v\in\mathbb R$ let
$g'(\lambda_v)=v$, and write
$\psi(v)=g(\lambda_v)-v\lambda_v$,
$a_v=G'(\lambda_v^2)>0$ and $b_v=g''(\lambda_v)>0$.
Let $p_v,q_v\gg0$ be the left and right Perron vectors of
$M+\lambda_v^2\mathcal D$, normalized by $p_v^\top q_v=1$.
For $0\ne\phi\in C_c^\infty(\mathbb R^n;\mathbb R_+^m)$, the solution
$Z_t=\mathcal D\Delta Z+MZ$, $Z(0,\cdot)=\phi$, satisfies
\[
 Z(t,vte)=\frac{e^{t\psi(v)}}
 {\sqrt{2\pi t b_v}\,(4\pi t a_v)^{(n-1)/2}}
 \left(q_v\int_{\mathbb R^n}e^{\lambda_v e\cdot y}
 p_v^\top\phi(y)\,dy+o(1)\right).
\]
Here $e$ is any unit vector, and the vector-valued remainder is uniform
for $v$ in compact subsets of $\mathbb R$ and $|e|=1$, including when
some $d_i$ vanish.
\end{enumerate}}
\eeL
\begin{proof}
{\color{revcolor}
(1) Since $M$ is Metzler, $M+sI\ge0$ for $s$ large, so
$e^{tM}=e^{-st}e^{t(M+sI)}\ge0$ for $t\ge0$. The Trotter product formula gives
$e^{t(M+iN)}=\lim_{k\to\infty}\bigl(e^{tM/k}e^{itN/k}\bigr)^k$, and each factor
$e^{itN/k}$ is diagonal with entries of modulus one, so
$\bigl|\bigl(e^{tM/k}e^{itN/k}\bigr)^k\bigr|\le\bigl(e^{tM/k}\bigr)^k=e^{tM}$
entrywise. For the second claim write
$M+\mu'\mathcal D=(M+\mu\mathcal D)+(\mu'-\mu)\mathcal D$ with
$(\mu'-\mu)\mathcal D\ge0$ for $\mu\le\mu'$, and use the variation-of-constants
formula, whose integrand is nonnegative.

(2) By Perron--Frobenius theory \cite{Seneta,Nussbaum1986} and
Cohen's diagonal convexity theorem \cite{Cohen81}, $G$ is real analytic, convex, and
$G'(\mu)=p(\mu)^\top\mathcal Dq(\mu)>0$, where the Perron vectors are
normalized by $p^\top q=1$. Hence
$g''(\lambda)=2G'(\lambda^2)+4\lambda^2G''(\lambda^2)
\ge2G'(0)>0$, and $g'$ maps $\mathbb R$ bijectively onto $\mathbb R$.
Fix a compact interval $K$ of velocities. All constants below are uniform
for $v\in K$ and $|e|=1$.
With $\widehat\phi(k)=\int e^{-ik\cdot y}\phi(y)\,dy$, Fourier inversion gives
\[
 Z(t,\zeta)=(2\pi)^{-n}\int_{\mathbb R^n}e^{ik\cdot\zeta}
 e^{t(M-(k\cdot k)\mathcal D)}\widehat\phi(k)\,dk.
\]

Write $\eta_1=e\cdot\eta$ and $\eta_\perp=\eta-\eta_1e$. On shifting
the Fourier contour to $k=\eta+i\lambda_ve$, the matrix exponent is
\[
 B_{v,e}(\eta)=M+
 (\lambda_v^2-|\eta|^2-2i\lambda_v\eta_1)\mathcal D.
\]
The simple eigenvalue $\nu_{v,e}(\eta)$ continuing $g(\lambda_v)$,
and its spectral projection $P_{v,e}(\eta)$, are analytic for
$|\eta|<\eta_0$, with a common $\eta_0>0$. Taylor expansion gives
\[
 \begin{aligned}
 \nu_{v,e}(\eta)&=g(\lambda_v)-iv\eta_1
       -\tfrac12 b_v\eta_1^2-a_v|\eta_\perp|^2+O(|\eta|^3),\\
 P_{v,e}(\eta)&=q_vp_v^\top+O(|\eta|).
 \end{aligned}
\]
Since $a_v,b_v$ are bounded away from zero on $K$, decreasing $\eta_0$
if necessary gives
$\Re\nu_{v,e}(\eta)\le g(\lambda_v)-c_0|\eta|^2$.
The remaining spectrum is separated from $g(\lambda_v)$ at $\eta=0$;
compactness and the spectral projection formula bound its semigroup
contribution by $Ce^{t(g(\lambda_v)-\delta_0)}$ for $|\eta|<\eta_0$.

For $|\eta|\ge\eta_0$, part~(1) and entrywise monotonicity give
\[
 \bigl|e^{tB_{v,e}(\eta)}\bigr|
 \le e^{t(M+(\lambda_v^2-|\eta|^2)\mathcal D)}
 \le e^{t(M+(\lambda_v^2-\eta_0^2)\mathcal D)}.
\]
The last semigroup is bounded by
$Ce^{tG(\lambda_v^2-\eta_0^2)}$ on $K$.
Strict positivity of $G'$ gives the uniform gap
\[
 \delta_1=\min_{v\in K}
 \bigl[G(\lambda_v^2)-G(\lambda_v^2-\eta_0^2)\bigr]>0.
\]
Thus the part of the Fourier integral with $|\eta|\ge\eta_0$ is
$O(e^{t(g(\lambda_v)-\delta_1)})$; the shifted transform
$\widehat\phi(\eta+i\lambda_ve)$ is integrable and rapidly decreasing,
uniformly on $K$ and the unit sphere. These estimates do not require
$\mathcal D$ to be invertible.

On $t^{-2/5}\le|\eta|<\eta_0$ the preceding local spectral bounds give
an exponentially smaller contribution
$O(e^{tg(\lambda_v)-c_1t^{1/5}})$.
On $|\eta|<t^{-2/5}$, Taylor's formula gives a relative error
$O(t|\eta|^3)+O(|\eta|)=O(t^{-1/5})$ after factoring out the
principal Gaussian. The linear imaginary term $-iv\eta_1$ cancels
$e^{itv\eta_1}$, and the leading Fourier amplitude is
\[
 q_vp_v^\top\widehat\phi(i\lambda_ve)
 =q_v\int_{\mathbb R^n}e^{\lambda_ve\cdot y}
 p_v^\top\phi(y)\,dy.
\]
The Gaussian integral, with the Fourier factor $(2\pi)^{-n}$, equals
\[
 \frac1{(2\pi)^n}\int_{\mathbb R^n}
 e^{-t(b_v\eta_1^2/2+a_v|\eta_\perp|^2)}\,d\eta
 =\frac1{\sqrt{2\pi t b_v}\,(4\pi t a_v)^{(n-1)/2}}.
\]
The discarded Gaussian tail is uniformly exponentially small. This
proves the formula with a uniform remainder (in fact $O(t^{-1/5})$).}
\end{proof}

\beP[Optimization formula for the linear spreading speed]
\label{p:spread}
Let $\mathcal D=\diag(d_i)$ with $d_i\ge0$, not all zero, let $M$ be an irreducible Metzler matrix with $s(M)>0$, and consider
\[
 Z_t=\mathcal D\Delta_\zeta Z+MZ,\qquad Z=(Z_1,\ldots,Z_m)^\top,\quad\zeta\in\mathbb R^n.
\]
For $\lambda>0$, let $A(\lambda)=\lambda^2\mathcal D+M$, with right Perron vector $q(\lambda)\gg0$, and let $g(\lambda)=s(A(\lambda))$, $c(\lambda)=g(\lambda)/\lambda$. Then:
\begin{enumerate}
\item The function $g$ is strictly convex, and $c$ has a unique minimizer $\lambda_*>0$. Put $c_*=c(\lambda_*)>0$.
\item The number $c_*$ is the spreading speed of this linear system: for nonzero $Z(0,\cdot)\in C_c(\mathbb R^n;\mathbb R_+^m)$,
\[
 c>c_*\implies\sup_{|\zeta|\ge ct}\|Z(t,\zeta)\|\to0,\qquad
 0<c<c_*\implies\inf_{|\zeta|\le ct}Z_i(t,\zeta)\to+\infty\quad\text{for every }i .
\]
\item For every unit vector $e$ and every $\lambda>0$, $Z=e^{-\lambda(e\cdot\zeta-c(\lambda)t)}q(\lambda)$ is a solution: an exponential traveling mode normal to $e$ with speed $c(\lambda)$, unbounded behind its leading edge, so not a front joining two equilibria. The speeds of such modes are exactly the $c\ge c_*$.
\item If $p_*,q_*$ are the left and right Perron vectors of $A(\lambda_*)$, with $\langle p_*,q_*\rangle=1$, then
\[
 c_*=2\lambda_*\langle p_*,\mathcal Dq_*\rangle.
\]
\end{enumerate}
\eeP

\beR[Next-generation reading, and direction]
\label{r:splitDir}
\rev{Proposition~\ref{p:spread} is the linear speed of Weinberger, Lewis and Li \cite{WeinbergerLewisLi2002}, restated here on a siphon face; items 1, 3 and 4 are theirs, and the proof of item 2 is included because the uniformity used later is not quoted elsewhere.} 1. The speed of exponential modes depends in principle on $e$, but since the Laplacian
models diffusion independent of the direction, here it does not depend on $e$ (if each species $i=1,\ldots,m$ diffuses as $\nabla\cdot(\mathcal A_i\nabla Z_i)$, with $\mathcal A_i$ a constant symmetric positive definite $n\times n$ matrix, then the $m\times m$ diagonal matrix with entries $e^\top\mathcal A_1e,\ldots,e^\top\mathcal A_me$ replaces $\mathcal D$).

2. For any regular splitting $M=F-V$, $s(M)>0\iff R_0=\rho(FV^{-1})>1$.
Thus, $c_*>0$ exactly when $R_0>1$; for $s(M)<0$ the infimum is $-\infty$, and for $s(M)=0$ it is zero. For $m=1$ and $M=r$ the proposition gives the famous Fisher--KPP speed $c_*=2\sqrt{\mathcal Dr}$,
and if $r=\mu(R_0-1)$ then $c_*=2\sqrt{\mathcal D\mu(R_0-1)}$. The review \cite{DavydovychDutkaCherniha2023} surveys this dichotomy for epidemic models, from the plague and rabies models of Noble and of K\"all\'en and Murray to spatial COVID-19 models simulated in the plane with data from Lombardy.

\eeR

\begin{proof}
(1) The Perron root of a Metzler matrix is a nondecreasing convex function of its diagonal entries \rev{(Cohen \cite{Cohen81}; see also \cite{Nussbaum1986})}, so $G(\mu)=s(M+\mu\mathcal D)$ is convex. Since $M$ is irreducible, $G$ is a simple eigenvalue, hence
smooth, with $G'=\langle p,\mathcal Dq\rangle>0$, because $p,q\gg0$ and $\mathcal D\ne0$. Thus $g(\lambda)=G(\lambda^2)$ has $g'=2\lambda G'$ and $g''=2G'+4\lambda^2G''\ge2G'(0)>0$. Since $g\to s(M)>0$ as $\lambda\to0^+$ and
$g\sim\lambda^2\max_id_i$, $c\to+\infty$ at both ends. A critical point of $c$ solves $h=0$ with $h=\lambda g'-g$, and $h'=\lambda g''>0$, so it is unique and is the minimizer.

(4) $c'(\lambda_*)=0$ reads $\lambda_*g'(\lambda_*)=g(\lambda_*)$.

(3) The first statement follows by substitution, and the second because $c$ is continuous with minimum $c_*$ and tends to $+\infty$.

\rev{(2) For $0\ne Z(0,\cdot)\in C_c(\mathbb R^n;\mathbb R_+^m)$, choose
$0\ne\phi\in C_c^\infty(\mathbb R^n;\mathbb R_+^m)$ with
$0\le\phi\le Z(0,\cdot)$. Positivity of the solution operator gives
$Z(t,\zeta)\ge Z_\phi(t,\zeta)$. Apply Lemma~\ref{l:trotter}(2)
to $Z_\phi$. Since $\psi'(v)=-\lambda_v$ and $g$ is even and
strictly convex, $\psi$ is even and strictly decreasing for $v>0$;
$\psi(c_*)=0$ because $g(\lambda_*)=c_*\lambda_*$ and
$g'(\lambda_*)=c_*$. Hence $\min_{|v|\le c}\psi(v)>0$ for
$0<c<c_*$. The leading coefficient in the lemma is componentwise
positive and has a positive lower bound when $|v|\le c$ and $|e|=1$.
Its uniform asymptotic proves the interior assertion. The exterior
assertion follows by comparison with the solution
$Kq(\lambda)e^{-\lambda(e\cdot\zeta-c(\lambda)t)}$, where
$e=\zeta/|\zeta|$, $c(\lambda)<c$, and $K$ can be chosen independently
of $e$ because the initial data have compact support.}
\end{proof}
The formula for the derivative of a Perron eigenvalue used in the proof is classical \cite{Seneta,Nussbaum1986}; it also underlies the reduction principle and the evolution of dispersal \cite{Karlin1982,Altenberg2012,Dockery1998}.

\beP[A moving-frame next-generation test and a drift threshold][moving-frame test and drift][drift]
\label{r:movNGM}
With the splitting $M=F-V$ of Remark~\ref{r:splitDir}, the speed can be read as a family of next-generation tests.
\begin{enumerate}
\item \emph{A moving observer; no transport is added to the model.} With $\xi=e\cdot\zeta-ct$, a perturbation $e^{-\lambda\xi}w(t)$ of the system of Proposition~\ref{p:spread} satisfies $w_t=(\lambda^2\mathcal D+F-V-c\lambda I)w$, so it decays, as seen by an observer moving at speed $c$, exactly when $s(\lambda^2\mathcal D+F-V-c\lambda I)<0$. When $B_c(\lambda)=V+c\lambda I-\lambda^2\mathcal D$ is a nonsingular $M$-matrix, the regular splitting $F-B_c(\lambda)$ turns this into $R_c(\lambda)=\rho(FB_c(\lambda)^{-1})<1$; otherwise $s(\lambda^2\mathcal D+F-V-c\lambda I)\ge0$. Hence
\[
 c>c_*\iff R_c(\lambda)<1\ \text{ for some admissible }\lambda>0:
\]
a speed exceeds $c_*$ exactly when an exponentially decaying bound of the leading edge exists in that moving frame. The value $R_c(\lambda)$ depends on the mode, so it serves as a test, with $R_0$ its value at $\lambda=0$ for every $c$.
\item \emph{Near the threshold.} For a smooth family $M\to M_0$, with $M_0$ irreducible Metzler, $s(M_0)=0$, and $\mathcal D\ne0$ fixed, $c_*=2\sqrt{\langle\mathcal D\rangle\,s(M)}\,(1+o(1))$, where $\langle\mathcal D\rangle=p^\top\mathcal Dq/p^\top q$ with $p,q$ the Perron vectors of $M_0$.
\item \emph{Physical drift.} Consider $Z_t=dZ_{\zeta\zeta}+qZ_\zeta+(b-v)Z$ with $d>0$, $b>v>0$, $q\ge0$: the transport is toward decreasing $\zeta$, $R_0=b/v$, the speed against the flow is $c_{\rm up}=2\sqrt{d(b-v)}-q$, and
\[
 R_{\rm abs}=\min_{\substack{\lambda\ge0\\v+q\lambda-d\lambda^2>0}}\frac b{v+q\lambda-d\lambda^2}=\frac b{v+q^2/(4d)},\qquad R_{\rm abs}>1\iff q<2\sqrt{d(b-v)}.
\]
In the range $R_{\rm abs}<1<R_0$ the perturbation grows but cannot spread upstream \cite{SpeirsGurney2001}. In general, for $Z_t=\mathcal D\Delta_\zeta Z+Q\,\nu\cdot\nabla_\zeta Z+MZ$, with $Q\ge0$ diagonal and $\nu$ a fixed unit vector, the matrix $\lambda^2\mathcal D-\lambda(e\cdot\nu)Q+M$ depends on $e$. For $e=\nu$, let $c_{\rm up}=\inf_{\lambda>0}s(\lambda^2\mathcal D-\lambda Q+M)/\lambda$ and let $R_{\rm abs}$ be the infimum of $\rho(F(V-\lambda^2\mathcal D+\lambda Q)^{-1})$ over $\lambda\ge0$ with $V-\lambda^2\mathcal D+\lambda Q$ a nonsingular $M$-matrix. Then $R_{\rm abs}>1$ implies $c_{\rm up}>0$, and $c_{\rm up}>0$ implies $R_{\rm abs}\ge1$, strictly when the infimum is attained. For $Q=qI$, $c_{\rm up}=c_*-q$.
\end{enumerate}
All this concerns the linear speed: when growth is stronger at higher density, a front can advance although $R_0<1$.
\eeP

\beR[Linear invasion speeds on the siphon lattice]
\label{ss:kpp-siphon}
Let $P_\Sigma$ be a spatially constant resident equilibrium on a siphon face, with irreducible transverse matrix $M=F_\Sigma(P_\Sigma)-V_\Sigma(P_\Sigma)$, $s(M)>0$ and $\mathcal D_\Sigma\ne0$. Proposition~\ref{p:spread} then gives the linear spreading speed
\[
 c_\Sigma(P_\Sigma)=\inf_{\lambda>0}\frac1\lambda\,s\bigl(\lambda^2\mathcal D_\Sigma+F_\Sigma(P_\Sigma)-V_\Sigma(P_\Sigma)\bigr).
\]
Reducible blocks and spatially varying residents require further analysis.
In multistrain systems this gives $c_1,c_2$ for the linearized invasion of the disease-free equilibrium by each strain, and $c_{2|1},c_{1|2}$ for the linearized invasion of the resident of the other strain.
Identifying these with the speeds of the nonlinear system requires linear-determinacy hypotheses \cite{WeinbergerLewisLi2002}; for partly diffusing cooperative systems see \cite{Li2012}.
\eeR

\subsection{Application to a two-strain cholera model}

\label{ss:kpp-two-strain}

\beF[Species-reaction graph of a two-strain cholera model]
\centering
\begin{tikzpicture}[>=stealth, scale=1.0,
  species/.style={rectangle, draw, thick,
                  minimum width=1.1cm, minimum height=0.65cm},
  rxn/.style   ={circle, draw, thick, fill=white,
                  font=\small, inner sep=3pt},
  sp2rx/.style ={->, dotted, thick},
  rx2sp/.style ={->, dashed, thick}
]

\node[anchor=north west, inner sep=0pt] at (0, 0.3) {%
  \parbox{6.3cm}{%
    \[
    \left\{
    \begin{aligned}
    S_t      &= \Lambda-\mu S-\beta_1SB_1-\beta_2SB_2,\\
    (I_1)_t  &= \beta_1SB_1-(\mu+\gamma_1)I_1,\\
    (I_2)_t  &= \beta_2SB_2-(\mu+\gamma_2)I_2,\\
    (B_1)_t  &= d_1\Delta_\zeta B_1+\xi_1I_1-\delta_1B_1,\\
    (B_2)_t  &= d_2\Delta_\zeta B_2+\xi_2I_2-\delta_2B_2.
    \end{aligned}
    \right.
    \]%
  }%
};

\begin{scope}[xshift=8.2cm, scale=0.90]

\node[species] (S)  at (0,    0)   {$S$};
\node[species] (I1) at (3,  1.8)   {$I_1$};
\node[species] (I2) at (3, -1.8)   {$I_2$};
\node[species] (B1) at (6,  1.8)   {$B_1$};
\node[species] (B2) at (6, -1.8)   {$B_2$};

\node[rxn] (R3) at (1.5,  1.8) {\scriptsize$\beta_1 SB_1$};
\node[rxn] (R4) at (1.5, -1.8) {\scriptsize$\beta_2 SB_2$};

\node[rxn] (R7) at (4.5,  2.7) {\scriptsize$\xi_1 I_1$};
\node[rxn] (R8) at (4.5, -2.7) {\scriptsize$\xi_2 I_2$};

\draw[->] (-1.6, 0)    -- (S)            node[midway, above] {$\Lambda$};
\draw[->] (S.north)    -- ++(0,+0.9)     node[above]         {$\mu S$};
\draw[->] (I1.north)   -- ++(0,+0.9)     node[above]         {$(\mu{+}\gamma_1)I_1$};
\draw[->] (I2.south)   -- ++(0,-0.9)     node[below]         {$(\mu{+}\gamma_2)I_2$};
\draw[->] (B1.north)   -- ++(0,+0.9)     node[above]         {$\delta_1 B_1$};
\draw[->] (B2.south)   -- ++(0,-0.9)     node[below]         {$\delta_2 B_2$};

\draw[sp2rx] (S.north east) to[bend left=20]  (R3);
\draw[sp2rx] (B1.west)                       -- (R3);
\draw[rx2sp] (R3)                             -- (I1);

\draw[sp2rx] (S.south east) to[bend right=20] (R4);
\draw[sp2rx] (B2.west)                       -- (R4);
\draw[rx2sp] (R4)                             -- (I2);

\draw[sp2rx] (I1.north) to[bend left=15]   (R7);
\draw[rx2sp] (R7)       to[bend left=15]   (I1.east);
\draw[rx2sp] (R7)       to[bend right=15]  (B1.north);

\draw[sp2rx] (I2.south) to[bend right=15]  (R8);
\draw[rx2sp] (R8)       to[bend right=15]  (I2.east);
\draw[rx2sp] (R8)       to[bend left=15]   (B2.south);

\end{scope}
\end{tikzpicture}

\vspace{0.8em}
\renewcommand{\arraystretch}{1.2}
\begin{tabular}{rllp{5.2cm}}
\toprule
No. & Reaction & Rate & Interpretation \\
\midrule
1  & $\varnothing\to S$       & $\Lambda$               & Host recruitment \\
2  & $S\to\varnothing$        & $\mu S$                 & Natural death of $S$ \\
3  & $S+B_1\to I_1+B_1$      & $\beta_1 SB_1$          & Strain-1 infection ($B_1$ catalytic) \\
4  & $S+B_2\to I_2+B_2$      & $\beta_2 SB_2$          & Strain-2 infection ($B_2$ catalytic) \\
5  & $I_1\to\varnothing$      & $(\mu+\gamma_1)I_1$     & Removal of infected host $I_1$ \\
6  & $I_2\to\varnothing$      & $(\mu+\gamma_2)I_2$     & Removal of infected host $I_2$ \\
7  & $I_1\to I_1+B_1$        & $\xi_1 I_1$             & Pathogen shedding by $I_1$ ($I_1$ catalytic) \\
8  & $I_2\to I_2+B_2$        & $\xi_2 I_2$             & Pathogen shedding by $I_2$ ($I_2$ catalytic) \\
9  & $B_1\to\varnothing$      & $\delta_1 B_1$          & Environmental decay of $B_1$ \\
10 & $B_2\to\varnothing$      & $\delta_2 B_2$          & Environmental decay of $B_2$ \\
\bottomrule
\end{tabular}
\lbl{fig:twostrain-env}
\caption{Two-strain environmental transmission model.
  \emph{Left}: ODE/PDE system (\textbf{diffusing}: $B_1, B_2$;
  \textbf{non-diffusing}: $S, I_1, I_2$).
  Disease siphons: $\sigma_k=\{I_k,B_k\}$, $k=1,2$,
  and the joint siphon $\sigma_{12}=\{I_1,B_1,I_2,B_2\}$.
  \emph{Right}: species-reaction (Volpert/DSR hybrid) graph.
  All multi-species reactions are shown as bipartite nodes.
  Infection reactions~3--4 ($S+B_k\to I_k+B_k$): dotted arrows from
  $S$ and $B_k$ into the node, dashed arrow out to $I_k$.
  Shedding reactions~7--8 ($I_k\to I_k+B_k$): dotted arrow from $I_k$
  into the node, dashed arrows out to $I_k$ (returned) and $B_k$
  (produced); nodes placed above/below the species row to avoid overlap.
  \textbf{Dotted}: reactant $\to$ node;
  \textbf{dashed}: node $\to$ product.}
\eeF
The disease-free equilibrium is
$$
X_0=
\left(
\frac{\Lambda}{\mu},
0,0,0,0
\right).
$$

The strain-$1$ transversal subsystem at the DFE is
$$
Z_1=(I_1,B_1)^T,
$$
with
$$
(Z_1)_t
=
\mathcal D_1\Delta_\zeta Z_1
+
\T M_1Z_1,
$$
where
$$
\mathcal D_1=
\bep
0&0\\
0&d_1
\eep
$$
and
$$
\T M_1=
\bep
-(\mu+\gamma_1)
&
\dfrac{\beta_1\Lambda}{\mu}
\\[6pt]
\xi_1
&
-\delta_1
\eep .
$$

The next-generation matrix is
$$
F_1V_1^{-1}
=
\bep
0
&
\dfrac{\beta_1\Lambda}{\mu\delta_1}
\\[6pt]
\dfrac{\xi_1}{\mu+\gamma_1}
&
0
\eep .
$$

Hence
$$
R_1
=
\sqrt{
\frac{
\beta_1\xi_1\Lambda
}{
\mu(\mu+\gamma_1)\delta_1
}
}.
$$

The spreading speed of strain $1$ into the DFE is therefore
$$
c_1
=
\inf_{\lambda>0}
\frac1\lambda
s(\lambda^2\mathcal D_1+\T M_1).
$$

Now
$$
\lambda^2\mathcal D_1+\T M_1
=
\bep
-(\mu+\gamma_1)
&
\dfrac{\beta_1\Lambda}{\mu}
\\[6pt]
\xi_1
&
d_1\lambda^2-\delta_1
\eep .
$$

Hence
$$
s(\lambda^2\mathcal D_1+\T M_1)
=
\frac12
\left[
d_1\lambda^2-(\mu+\gamma_1+\delta_1)
+
\sqrt{
(d_1\lambda^2+\mu+\gamma_1-\delta_1)^2
+
4\dfrac{\beta_1\xi_1\Lambda}{\mu}
}
\right].
$$

Thus
$$
R_1>1
\implies
c_1>0.
$$

Exactly the same formula holds for strain $2$.

\subsection{Application to transversal reinvasion}
\label{ss:kpp-reinvasion}

Assume now that strain $1$ has established the semitrivial equilibrium
$$
X_1^*
=
(S_1^*,I_1^*,0,B_1^*,0).
$$

The transversal subsystem for strain $2$ is
$$
Z_{2|1}=(I_2,B_2)^T,
$$
with
$$
(Z_{2|1})_t
=
\mathcal D_2\Delta_\zeta Z_{2|1}
+
\T M_{2|1}Z_{2|1},
$$
where
$$
\mathcal D_2=
\bep
0&0\\
0&d_2
\eep
$$
and
$$
\T M_{2|1}
=
\bep
-(\mu+\gamma_2)
&
\beta_2S_1^*
\\[6pt]
\xi_2
&
-\delta_2
\eep .
$$

The transversal reproduction number is
$$
R_{2|1}
=
\sqrt{
\frac{
\beta_2\xi_2S_1^*
}{
(\mu+\gamma_2)\delta_2
}
}.
$$

The transversal spreading speed becomes
$$
c_{2|1}
=
\inf_{\lambda>0}
\frac1\lambda
s(\lambda^2\mathcal D_2+\T M_{2|1}).
$$

Hence:
$$
R_{2|1}>1
\Longrightarrow
c_{2|1}>0.
$$

Thus the linearized transverse equation has positive spreading speed; since $S$ makes the system non-cooperative, identifying this speed with that of the nonlinear replacement requires further work. \rev{The two are numerically consistent in Figure~\ref{fig:choFront}(h); linear determinacy is not proved here.}
Since $S_1^*=(\mu+\gamma_1)\delta_1/(\beta_1\xi_1)$,
\[
 R_{2|1}=\frac{R_2}{R_1},\qquad R_{1|2}=\frac{R_1}{R_2},\qquad R_{2|1}R_{1|2}=1,
\]
so the two strains cannot both invade the resident equilibrium of the other. \rev{This identity uses the hypotheses of this subsection and no others: a spatially constant resident, one susceptible compartment with its constant value $S_1^*$, and no cross-immunity, so that each strain sees the other only through $S$.} Speed into unoccupied territory and invasion of the resident are different competitive advantages: in Figure~\ref{fig:choFront}, strain $1$ has $c_1>c_2$, while strain $2$ has $R_{2|1}>1$.

\begin{figure}[H]
\centering
\includegraphics[width=\linewidth]{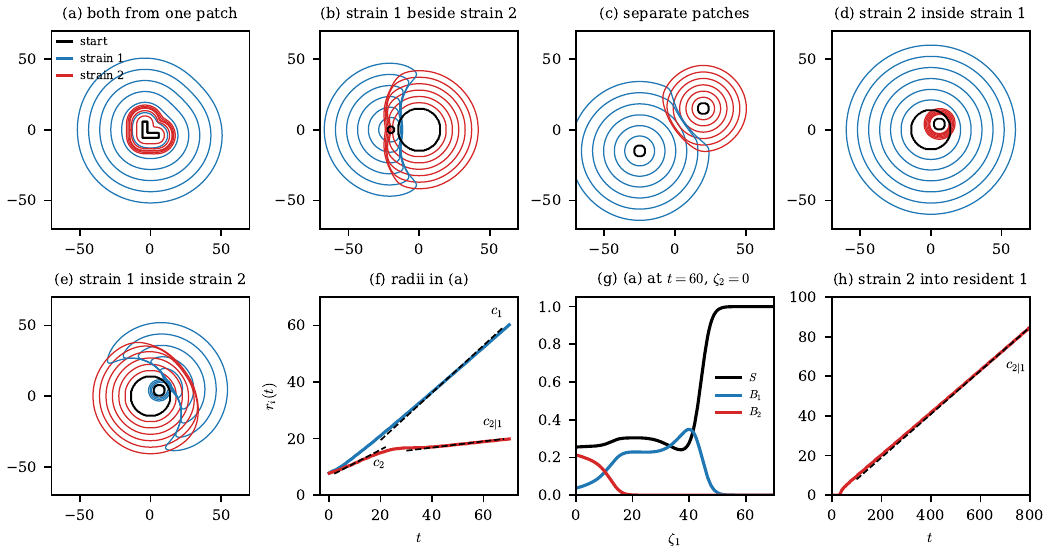}
\caption{The two-strain cholera model with spatial diffusion $d_1=1$, $d_2=1/4$, in the plane, with illustrative parameters: $\Lambda=\mu=0.1$, $\beta_1=1$, $\beta_2=1.2$, $\xi_i=\delta_i=1$, $\gamma_i=0.2$,  so $R_1\approx1.83$, $R_2=2$, $R_{2|1}\approx1.10$.
Panels (a)--(e): the fronts, meaning the level curves $B_1=0.01$ (blue) and $B_2=0.01$ (red), at $t=10,20,\ldots,60$; the black curves are the starting patches. (a) Both strains start from the same $L$-shaped patch. (b) Strain $1$ starts from a small disc beside a large disc of strain $2$. (c) The strains start from separate small discs, and each spreads into the disease-free state until they meet. (d) Strain $2$ starts from a small disc inside a large disc of strain $1$, and stays behind it. (e) Strain $1$ starts from a small disc inside a large disc of strain $2$, and escapes.
Panel (f): in (a), the front radius $r_i(t)=\sup\{|\zeta|:B_i(t,\zeta)\ge0.01\}$; the dashed lines have the slopes of the linear speeds $c=\inf_{\lambda>0}s(\lambda^2\mathcal D+M)/\lambda$, namely $c_1\approx0.84$, $c_2\approx0.47$, $c_{2|1}\approx0.11$. Strain $2$ moves into the disease-free state, at a slope approaching $c_2$, until strain $1$, which diffuses faster although it has the smaller reproduction number, has passed and lowered $S$; near $t=25$ the slope changes to about $c_{2|1}$, since strain $2$ now invades the resident equilibrium of strain $1$. Panel (g): the profiles in (a) at $t=60$ along $\zeta_2=0$, $\zeta_1\ge0$, show this: ahead $S=1$, behind the front of strain $1$ $S$ settles at $S_1^*=0.3$, and behind that of strain $2$ it falls toward $S_2^*=0.25$. Panel (h): the invasion in its pure form, on a half line, with the resident of strain $1$ everywhere and a small localized $B_2$; the position of $B_2=0.01$ has slope $0.1076$ on $400<t<800$, against the linear $c_{2|1}=0.1097$, and halving the mesh changes this by less than $0.001$.
Panels (a)--(e): $\zeta\in[-75,75]^2$, no-flux boundary, mesh $0.5$; (h): $\zeta_1\in[0,160]$, mesh $0.25$; explicit Euler.}
\label{fig:choFront}
\end{figure}

\subsection{Computational realization by conic optimization}
\label{s:spread-conic}

\beT[Convex feasibility certificates for the speed][][SpdQC]
For the linear cooperative invasion equation of Proposition~\ref{p:spread}, for any $c>0$,
\[
c_*\le c
\qquad\Longleftrightarrow\qquad
\exists\,\lambda>0,\ w\gg0:\quad
\bigl(M+\lambda^2\mathcal D-c\lambda I\bigr)w\le0 .
\]
With $w_i=e^{u_i}$ and $u_1=0$, this is the convex feasibility problem in
$(\lambda,u)$
\[
M_{ii}+d_i\lambda^2-c\lambda+\sum_{j\ne i}M_{ij}e^{u_j-u_i}\le0,
\qquad i=1,\ldots,m ,
\]
and each accepted trial value comes with a certificate $(\lambda,w)$
verifiable by $m$ inner products.
\eeT
\begin{proof}
$c_*\le c$ iff $g(\lambda)\le c\lambda$ for some $\lambda>0$, i.e.\
$s(M+\lambda^2\mathcal D)\le c\lambda$, i.e.\
$s\bigl(M+\lambda^2\mathcal D-c\lambda I\bigr)\le0$, the shift by
$-c\lambda I$ moving the spectrum rigidly. The shifted matrix is Metzler,
and for irreducible Metzler matrices $s(A)\le0$ is equivalent to the
existence of $w\gg0$ with $Aw\le0$ by Collatz--Wielandt. Dividing the $i$th
inequality by $w_i=e^{u_i}>0$ gives the displayed constraints, which are
convex in $(\lambda,u)$ because $M_{ij}\ge0$ for $j\ne i$; the normalization
$u_1=0$ removes the scaling freedom.
\end{proof}

Deciding whether a proposed
speed is achievable is a convex feasibility problem in $(\lambda,u)$, a
linear program in $w$ only once $\lambda$ is fixed, and a positive answer is
witnessed by a pair $(\lambda,w)$ that can be checked without any
eigenvalue computation: one forms $M+\lambda^2\mathcal D-c\lambda I$,
multiplies by $w$, and checks $m$ signs. Bisecting on $c$ then converges
to $c_*$; the problem is convex for fixed $c$ only, not jointly in $c$. Each verified feasible iterate provides a certified upper bound, whereas a failure of the solver certifies nothing, so a certified lower bound needs an infeasibility certificate.

\medskip\noindent\textbf{Numerical check.} On $300$ random pairs
$(M,\mathcal D)$ with $M$ Metzler, $s(M)>0$, and sizes $2$ to $4$, the
second differences of $g$ on a grid of $400$ values of $\lambda$ were
nonnegative throughout, and every sampled sublevel set of $c$ was an
interval, as Proposition~\ref{p:spread} requires. On four further instances the
bisection of the theorem, which never computes a spectral radius, was
compared with direct minimization of
$c(\lambda)$ over a fine grid; the two agreed to $3\cdot10^{-4}$ or
better, the residual being the grid spacing of the direct method instead of
error in the bisection.

\beR[]
\label{r:spread-conic-scope}
The strong convexity of $g$ in Proposition~\ref{p:spread} is what makes the one-dimensional search
safe, and it is not automatic: it uses the convexity of the Perron root
in the diagonal entries, which is a property of Metzler matrices and
would fail for a general linearization. This is the same structural
hypothesis --- irreducible Metzler transversal block --- that underlies
the invasion certificates elsewhere in this paper, appearing here in the
different guise of making a wave-speed optimization tractable.
\eeR

\beR[What the theorem says, and a conjecture on attaining the bound]
[attaining the speed bound][SpdAtt]
We specialize to $B=I$ in Theorem~\ref{t:SpdBd}(3). The following calculation provides evidence for Conjecture~\ref{cj:spdatt}. The response $\Phi(q)=Lq/(1+Lq)$, the odds polynomial of degree one with $a_1=L$, belongs to $\mathcal C_L$, with $\Phi'(0)=L$, $\Phi''<0$ and $0\le\Phi(q)\le Lq$. Thus the linearized system at the origin is the one defining $\bar c(L)$, and $f(x)=-Vx+\Phi(\beta Wx)$ satisfies $f(x)\le(\beta LW-V)x$. \rev{If $\beta LW-V$ is irreducible, which no zero row of $W$ does not by itself give,} $f$ is cooperative and irreducible; if $W$ has no zero row, it is strictly subhomogeneous. Under this additional assumption, with bounded production and an unstable origin, it has exactly one positive equilibrium, and no other equilibrium in the corresponding order interval except the origin. These are standard ingredients of linear-determinacy arguments \cite{WeinbergerLewisLi2002,Li2012}, but a direct application of the theorem of \cite{Li2012} additionally requires verification of its clipped-vector inequality: since $\beta LW-V$ is Metzler without being entrywise nonnegative, $f(x)\le(\beta LW-V)x$ does not by itself give the clipped inequality $f(\min\{\tau\xi,\bar x\})\le(\beta LW-V)\tau\xi$, $\tau>0$, of \cite{Li2012}, with $\xi$ the relevant positive eigenvector and the minimum taken componentwise.

\rev{Proposition~\ref{p:spaCrit} gives an exact criterion for this inequality, and shows where it fails.}

As a numerical illustration, take $\mathcal D=\operatorname{Diag}(1,1/10)$, $V=\begin{pmatrix}1&0\\-1/2&7/10\end{pmatrix}$,
$W=\begin{pmatrix}1/5&9/10\\3/10&1/10\end{pmatrix}$, $\beta=1$, so $\bar R_0\approx1.267$. An explicit finite-difference integration of planar fronts, which depend on $\zeta$ through one coordinate, to $t=120$ gives the speeds below;
we did not separate the part of the gap to $\bar c(L)$ due to finite time, to the grid and to the domain.
\begin{center}
\begin{tabular}{llll}
\hline
Response & $L$ & front speed & $\bar c(L)$\\\hline
$Lq/(1+Lq)$ & $2$ & $1.434$ & $1.460$\\
$1-e^{-Lq}$ & $2$ & $1.435$ & $1.460$\\
$Lq/(1+Lq)$ & $4$ & $2.390$ & $2.442$\\
$1-e^{-Lq}$ & $4$ & $2.387$ & $2.442$\\\hline
\end{tabular}
\end{center}
For these data, with $(\xi_*)_1=1$, numerically $\lambda_*\approx1.0598$ and $A\xi_*\approx(0.4236,0.8157)^\top$ for $L=2$, and $\lambda_*\approx1.6251$ and $A\xi_*\approx(1.3281,1.5727)^\top$ for $L=4$, so the clipped inequality holds in both examples.
\eeR

\beCj[The uniform speed bound is attained over the response class]
\label{cj:spdatt}
{\color{revcolor}
In the setting of Theorem~\ref{t:SpdBd}(3) with $B=I$, $\bar c(L)$ is attained over
$\mathcal C_L$: there is a response $\Phi\in\mathcal C_L$ whose fronts spread at
the speed $\bar c(L)$. The evidence is the computation reported in
Remark~\ref{r:SpdAtt}; the obstruction to a proof is identified in
Proposition~\ref{p:spaCrit}.}
\eeCj

\beP[The clipped-vector inequality of the linear-determinacy theorem has an exact criterion][clipped inequality criterion][spaCrit]
{\color{revcolor}
In the setting of Remark~\ref{r:SpdAtt}, with $\bar x$ the positive equilibrium, put $A=\beta LW-V$. For $\xi\gg0$,
\[
 f(\min\{\tau\xi,\bar x\})\le\tau A\xi\ \text{ for every }\tau>0\iff A\xi\ge0.
\]}
\eeP
\begin{proof}
{\color{revcolor}
Indeed, let $u=\min\{\tau\xi,\bar x\}$. If $u_i=\tau\xi_i$, then $f_i(u)\le(Au)_i\le\tau(A\xi)_i$, because $A$ is Metzler and the $i$th coordinates agree. If $u_i=\bar x_i$, cooperativity gives $f_i(u)\le f_i(\bar x)=0$. Conversely, for large $\tau$ all coordinates are clipped, and the inequality reads $0\le\tau A\xi$. For the Perron vector $\xi_*$ at the minimizer $\lambda_*$, $(\lambda_*^2\mathcal D+A)\xi_*=g(\lambda_*)\xi_*$, so the criterion becomes $g(\lambda_*)\ge\lambda_*^2\max_id_i$. It holds automatically for equal diffusivities, since then $A\xi_*=s(A)\xi_*\gg0$. For unequal diffusivities it can fail, so it does not establish the conjecture in general; the other hypotheses of the theorem of \cite{Li2012} also remain to be verified, and linear determinacy by another theorem is not excluded.}
\end{proof}

\beR[Related work]
[diffusion literature][SpdLit]
Adding diffusion to threshold networks is not new. The preservation of contraction under diagonal diffusion is established in \cite{AminzareSontag2012} on bounded domains with homogeneous Neumann boundary conditions, for constant diagonal weights; the comparison argument of item 1 gives the weighted maximum-norm estimate on $\mathbb R^n$. Reaction--diffusion neural networks with activation functions, $M$-matrix criteria and exponential estimates are studied in \cite{SongCaoZhao2006}, and activation-driven front initiation in nonlocal neural fields in \cite{FayeKilpatrick2018}. What is added here is that the diffusion extension preserves the sharp response-class threshold of the constructive multistability theorem, while the siphon decomposition identifies the resident-dependent transverse matrices governing linear invasion.
\eeR

\begin{figure}[H]
\centering
\includegraphics[width=\linewidth]{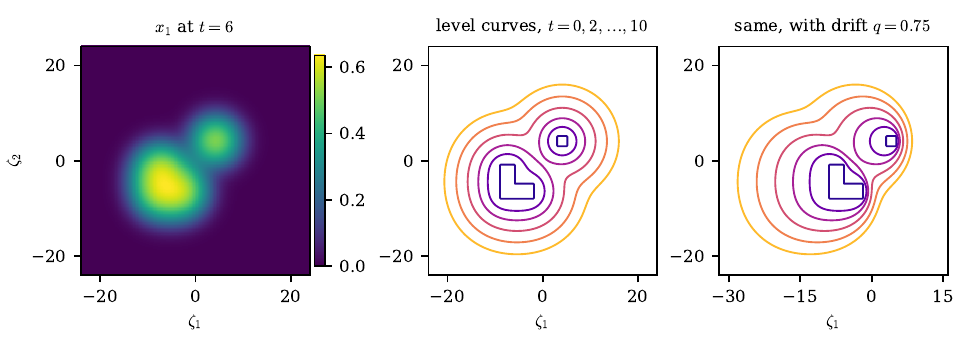}
\caption{Fronts in the plane $\zeta=(\zeta_1,\zeta_2)$, shown without reducing to a distance from one point, by level curves of the state at successive times, for $x_t=\mathcal D\Delta_\zeta x-Vx+\Phi(Wx)$ with $\mathcal D,V,W$ of Remark~\ref{r:SpdAtt} and $\Phi(q)=2q/(1+2q)$, started from an $L$-shaped patch and a separate square.
Left: $x_1$ at $t=6$. Middle: the level curves $x_1=0.1$ at $t=0,2,\ldots,10$; the two fronts merge and the outline rounds, since the irregular initial shape fades relative to the distance travelled and the diffusion has no preferred direction. Right: the same with drift $Q=0.75\,I$
transporting toward decreasing $\zeta_1$, computed from the exact shift $x(t,\zeta)=x_0(t,\zeta+0.75\,te_1)$. At these times the outline advances at about $1.2$--$1.3$ per unit
time along the axes, below the upper bound $\bar c(2)\approx1.46$, conjectured to equal the asymptotic spreading speed.}
\label{fig:spaFront}
\end{figure}
\subsection{Fronts for odds responses and non-disjoint siphons}\label{s:spaOv}
\beC[At a fixed resident the linear speed depends on the odds response only through its onset coefficient]
[odds speed][OddSpd]
Let $\Phi_a$ be an odds-polynomial response, with $h_i=0$ and nonnegative polynomials $P_i$, $P_i(0)=0$, and let $y$ be a spatially constant resident equilibrium on the face of a siphon $\Sigma$.
Put $G(y)=\operatorname{Diag}(b_i)\partial_{x_\Sigma}P_\Sigma(0,y)$, and suppose that $M=a_1G(y)-V_\Sigma$ is irreducible with $s(M)>0$ and that $\mathcal D_\Sigma$ is nonnegative diagonal and nonzero. Then
\[
 c_*(a_1,y)=\inf_{\lambda>0}\frac1\lambda\,s\bigl(\lambda^2\mathcal D_\Sigma-V_\Sigma+a_1G(y)\bigr),\qquad
 \frac{\partial c_*}{\partial a_1}=\frac{\langle p_*,G(y)q_*\rangle}{\lambda_*\langle p_*,q_*\rangle},
\]
with $y$ held fixed in the derivative, where $p_*,q_*$ are the Perron vectors at the minimizer $\lambda_*$. This linear spreading speed is positive exactly when $R_\Sigma=a_1\rho(G(y)V_\Sigma^{-1})>1$. At a fixed resident
state only $a_1$ enters; the higher coefficients act only by moving $y$, and along a resident branch $y=y(a_1)$ the derivative has further terms from the motion of $y$.
\eeC
\begin{proof}
The face stays invariant under diffusion, since $\mathcal D$ is diagonal. Proposition~\ref{p:Odd} gives the transverse matrix, and Proposition~\ref{p:spread} the formula. The minimizer is
unique by Proposition~\ref{p:spread}, so the envelope theorem and the derivative of the Perron eigenvalue, $\partial_{a_1}s=\langle p,Gq\rangle/\langle p,q\rangle$ give the derivative.
\end{proof}

\beR[The double relay gives two speeds, and stable residents can be invaded]
[relay speeds][RelSpd]
For the three-population network of Theorem~\ref{t:KBoverlapCusp} with equal diffusion $d$, the transverse matrices have the common diagonal $a_1-1$, so
$s(\lambda^2dI+J_{12}^\perp)=d\lambda^2+R_{\Sigma_{12}}-1$. Hence the double boundary transcritical relay of Remark~\ref{r:KBrelay} switches on the linear speeds
\[
 c_{12}=2\sqrt{d(R_{\Sigma_{12}}-1)}\quad(\beta>\beta_{12}),\qquad c_{23}=2\sqrt{d(R_{\Sigma_{23}}-1)}\quad(\beta>\beta_{23})
\]
in turn. These need not be the speeds of actual fronts. For the reference coefficients and the couplings of Remark~\ref{r:KBrelay}, with $d=1$, integrating from the upper state of
Figure~\ref{fig:relay} as a planar front into $E_{12}$ gives the following.
\begin{center}
\begin{tabular}{llll}
\hline
$\beta/\beta_{12}$ & $R_{\Sigma_{12}}$ & front speed & $c_{12}$\\\hline
$0.5$ & $0.99915$ & $0.059$ & none\\
$0.99$ & $0.99999$ & $0.074$ & none\\
$1.2$ & $1.00028$ & $0.078$ & $0.033$\\\hline
\end{tabular}
\end{center}
Thus, numerically, a front invades $E_{12}$ while $E_{12}$ is linearly stable, and the spatial transition does not occur at $\beta_{12}$.
\eeR

\beR[The scalar quartic, and an upper bound by the linear speed]
[scalar odds fronts][OddFro]
For planar fronts of $x_t=\Delta_\zeta x-x+\Phi_a(x)$ with the states $0$ and $s_a$ stable and one zero of $f=-x+\Phi_a(x)$ between them, the bistable theory \cite{FifeMcLeod1977} gives a unique front, whose
speed has the sign of $\int_0^{s_a}f$. For the reference coefficients this integral is approximately $+2.5\cdot10^{-4}$, and shooting the profile equation gives speed
approximately $+0.0535$ with $R=a_1\approx0.997$: the stable origin is invaded although $R<1$, an example of the warning in Remark~\ref{r:movNGM}.

Conversely, if $a_k\le a_1a_{k-1}$ for $k\ge2$, then comparing coefficients in $\varphi(q)\le a_1q(1+\varphi(q))$ gives $\Phi_a(q)\le a_1q$. For cooperative systems in which $\Phi_a$ has a linear argument this gives a KPP-type upper bound by the linear speed; it does not prove linear determinacy for nonlinear arguments such as $x_2+\beta x_1x_3$, or for invasion of a nonzero resident. For $k=2$ this is $a_2\le a_1^2$, the opposite of the backward condition of Remark~\ref{r:KBrelay}. It is sufficient but not necessary: with the data of
Remark~\ref{r:SpdAtt} and $(a_1,a_2)=(1,6/5)$ the front speed is $0.550$, close to the linear value $0.572$, whereas $(a_1,a_2)=(1,12/5)$ gives $0.613$.
\eeR

\beR[A field example of the bistable regime: Wolbachia in mosquitoes]
[Wolbachia fronts][Wolb]
The frequency $p$ of the bacterium \emph{Wolbachia} in \emph{Aedes aegypti} is modelled in \cite{TurelliBarton2017} by $p_t=\frac{\sigma^2}2\Delta_\zeta p+s_hp(1-p)(p-\hat p)$, with $\sigma$ the dispersal scale and $\hat p$ the unstable frequency. This is
the regime of Remark~\ref{r:OddFro}. The state $p=0$ is linearly stable, so the formula of Proposition~\ref{p:spread} gives no speed, yet for $\hat p<1/2$ a front advances, with the single speed $c=\sigma\sqrt{s_h}(\frac12-\hat p)$; its sign is that of
$\int_0^1s_hp(1-p)(p-\hat p)\,dp$. The same speed describes the asymptotic radial expansion in the plane. Unlike the linear case, a compact start must be large enough: below a critical radius, about $3.5\sigma$ at $\hat p=0.3$, a release fails. These predictions were compared with
releases in Cairns \cite{SchmidtEtAl2017}: the predicted speed is of order $150$--$300$ m per year, the two larger releases spread at about $100$--$200$ m per year, and the smallest showed little spread after two years. For a two-component model of infected and uninfected mosquitoes, \cite{ChanKimMarangell2016} find a unique wave speed, numerically, and show by an Evans function computation that the wave is linearly stable; the fronts in our figures are discussed in Section~\ref{s:spaSta}.
\eeR

\beR[Where two residents meet, the shared species ignites]
[ignition at contact][OvlIgn]
Let the residents $E_{23}$ and $E_{12}$ occupy two territories with a common border, and let $x_2=0$ initially. This state lies on no single siphon face. In each territory the production of the shared species $2$ vanishes, because $x_1x_3=0$; at the
border, diffusion makes $x_1x_3$ positive, so $x_2>0$ at once, and by ignition we mean growth of $x_2$ to a macroscopic level. As a local kinetic estimate, with no threshold proved for the PDE: for fixed $p=x_1x_3$, the small root of $x_2=\Phi_a(x_2+\beta p)$ disappears when $\beta p$ exceeds $\max_{0<x<1/8}(\Phi_a^{-1}(x)-x)\approx3.5\cdot10^{-5}$ for the reference coefficients,
so that with $p\approx1/4$ the estimated threshold is $\beta\approx1.4\cdot10^{-4}$. Below it $x_2$ stays small: at the table coupling $\beta=10^{-6}$, $\max x_2\approx1.7\cdot10^{-6}$ at $t=40$. Above it $x_2$ ignites and approaches one
(Figure~\ref{fig:overFront}). At $\beta=1.2\,\beta_{12}$, across a straight border, the half-width of $\{x_2>1/4\}$ fits $5.55\sqrt t-7$ for $40\le t\le800$, as the diffusive tails of $x_1$ and $x_3$, which decay slowly because $1-a_1$ is small, carry $\beta x_1x_3$ outward (Figure~\ref{fig:overFront}(b,c)).
This is a transient, not a law: the slope falls to $0.083$ on $800<t<1600$ and to $0.066$ on $3200<t<4000$, still decreasing, while the fronts of $x_1$ and $x_3$ pass each other at about $0.0535$, the speed of the scalar quartic of Remark~\ref{r:OddFro}. The eventual speed of $x_2$ remains undetermined. These statements are numerical, at one parameter choice and with a sharp initial border.
\eeR

\begin{figure}[H]
\centering
\includegraphics[width=\linewidth]{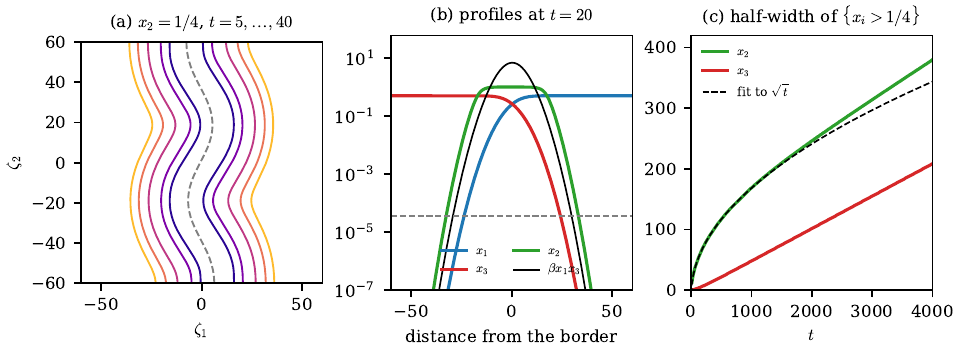}
\caption{The three-population network in the plane, $(x_1)_t=\Delta_\zeta x_1-x_1+\Phi_a(x_1+\alpha x_2)$, $(x_2)_t=\Delta_\zeta x_2-x_2+\Phi_a(x_2+\beta x_1x_3)$, $(x_3)_t=\Delta_\zeta x_3-x_3+\Phi_a(x_3+\gamma x_2)$, with the reference coefficients, the couplings $\alpha,\gamma$ of Table~\ref{tab:KBoverlapNum} and $\beta=1.2\,\beta_{12}$. The residents occupy
the two sides of a border, with $x_1=s_a$ on the right, $x_3=s_a$ on the left and $x_2=0$. (a) A curved border in the plane: the level curves $x_2=1/4$ at $t=5,10,20,30,40$ open on both sides of the level curve $x_1=1/4$ at $t=40$ (dashed), which has hardly moved. (b) A straight border, at $t=20$: the forcing $\beta x_1x_3$ exceeds the local estimate $3.5\cdot10^{-5}$ (dashed) well outside the region where $x_1,x_3$ are both large, and $x_2$ has ignited there. (c) The half-widths of $\{x_2>1/4\}$ and $\{x_3>1/4\}$, measured from the border: $x_2$ first follows the $\sqrt t$ fit, then departs from it, while $x_3$ moves at constant speed. Plane $[-60,60]^2$ to $t=40$, line $[-800,800]$ to $t=4000$; no-flux boundary, mesh $0.5$, explicit Euler with step $0.05$.}
\label{fig:overFront}
\end{figure}

\subsection{What else survives diffusion}\label{s:spaSur}
\beC[Reductions, extinction certificates, box convergence and homogeneous cusps extend to diffusion]
[extensions to diffusion][DifAll]
Let $\mathcal D$ be a positive diagonal matrix. In items 1--3 the spatial domain is $\mathbb R^n$, with bounded solutions, or a bounded domain with homogeneous Neumann or periodic conditions.
\begin{enumerate}
\item \emph{Reductions.} If $f(Hy;\theta)=H\widehat f(y;\theta)$ as in Theorem~\ref{t:Spe} and $\mathcal DH=H\widehat{\mathcal D}$, that is, $d_i=\widehat d_a$ for $i\in C_a$, then $y_t=\widehat{\mathcal D}\Delta_\zeta y+\widehat f(y;\theta)$ gives the solution $x=Hy$ of $x_t=\mathcal D\Delta_\zeta x+f(x;\theta)$. If two species of one group have different diffusivities, a spatially varying group-constant state generally leaves this subspace at once.
\item \emph{Extinction.} In Theorem~\ref{t:Cmp}, let $z\ge0$ be bounded with $z_t\le\mathcal D\Delta_\zeta z+(C-V)z$ componentwise. With the same $p,c,w,c'$,
\[
 \sup_\zeta\max_i\frac{z_i(t,\zeta)}{p_i}\le e^{-ct}\sup_\zeta\max_i\frac{z_i(0,\zeta)}{p_i},
\]
and on a bounded domain $\Omega$ the function $\mathcal L(t)=\int_\Omega w^\top z\,d\zeta$ satisfies $\dot{\mathcal L}\le-c'\mathcal L$. No cooperativity of the resident equations is needed.
\item \emph{Global convergence on a box.} In Theorem~\ref{t:EqWgt}, let $K=[a,b]$ and let its equivalent conditions hold. Every solution with values in $K$ converges to $e$ uniformly and exponentially, for every $\mathcal D$.
\item \emph{Homogeneous saddle-nodes and cusps.} The equilibrium equations of constant states do not change, so the constructed saddle-nodes and cusps stay at the same parameter values. On a bounded connected domain with Neumann conditions, if $f$ is smooth and cooperative, and $J$ has a simple zero eigenvalue with $s(J)=0$ and all other eigenvalues of negative real part, then only the spatially constant mode is critical, and the saddle-node or cusp of the constant states is one of the reaction--diffusion system, with the same center dynamics.
\end{enumerate}
\eeC
\begin{proof}
(1) $Hy_t=H\widehat{\mathcal D}\Delta_\zeta y+H\widehat f(y)=\mathcal D\Delta_\zeta(Hy)+f(Hy)$.

(2) $(C-V)p=-\mathbf1\le-cp$, and comparison as in the proof of Corollary~\ref{c:DifSta} gives the first estimate. Integrating, the diffusion terms vanish, and $\dot{\mathcal L}\le\int_\Omega w^\top(C-V)z\,d\zeta=-\int_\Omega\mathbf1^\top z\,d\zeta\le-c'\mathcal L$.

(3) By comparison with the spatially constant solutions, $\varphi_t(a)\le x(t,\zeta)\le\varphi_t(b)$, and Theorem~\ref{t:EqWgt}(2) gives the exponential convergence of both bounds to $e$.

(4) With $0=\mu_0<\mu_1\le\cdots$ the eigenvalues of $-\Delta_\zeta$, the linearization splits into the blocks $J-\mu_k\mathcal D$. Since these are Metzler and $J-\mu_k\mathcal D\le J-\mu_k\min_id_i\,I$, one has $s(J-\mu_k\mathcal D)\le-\mu_k\min_id_i<0$ for $k\ge1$, a spectral gap. The constants are invariant, so the center manifold of the ODE is that of the reaction--diffusion system.
\end{proof}

\beR[Limits of these extensions]
[limits of the extensions][DifLim]
In item 3, what transfers is the global exponential attraction: the state-dependent contraction weight of Theorem~\ref{t:EqWgt} does not automatically give a contraction metric for the reaction--diffusion system, since differentiation produces terms in the spatial gradients of the weight. Item 4 does not hold in this form on $\mathbb R^n$, where arbitrarily small wave numbers leave no spectral gap. Without cooperativity, arbitrary diagonal diffusion may destabilize the complementary spatial modes; equal diffusion $\mathcal D=dI$ avoids this particular difficulty.
\eeR

\beP[Integrated Lyapunov functions acquire a diffusion term, which has a sign for Volterra functions and gives a spatial basin certificate][][DifLya]
Let $\Omega$ be a bounded domain with smooth boundary, and let $x$ be a classical solution of $x_t=\mathcal D\Delta_\zeta x+f(x)$ with Neumann or periodic conditions, $\mathcal D$ positive diagonal.
\begin{enumerate}
\item If $H$ is $C^2$ on a convex set containing the values of $x$, then
\[
 \frac{d}{dt}\int_\Omega H(x)\,d\zeta=\int_\Omega\nabla H(x)^\top f(x)\,d\zeta-\sum_\alpha\int_\Omega(\partial_{\zeta_\alpha}x)^\top\nabla^2H(x)\,\mathcal D\,(\partial_{\zeta_\alpha}x)\,d\zeta.
\]
The last term is nonpositive when $\nabla^2H\,\mathcal D+\mathcal D\nabla^2H\succeq0$ on that set.
\item For positive solutions and the Lyapunov--Volterra function $L=\sum_ic_i(x_i-x_i^*-x_i^*\log(x_i/x_i^*))$ of Subsection~\ref{s:FenBre},
\[
 \frac{d}{dt}\int_\Omega L(x)\,d\zeta=\int_\Omega\nabla L(x)^\top f(x)\,d\zeta-\sum_ic_id_ix_i^*\int_\Omega\frac{|\nabla_\zeta x_i|^2}{x_i^2}\,d\zeta.
\]
\item Let $f$ be the right side of \eqref{eq:KBoverlapCusp}, under the hypotheses and with the notation of Proposition~\ref{p:ComBas}, and put $\sigma=\eta/(2\sqrt{\alpha\beta})$ and
\[
 R=K\cap\{x_1\le\sigma\sqrt\alpha,\ x_2\le\sigma\sqrt\beta\}.
\]
Then $R$ is invariant for every $\mathcal D$, and $K$ is invariant if $d_1=d_2$. Put
\[
 a_*=\min_{y\in[y_-,y_+]}\frac{y-\Phi_a(y)}{y-s_a},\qquad
 b_*=\frac{\gamma M_\Phi}{\sqrt\alpha}\max_{y\in[y_-,y_+]}|y-s_a|,
\]
with value $1-\Phi_a'(s_a)$ of the quotient at $s_a$; then $a_*>0$. For a solution with values in $K$ and $\kappa>b_*/c$, the function
\[
 \mathcal H(t)=\int_\Omega\Bigl[\tfrac12(x_3-s_a)^2+\kappa\bigl(\sqrt\beta\,x_1+\sqrt\alpha\,x_2\bigr)\Bigr]d\zeta
\]
satisfies
\[
 \dot{\mathcal H}\le-a_*\int_\Omega(x_3-s_a)^2d\zeta-(\kappa c-b_*)\int_\Omega\bigl(\sqrt\beta\,x_1+\sqrt\alpha\,x_2\bigr)d\zeta-d_3\int_\Omega|\nabla_\zeta x_3|^2d\zeta,
\]
so $\mathcal H(t)\le e^{-\omega t}\mathcal H(0)$ with $\omega=\min\{2a_*,c-b_*/\kappa\}$.
\item The Lyapunov--Volterra function $H_V$ of Proposition~\ref{p:LVboundary} applies to the same solutions: for $\kappa>b_V/c$, the function $\mathcal H_V(t)=\int_\Omega H_V(x(t,\zeta))\,d\zeta$ satisfies
\[
 \dot{\mathcal H}_V=\int_\Omega\nabla H_V(x)^\top f(x)\,d\zeta-d_3s_a\int_\Omega\frac{|\nabla_\zeta x_3|^2}{x_3^2}\,d\zeta\le-\omega_V\mathcal H_V,
\]
so $\mathcal H_V(t)\le e^{-\omega_Vt}\mathcal H_V(0)$.
\end{enumerate}
\eeP

\beR[Which energies are compatible with diffusion]
[energies and diffusion][DifEne]
Item 3 replaces the storage $\mathcal E_a-\mathcal E_a(s_a)$ of Proposition~\ref{p:ComBas} by a quadratic one, because $\mathcal E_a''(y)=1/\Phi_a'(\Phi_a^{-1}(y))-1$ changes sign\rev{: at the values following that proposition the minimum of $\mathcal E_a''$ on $[1/5,7/10]$ is about $-0.0099$, near $y=0.33$. Item 4 uses the Lyapunov--Volterra term, which is convex.} At these values $a_*\approx2.0\cdot10^{-4}$, $M_\Phi\approx1.01$, $b_*\approx3.4\cdot10^{-5}$, $b_*/c\approx2.07$, and $R$ has sides $x_1\le0.0224$, $x_2\le0.05$; the rate $\omega$ is at most $c\approx1.6\cdot10^{-5}$.

For the scalar equation $x_t=d\Delta_\zeta x-x+\Phi_a(x)$ the appropriate functional contains the gradient: with $G(x)=\int_0^x(s-\Phi_a(s))\,ds$,
\[
 \frac{d}{dt}\int_\Omega\Bigl[\frac d2|\nabla_\zeta x|^2+G(x)\Bigr]d\zeta=-\int_\Omega x_t^2\,d\zeta.
\]

For the energy of Proposition~\ref{p:SymEne}, $\nabla^2\mathcal E=\operatorname{Diag}(h''(x_i))-\beta W$. The sign condition of item 1 fails for $\mathcal D=dI$ wherever this Hessian has a negative eigenvalue, and for other $\mathcal D$ it can fail even where the Hessian is positive definite. This does not show that $\int_\Omega\mathcal E$ increases along any solution; we simply have no spatial energy for these networks.

Convexity of the generator does not settle the question. In $x$, $D_h(x,p)$ has second derivative $h''(x)\ge0$, while $D_h(p,x)$, the order used in $L_\Phi$ of Subsection~\ref{s:FenBre}, has second derivative $h''(x)+(x-p)h'''(x)$. For $h=x\log x-x$ this equals $p/x^2\ge0$, which is item 2; for other responses it has to be checked.
\eeR
\begin{proof}
(1) Differentiate under the integral and integrate by parts; the boundary terms vanish.

(2) $\nabla^2L=\operatorname{Diag}(c_ix_i^*/x_i^2)$; apply item 1.

(3) On the face $x_1=\sigma\sqrt\alpha$ of $R$, the bounds of the proof of Proposition~\ref{p:ComBas} give $-x_1+\Phi_a(x_1+\alpha x_2)\le-\sigma\sqrt\alpha+\rho\sigma(\sqrt\alpha+\alpha\sqrt\beta)=-c\sigma\sqrt\alpha<0$, and likewise $-x_2+\Phi_a(x_2+\beta x_1x_3)\le-c\sigma\sqrt\beta<0$ on $x_2=\sigma\sqrt\beta$. Also $R\subset K$, since $\sqrt\beta x_1+\sqrt\alpha x_2\le2\sigma\sqrt{\alpha\beta}=\eta$, and the other faces are faces of $K$, where that proof gives the direction of the field. The faces of a rectangle are invariant under any diagonal diffusion, by the invariant-region principle \cite{Smoller1994}. If $d_1=d_2$, the normal $\ell$ of the oblique face of $K$ is a left eigenvector of $\mathcal D$, and the same principle applies to $K$.

The sign conditions on $\Phi_a$ and $\Phi_a'(s_a)<1$ give $a_*>0$. By item 1 with the convex $H$, the diffusion terms contribute $-d_3\int_\Omega|\nabla_\zeta x_3|^2$ and nothing from the linear part. For the reaction, $(y-s_a)(\Phi_a(y)-y)\le-a_*(y-s_a)^2$; next $(y-s_a)(\Phi_a(y+\gamma x_2)-\Phi_a(y))\le b_*(\sqrt\beta x_1+\sqrt\alpha x_2)$, using $x_2\le\ell^\top z/\sqrt\alpha$; and $\ell^\top\dot z\le-c\,\ell^\top z$ as in that proof. Finally, $a_*\int_\Omega(x_3-s_a)^2=2a_*\int_\Omega\frac12(x_3-s_a)^2$ and $\kappa c-b_*=\kappa(c-b_*/\kappa)$ give the rate.

(4) Apply item 1 with $\nabla^2H_V=\operatorname{Diag}(0,0,s_a/x_3^2)$, and integrate the pointwise inequality $\nabla H_V^\top f\le-\omega_VH_V$ of Proposition~\ref{p:LVboundary}.
\end{proof}

\subsection{Stability of the fronts}\label{s:spaSta}
The exponential-weight method is classical \cite{Sattinger1976}; the calculations below concern the linear operator ahead of the front.
\beP[Exponential weights stabilize the state ahead exactly above $c_*$]
[weighted spectrum][WeiSpe]
In the setting of Proposition~\ref{p:spread}, let $z=e\cdot\zeta-ct$ and $w=e^{\lambda z}Z$, with $\lambda>0$. Then $w_t=\mathcal Dw_{zz}+(cI-2\lambda\mathcal D)w_z+(\lambda^2\mathcal D-c\lambda I+M)w$, with symbol
\[
 T(k)=-k^2\mathcal D+ik(cI-2\lambda\mathcal D)+\lambda^2\mathcal D-c\lambda I+M,\qquad\text{and}\qquad \sup_ks(T(k))=s(T(0))=\lambda\bigl(c(\lambda)-c\bigr).
\]
Hence the spectrum of this operator on $L^2(\mathbb R;\mathbb C^m)$, for planar perturbations, lies in an open left half-plane for some $\lambda>0$ exactly when $c>c_*$, and at $c=c_*$ the best choice is $\lambda=\lambda_*$, with spectral abscissa $0$. Without the weight, $\lambda=0$, the spectral abscissa is $s(M)>0$. Transverse wave numbers add $-|k_\perp|^2\mathcal D$, which does not raise the bound.
\eeP
\begin{proof}
As in the proof of Proposition~\ref{p:spread}, entrywise $|e^{tT(k)}|\le e^{t(\lambda^2\mathcal D-c\lambda I+M-k^2\mathcal D)}\le e^{t(\lambda^2\mathcal D-c\lambda I+M)}$, and $T(0)$ has spectral abscissa $g(\lambda)-c\lambda$. The conclusion follows from $\min_\lambda c(\lambda)=c_*$.
\end{proof}

\beR[The Evans function]
[Evans function][Evans]
Linearize a front $U(z)$, $z=e\cdot\zeta-ct$, joining equilibria at $z=\pm\infty$: $Lp=\mathcal Dp''+cp'+f'(U)p$. The equation $Lp=\mu p$ is a first-order system $Y'=A(z,\mu)Y$, with $A\to A_\pm(\mu)$, where components with $d_i=0$ contribute one equation instead of two, provided $c\ne0$; at $c=0$ they give algebraic equations. To the right of the essential spectrum, $A_\pm(\mu)$ are hyperbolic, and the solutions $Y^-_j$ decaying at $-\infty$ and $Y^+_j$ decaying at $+\infty$, normalized by the asymptotic eigenvectors, depend analytically on $\mu$, and their numbers add up to the dimension. The Evans function is
\[
 D(\mu)=\det\bigl(Y^-_1,\ldots,Y^-_k,Y^+_{k+1},\ldots,Y^+_N\bigr)(0,\mu).
\]
It vanishes exactly at the eigenvalues, with order equal to the algebraic multiplicity \cite{AlexanderGardnerJones1990,Sandstede2002}; $(2\pi i)^{-1}\oint D'/D$ counts the eigenvalues inside a contour, and $D(0)=0$ because $LU'=0$, when $U'$ belongs to the chosen space and $0$ is an isolated eigenvalue. At a critical pulled front $0$ lies at the edge of the essential spectrum, and this argument at $0$ is not available in this form. For a scalar equation, $D$ is the Wronskian of $p_\pm$, with $p_\pm\sim e^{r_\pm z}$, $r_+=-\frac c2-\sqrt{\frac{c^2}4+\mu-f'(0)}$ and $r_-=-\frac c2+\sqrt{\frac{c^2}4+\mu-f'(s_a)}$.
\eeR

\beR[The fronts of the figures]
[stability of our fronts][FroSta]
\begin{enumerate}
\item \emph{The bistable quartic.} For $x_t=x_{\zeta\zeta}-x+\Phi_a(x)$ with the reference coefficients, $f<0$ on $(0,1/8)$, $f>0$ on $(1/8,1/2)$, with a triple zero at $1/8$. The front of Remark~\ref{r:OddFro} attracts, up to a shift, every front-like initial state $u_0$ with values in $[0,1/2]$, $\liminf_{z\to-\infty}u_0(z)>1/8$ and $\limsup_{z\to+\infty}u_0(z)<1/8$: after rescaling $[0,1/2]$ to $[0,1]$ and reflecting $z$, existence follows from \cite[Theorem 2.4(c)]{FifeMcLeod1977}, and convergence to a translate from \cite[Theorem 3.1]{FifeMcLeod1977}. These results do not require a nonzero derivative at the intermediate equilibrium. The two constant stable states show that some such restriction is needed. Spectrally, $v=e^{cz/2}p$ gives the self-adjoint problem $v''+(f'(U)-c^2/4)v=\mu v$, and since $U'$ has one sign, $\mu=0$ is its simple top eigenvalue. The essential spectrum of $L$ reaches $f'(0)=a_1-1=-1/343$, so the stability is very slow. Numerically, shooting gives $c\approx0.0535$, the next eigenvalue of the self-adjoint problem is about $-0.0034$, and $D$ has winding number one on the circles $|\mu|=0.0015$, $0.0025$ and on the half-disc $\Re\mu\ge-0.0015$, $|\mu|\le2$. The relaxation rate of pulled fronts is described in \cite{EbertVanSaarloos2000}. The value approximately $-0.0034$ concerns the self-adjoint transformed operator; it is not identified here as an isolated eigenvalue of the original operator on unweighted $L^2$. The fronts of $x_1,x_3$ in Figure~\ref{fig:overFront}(c) move at the speed $c$.
\item \emph{Concave threshold fronts (Figure~\ref{fig:spaFront}).} The system is cooperative, and the computed profiles are monotone, the setting of the principal eigenvalue theory for monotone fronts \cite{Volpert-alt2}. The fronts started from compact data are measured below $c_*$ (Remark~\ref{r:SpdAtt}), which we conjecture they approach, and at $c_*$ Proposition~\ref{p:WeiSpe} leaves no spectral gap.
\item \emph{Cholera (Figure~\ref{fig:choFront}).} The strain-$1$ front is a front of $(S,I_1,B_1)$; it is not monotone, $S$ overshoots to about $0.24$ below $S_1^*=0.3$, and the system is not cooperative. The resident behind it has $\sup_ks(J-k^2\mathcal D)\approx-0.141$, and the state ahead is covered by Proposition~\ref{p:WeiSpe}. Transversally, $R_{2|1}>1$, so the resident is unstable to strain $2$; a perturbation in $(I_2,B_2)$ grows, but since $c_2<c_1$ and $c_{2|1}<c_1$ the linear speeds leave it behind the front of strain $1$, as in panel (d).
\item \emph{The ignition region (Figure~\ref{fig:overFront}).} Its speed is still changing at $t=4000$ (Remark~\ref{r:OvlIgn}), so we make no stability statement.
\end{enumerate}
For prescribed $c>c_*$ we computed the profiles and the rightmost eigenvalues of the operator, on $z\in[-100,40]$ with step $0.05$, central differences for diffusing components and upwind differences for the others, no flux at $z=-100$, the end value at $z=40$ taken from the exponential traveling mode, and perturbations vanishing beyond $z=40$; the script \texttt{frontSta.py} gives the remaining details. On a bounded interval the weight does not change the exact eigenvalues, so these do not distinguish the weighted operator from the unweighted one. The discretized operator is strongly nonnormal, however, so we computed them after the diagonal similarity $e^{\lambda_*\max\{z,0\}}$; without it, round-off produces complex values that move with the interval. All eigenvalues found near the origin are real and negative, close to the edge $\lambda_*(c_*-c)$ of Proposition~\ref{p:WeiSpe}. Halving the step changes the entries of the table by up to $0.005$, and extending the interval to $z=60$ by up to $0.004$, while the next eigenvalues move toward the first, so we read these values as the edge of the discretized essential spectrum instead of isolated eigenvalues of the front. This is numerical evidence: stabilizing the state ahead does not by itself exclude unstable eigenvalues generated by the interior of the front, and the fronts at $c_*$ itself are not covered.
\begin{center}
\begin{tabular}{llll}
\hline
$c/c_*$ & $1.02$ & $1.1$ & $1.2$\\\hline
cholera, rightmost eigenvalue & $-0.012$ & $-0.070$ & $-0.145$\\
cholera, $\lambda_*(c_*-c)$ & $-0.014$ & $-0.070$ & $-0.140$\\
Figure~\ref{fig:spaFront}, rightmost eigenvalue & $-0.040$ & $-0.168$ & $-0.339$\\
Figure~\ref{fig:spaFront}, $\lambda_*(c_*-c)$ & $-0.031$ & $-0.155$ & $-0.310$\\\hline
\end{tabular}
\end{center}
\eeR
\section{Additional results and examples}

{\color{newcolor}
\paragraph{The response function survives in the fluctuation limit.}
For example, if a population of size $n$ has active fraction
$z$ and transition rates
$n(1-z)\Phi(h+\beta z)$ and $nz[1-\Phi(h+\beta z)]$,
its large-population limit satisfies
\[
\dot x=\Phi(h+\beta x)-x.
\]
Under the usual smoothness assumptions and convergence of the
initial fluctuations, the central limit approximation
$\sqrt n(z-x)\Rightarrow U$ satisfies
\[
dU_t=
[\beta \Phi'(h+\beta x_t)-1]U_t\,dt
+\sqrt{(1-x_t)\Phi(h+\beta x_t)
+x_t[1-\Phi(h+\beta x_t)]}\,dB_t.
\]
Thus the original response $\Phi$ determines the limiting dynamics
and the noise amplitude, while its derivative determines the
linear fluctuation drift. Explicit limit results retaining nonlinear
responses include the martingale central limit theorem for
stochastic Wilson--Cowan neural fields of
\cite[Theorem~2.3]{riedler2013}, and the functional central
limit theorem for age-dependent Hawkes networks of
\cite{chevallier2017}. An interesting direction for further work is
to establish response-dependent fluctuation limits under independent
inverse-subordinator time changes, connecting the present models with
the work of Leonenko, Macci and Pacchiarotti
\cite{LeonenkoMacciPacchiarotti2021}.}

{\color{newcolor}
\subsection{Equilibria are partitioned by their exact supra-siphons}
Counting equilibria, and counting the stable ones, needs two ingredients from
above: the exact zero set (Definition~\ref{d:ExZ} and
Proposition~\ref{p:zeroset}), which says \emph{where} an equilibrium is counted,
and the block-triangular Jacobian (Proposition~\ref{p:Tri}), which says
\emph{whether} it is stable. Throughout, $\dot x=f(x;\theta)$ depends on
parameters $\theta$ ranging over a region on which the family of siphons does
not change.

\beD[Exact-support counts and the full-stability weight][full-stability weight][FSW]
For a siphon $\Sigma$, or for $\Sigma=\varnothing$, write
\[
\mathcal F_\Sigma^\circ=\{x\in \mathcal F_\Sigma:\ x_j>0\ \text{for } j\notin \Sigma\},
\]
the interior of the siphon face $\mathcal F_\Sigma$; $\mathcal F_\varnothing^\circ$ is the interior
of the orthant, and $\varnothing$ serves as an index and never as a siphon.
The points of $\mathcal F_\Sigma^\circ$ are exactly the points whose exact zero set is $\Sigma$.
Let
\begin{enumerate}
\item $N_\Sigma(\theta)$ be the number of equilibria in $\mathcal F_\Sigma^\circ$, and
$N_\Sigma^{\mathrm{st}}(\theta)$ the number of those that are locally exponentially
stable;
\item $N_{\ge0}(\theta)$ and $N_{\ge0}^{\mathrm{st}}(\theta)$ be the number of
all nonnegative equilibria and of the locally exponentially stable ones.
\end{enumerate}
At a point $(0,y)\in \mathcal F_\Sigma^\circ$, with the blocks of Proposition~\ref{p:Tri}
for the siphon $\Sigma$, the \emph{full-stability weight} is
\begin{equation}\label{eq:weight}
W_\Sigma(y;\theta)=
\ind{J_{\mathrm{tan},\Sigma}(y;\theta)\text{ Hurwitz}}\,
\ind{\JtrS{\Sigma}(y;\theta)\text{ Hurwitz}},
\end{equation}
where an empty block counts as Hurwitz. The first factor tests stability
within the siphon face, the second tests that the absent variables cannot
invade.
\eeD

\beP[The exact supra-siphons partition the equilibria, and the full-stability weight selects the stable ones][counting by exact supra-siphon][Part]
Assume that every boundary equilibrium has a siphon as its exact zero set; by
Proposition~\ref{p:zeroset} this holds for reaction networks under its rate
assumptions. Then:
\begin{enumerate}
\item every nonnegative equilibrium lies in exactly one set $\mathcal F_\Sigma^\circ$, namely
the one indexed by its exact supra-siphon, or by $\varnothing$ if it is
positive; consequently
\[
N_{\ge0}=\sum_\Sigma N_\Sigma,\qquad N_{\ge0}^{\mathrm{st}}=\sum_\Sigma N_\Sigma^{\mathrm{st}},
\]
the sums running over $\varnothing$ and all siphons, not only the minimal
ones;
\item an equilibrium $(0,y)\in \mathcal F_\Sigma^\circ$ is locally exponentially stable if
and only if $W_\Sigma(y;\theta)=1$, so that
$N_\Sigma^{\mathrm{st}}=\sum_{(0,y)\in \mathcal F_\Sigma^\circ,\ f(0,y)=0}W_\Sigma(y;\theta)$.
\end{enumerate}
\eeP
\begin{proof}
For item~1, an equilibrium $e$ lies in $\mathcal F_\Sigma^\circ$ exactly when $\Sigma=Z(e)$, and
the sets $\mathcal F_\Sigma^\circ$ for distinct $\Sigma$ are therefore disjoint; by assumption
$Z(e)$ is a siphon whenever $e$ is on the boundary. Item~2 is item~2 of
Proposition~\ref{p:Tri} applied to the siphon $\Sigma=Z(e)$.
\end{proof}

\beR[Why minimal siphons cannot be used as counting classes][][noDouble]
Summing over the minimal siphons instead of the exact supra-siphons counts some
equilibria several times when minimal siphons overlap. For example, with the
two overlapping minimal siphons $\{1,2\}$ and $\{2,3\}$, the origin lies on the
faces of both, but its exact supra-siphon is $\{1,2,3\}$, so it is counted once.
Likewise, the first factor of $W_\Sigma$ alone counts the attractors of the dynamics
restricted to the siphon face, some of which are invaded by the absent
variables. Without the assumption of Proposition~\ref{p:Part}, a boundary
equilibrium may lie in no counting class at all, as in the example of a
positive ODE whose boundary limit point has a zero set that is not a siphon.
\eeR
\bookref{t:Ric}{These counts are used in Theorem~\ref{t:Ric}, whose stability
weight is \eqref{eq:weight}, and in the sum~\eqref{eq:sum} over exact
supra-siphons.}{}}

{\color{newcolor}
For $\dot x=f(x)$ of class $C^1$, take $g=f$ and
$A=J^{[k]}$, where $J=Df$ and $J^{[k]}$ is the additive compound on
$\bigwedge^k\mathbb R^n$, of dimension $\binom nk$.
If $\Psi$ is the state variational matrix, differentiation of its minors gives
\[
 \frac{d}{dt}\Psi^{(k)}=J^{[k]}\Psi^{(k)},\qquad \Psi^{(k)}(0)=I.
\]
Proposition~\ref{p:VarCert} therefore bounds infinitesimal
$k$-dimensional volumes. Linearity of $J\mapsto J^{[k]}$ carries a
convex enclosure by $J_1,\ldots,J_q$ to the enclosure by
$J_1^{[k]},\ldots,J_q^{[k]}$.

\emph{Orlando's identity ($n\ge2$).}
For $p(\lambda)=\det(\lambda I-J)=\lambda^n+a_1\lambda^{n-1}+\cdots+a_n$,
let $\Delta_{n-1}(p)=\det(a_{2j-i})_{i,j=1}^{n-1}$, with $a_0=1$
and $a_\ell=0$ outside $0\le\ell\le n$. Since the eigenvalues of
$J^{[2]}$ are $\lambda_i+\lambda_j$, $i<j$, Orlando's identity reads
\begin{equation}\label{HG:orlandoeq}
 \Delta_{n-1}(p)=\det(-J^{[2]})
 =(-1)^{\binom n2}\prod_{i<j}(\lambda_i+\lambda_j).
\end{equation}
At a Hopf point a pair sums to zero. Thus a contraction certificate
for $J^{[2]}$ excludes a Hopf point at any equilibrium where it holds.}

{\color{newcolor}
\paragraph{Biologically labelled splittings and the exit test.}
If a biological convention has already specified a nonnegative matrix
$F$ of new infections, put $V=F-\JtrS{\Sigma}$. This is a meaningful
splitting only after its matrix conditions have been checked. The
labels selecting $F$ are extra model data; stoichiometry alone need
not distinguish a newly infected individual from an internal transfer.

\beT[When the labelled loss matrix is admissible]\label{p:toolExit}
For a transverse Jacobian on a siphon face and a chosen $F\geq0$, the matrix $V=F-\JtrS{\Sigma}$ is a nonsingular
$M$-matrix if and only if
\[
 V_{ij}\leq0\quad(i\ne j),\qquad
 \exists q\gg0:\ q^\top V\gg0.
\]
In the special transfer-and-exit representation
\[
 Q=-V,\quad Q_{ij}=t_{ij}\geq0\ (i\ne j),\quad
 Q_{jj}=-d_j-\sum_{i\ne j}t_{ij},\qquad d_j\geq0,
\]
this holds exactly when every compartment can reach one with a
strictly positive exit rate $d_j$ along transfer edges $j\to i$.
Equivalently, every terminal strongly connected component contains
a positive exit.
\eeT
\begin{proof}
The first equivalence is the preceding Metzler criterion applied to
$-V$. For the second, $Q$ is the generator, in column convention, of
a finite transfer process killed at rates $d_j$. If every state can
reach an exit, choose a common finite time and a positive lower bound
on the exit probability from every initial state. Iterating gives
exponential decay of the surviving mass, hence $Q$ is Hurwitz.
Conversely a terminal component without exit conserves its mass;
its generator has eigenvalue zero, so $Q$ cannot be Hurwitz.
\end{proof}
This is a standard compartmental, or $M$-matrix, test
\cite{Berman,Van}.
It turns the question of an epidemiologically meaningful splitting
into a labelled modelling choice followed by a precise admissibility
test.}

\subsection{When output-complete buffering simplifies threshold networks}
\label{ss:OMthreshold}

{\color{newcolor}

We record a rank test for the Okada--Mochizuki buffering index. It
concerns a specified kinetic/stoichiometric representation: equivalent
representations of one ODE can have different output completions.

\beP[Unit-column realizations identify buffered sets by outgoing influences]
[][OMunit]
Let
\[
 \dot x_i=-d_i x_i+P_i(x),\qquad d_i>0,\quad P_i\ge0,
 \qquad i=1,\ldots,n,
\]
with differentiable rates. Choose a realization consisting of the loss
columns $-e_i$, with rates $d_ix_i$, and production columns $e_i$,
with rates $P_i(x)$. Define
\[
 D_i=\{j:\partial_jP_i\not\equiv0\},\quad
 T(M)=\{i:D_i\cap M\ne\varnothing\},\quad
 \lambda_{\rm OM}(M)=\operatorname{rank}\Gamma_{\operatorname{Out}(M)}-|M|.
\]
Then, for every $M\subseteq\{1,\ldots,n\}$,
\[
 \operatorname{rank}\Gamma_{\operatorname{Out}(M)}
       =|M\cup T(M)|,\qquad
 \lambda_{\rm OM}(M)=|T(M)\setminus M|.
\]
Consequently, $M$ is buffered precisely when no kinetic dependence
edge $j\to i$, $j\in M$, terminates outside $M$. If the dependence
graph is strongly connected, only the empty set and the full species
set are buffered.
\eeP
\begin{proof}
Every $j\in M$ places the loss column $-e_j$ in
$\operatorname{Out}(M)$. Every $i\in T(M)$ places the production
column $e_i$ there. These are all the selected nonzero columns,
so their span is the coordinate subspace indexed by $M\cup T(M)$.
\end{proof}

For $\dot x=-Vx+\Pi(x)$ with a general nonsingular $M$-matrix $V$,
the same test applies to the unit-column realization which represents
$-v_{ii}x_i$ by $X_i\to0$ and each positive cross-production
$-v_{ij}x_j$ ($i\ne j$) by $X_j\to X_j+X_i$. In that realization
one adds the transfer edges $j\to i$ with $-v_{ij}>0$ to the
kinetic influence graph, and the same outgoing-edge count gives
$\lambda_{\rm OM}(M)$. This representation is a choice, not an
ODE-invariant characterization of buffering.

The production columns in this proposition may be read as the formal
reactions $0\to X_i$ with state-dependent rates $P_i(x)$. This formal
realization is suitable for the OM rank calculation, but its formal
reactant sets need not reproduce the siphons of the ODE. When a
reactant-positive catalytic decomposition exists, it must be used
for simultaneous siphon and buffering calculations.

\beR[The coupled odds examples have no proper buffering sets]
In the two-population model~\eqref{eq:oddsNet}, for
$\gamma,\varepsilon>0$, the two production rates depend on both
$x$ and $y$, so the influence graph is strongly connected. Its
unit-column OM lattice is $\{\varnothing,\{x,y\}\}$. At zero coupling,
both singleton sets become buffered, although both singleton siphons
remain present throughout the family.

For the three-population odds system~\eqref{eq:KBoverlapCusp} with
$\alpha,\beta,\gamma>0$, the dependence edges include
$2\to1$, $1\to2$, $3\to2$ and $2\to3$. Thus its unit-column OM
lattice is $\{\varnothing,\{1,2,3\}\}$. In particular,
the overlapping minimal siphons $\{1,2\}$ and $\{2,3\}$ are not
buffered sets. These conclusions also hold for the following
reactant-positive catalytic realization when $\Phi_a$ is strictly
increasing: use the three loss reactions $X_i\to0$ at rates $x_i$,
the self-production reactions $X_i\to2X_i$ at rates $\Phi_a(x_i)$,
and the additional reactions
\[
\begin{array}{ll}
 X_2\to X_2+X_1,&
 \Phi_a(x_1+\alpha x_2)-\Phi_a(x_1),\\
 X_1+X_3\to X_1+X_3+X_2,&
 \Phi_a(x_2+\beta x_1x_3)-\Phi_a(x_2),\\
 X_2\to X_2+X_3,&
 \Phi_a(x_3+\gamma x_2)-\Phi_a(x_3).
\end{array}
\]
All listed catalytic rates are nonnegative and are strictly positive
when their reactants are positive. Their stoichiometric columns are
$e_1,e_2,e_3$, so the rank formula applies without changing the
siphon lattice.
\eeR

\beR[The displayed two-strain environmental-transmission network]
For the ten reactions of Figure~\ref{fig:choFront}, take the
species set $\{S,I_1,B_1,I_2,B_2\}$. Each species has an independent
loss column $-e_i$ with rate depending on that species. The kinetic
influence graph includes
\[
 S\to I_k,\qquad I_k\to B_k,\qquad
 B_k\to S,\ I_k\qquad(k=1,2).
\]
It is strongly connected. Every proper nonempty $M$ has a reaction
whose rate depends on a species in $M$ and whose stoichiometric column
has a nonzero coordinate outside $M$. Together with the independent
loss columns, this shows
$\operatorname{rank}\Gamma_{\operatorname{Out}(M)}>|M|$.
Thus its OM lattice also consists only of the empty and full sets.
The two-strain invasion blocks are therefore obtained from the
siphon-face transverse Jacobian, not from a nontrivial OM buffer.
\eeR

For a reducible dependence graph the rank test is useful: in a
feed-forward network with $X_1\to X_2$ as a kinetic influence and no
influence from $X_2$ to $X_1$, the set $\{X_2\}$ is buffered. The
steady-state equation for $X_1$ is independent of rates confined to
this downstream set. More generally, buffered sets are unions of
sink-closed blocks of the dependence graph. This can organize the
solution of boundary resident equations before evaluating invasion
matrices, but the actual reaction representation must be fixed.}

\fi
\appendix

\section{Positive systems, siphons, triangularity of the Jacobian on siphon faces and the Metzler property of its transversal block, and invasion exponents}\label{ss:siphons}

\beD[Positive ODE][][pos]
An ODE $\dot x=f(x)$, with $f$ locally Lipschitz on a neighbourhood of
$\R^{n_S}_{\ge0}$, is called \emph{positive} \cite{rantzer2015scalable} or
\emph{non-negative} \cite{haddad2010} if the orthant $\R^{n_S}_{\ge0}$ is forward
invariant. Equivalently, $f_i(x)\ge0$ whenever $x\ge0$ and $x_i=0$.
\eeD

\beD[Siphon and minimal siphon][siphon][Siphon]
For a set $\Sigma$ of coordinates, let
$\mathcal F_\Sigma=\{x\in\R^{n_S}_{\ge0}:\ x_i=0\ \text{for all } i\in\Sigma\}$
be the corresponding face of the orthant. A nonempty set $\Sigma$ is a
\emph{siphon} of a positive ODE if $\mathcal F_\Sigma$ is forward invariant, and
a \emph{minimal siphon} if no nonempty proper subset of $\Sigma$ is a siphon.
\eeD

Since the coordinates in $\Sigma$ vanish on $\mathcal F_\Sigma$, the face is
forward invariant if and only if $f_i$ vanishes on $\mathcal F_\Sigma$ for
every $i\in\Sigma$. A siphon is thus a set of variables that cannot be created
from nothing: once absent, they stay absent. For a reaction network with
mass-action kinetics and positive rate constants, or more generally with nonnegative rates that vanish when a reactant is absent
and are positive when all reactants are present, this is the
combinatorial condition that every reaction producing a species of $\Sigma$
consumes one; for more general rates, additional invariant faces can arise from
rates that vanish for other reasons.

\beD[Exact zero set, or exact supra-siphon, of a point][exact supra-siphon][ExZ]
The \emph{exact zero set} of a point $w\in\R^{n_S}_{\ge0}$ is $Z(w)=\{i:\ w_i=0\}$. When $Z(w)$ is
a siphon it is also called the \emph{exact supra-siphon} of $w$. Since $w\in\mathcal F_\Sigma$
exactly when $\Sigma\subseteq Z(w)$, the exact supra-siphon is the largest siphon on whose face
$w$ lies, it contains every siphon whose face contains $w$, and $w$ lies in the relative interior of that
siphon face,
\[
\mathcal F^\circ_{Z(w)}=\{x\in\mathcal F_{Z(w)}:\ x_j>0\ \text{for } j\notin Z(w)\}.
\]
For $w\gg0$, $Z(w)=\varnothing$ and $w$ lies in the interior of the orthant.
\eeD

\beP[The exact zero set of a RN boundary omega-limit point is a siphon][exact zero sets are siphons][zeroset]\label{p:AdLS}
Consider the reaction system
\[
\dot x=f(x)=\sum_r v_r(x)(b_r-a_r),
\]
where $a_r,b_r\in\R_{\ge0}^{n_S}$ are the reactant and product vectors. Assume that the rates are $C^1$
on a neighbourhood of the nonnegative orthant and, for $x\ge0$, satisfy
\[
v_r(x)\ge0,\qquad
v_r(x)>0
\quad\Longleftrightarrow\quad
x_j>0\ \text{for every }j\text{ with }a_{rj}>0
\]
(these assumptions hold for mass-action kinetics with positive rate constants). Let $x(t)$ be a
forward-complete solution with $x(0)\gg0$, and let $w$ be a boundary $\omega$-limit point. Then
its exact zero set $Z(w)$ is a siphon: every reaction producing a species in $Z(w)$ consumes at
least one species in $Z(w)$. The same holds for every boundary equilibrium $w$.
\eeP
\begin{proof}
First we show that $f_i(w)=0$ for every $i\in Z(w)$; for an equilibrium this holds by assumption.
Choose $t_k\to\infty$ with $x(t_k)\to w$. Local existence, uniqueness and continuous dependence
give a solution $y:(-\delta,\delta)\to\R^{n_S}$ through $y(0)=w$ such that, for each
$s\in(-\delta,\delta)$,
\[
y(s)=\lim_{k\to\infty}x(t_k+s)\ge0.
\]
Thus, for $i\in Z(w)$, the differentiable function $y_i$ has a local minimum at $s=0$, and
$f_i(w)=y_i'(0)=0$.
Fix $i\in Z(w)$. Every reaction consuming species $i$ has zero rate at $w$, because $w_i=0$. Hence
\[
0=f_i(w)=\sum_{r:\,a_{ri}=0}b_{ri}v_r(w).
\]
All summands are nonnegative. Therefore every reaction producing $i$ without consuming it has zero
rate at $w$. By the rate assumption, such a reaction consumes some species whose concentration at $w$
is zero, hence a species in $Z(w)$. A reaction that consumes $i$ already consumes a species in
$Z(w)$. This proves the siphon property.
\end{proof}

\beXa[For a general positive ODE, an exact zero-set of  a boundary omega-limit point  may  fail to be siphon][][ZsNsip] On $[0,1]\times[0,2]$ the system
$\dot x=(y-1)^2-x$, $\dot y=1-y$ sends every trajectory to $(0,1)$, whose exact zero set $\{x\}$ is not a siphon,
since $\dot x=(y-1)^2>0$ on $\{x=0\}$ for $y\ne1$. The exact zero set of a point is distinguished from the larger invariant
faces containing it.
\eeXa

For a   siphon $\Sigma$, $z\in\R^m_{\ge0}$ collects the coordinates in
$\Sigma$ and $y\in\R^n_{\ge0}$ the remaining, \emph{resident}, ones; we write
$\dot z=f_z(z,y)$, $\dot y=f_y(z,y)$, with $f$ of class $C^1$, and at a point
$(0,y)$ of the face we call
\[
\Jtr(y):=D_zf_z(0,y),\qquad J_{\mathrm{tan}}(y):=D_yf_y(0,y),\qquad
J_{\mathrm{mix}}(y):=D_zf_y(0,y)
\]
the \emph{transversal}, \emph{tangential} and \emph{mixed} blocks of the
Jacobian.
 The transversal block is also written $J^\perp$ or
$\JtrS{\Sigma}$; these are the same matrix. Invariance of the face reads
\begin{equation}\label{eq:faceinv}
f_z(0,y)=0\qquad\text{for all } y\in\R^n_{\ge0}.
\end{equation} $J_{\mathrm{tan}}$ is the Jacobian of the ODE restricted to the face,
while $\Jtr$ governs the linearized dynamics of the absent variables.

Below, $s(A)=\max\{\operatorname{Re}\mu:\ \mu\in\operatorname{spec}A\}$ denotes the
spectral abscissa of a square matrix $A$.

\beP[The Jacobian on a siphon face is block triangular, with a Metzler transversal block][][Tri]
\label{p:Sta}
Let $\Sigma$ be a siphon of a positive $C^1$ system. Then:
\BEN
\im At every point of its face, $D_yf_z(0,y)=0$ and
\begin{equation}\label{eq:blocks}
J(0,y)=
\begin{pmatrix}
\Jtr(y)&0\\
J_{\mathrm{mix}}(y)&J_{\mathrm{tan}}(y)
\end{pmatrix}.
\end{equation}
Its spectrum is the union of the spectra of the diagonal blocks.
\im An equilibrium $E=(0,y^*)$ is locally exponentially stable if and only if
$J_{\mathrm{tan}}(y^*)$ and $\Jtr(y^*)$ are Hurwitz. It is unstable if either
block has an eigenvalue with positive real part.
\im The transversal block $\Jtr(y)$ is Metzler at every point of the face;
the tangential block need not be Metzler.
\EEN
The same conclusions hold on a face of an invariant box $[0,\bar x]$, $\bar x\gg0$,
with derivatives understood in admissible directions when the vector field
is only $C^1$ up to the boundary. Empty blocks impose no condition.
\eeP
\begin{proof}
Differentiating~\eqref{eq:faceinv} in the resident variables gives
$D_yf_z(0,y)=0$ in the relative interior and, by continuity, on the entire
face. The characteristic polynomial therefore factors over the two diagonal
blocks, giving the stability and instability assertions. Finally, for distinct $i,j\in\Sigma$, positivity at $(te_j,y)$ gives
\[
(f_z)_i(te_j,y)\ge0=(f_z)_i(0,y)\qquad(t>0).
\]
Dividing by $t$ and letting $t\downarrow0$ proves
$(\Jtr)_{ij}(y)\ge0$. These arguments use only admissible directions and
apply equally on a box.
\end{proof}

\beD[Transversal block and invasion exponent][invasion exponent][Inv]
At a fixed point $E=(0,y^*)$ of the face of a siphon $\Sigma$, the
\emph{transversal block} is $M_\Sigma(E):=\Jtr(y^*)$, and the \emph{invasion
exponent} of $\Sigma$ at $E$ is
\[
\lambda_\Sigma(E):=s\bigl(M_\Sigma(E)\bigr).
\]
\eeD

By Proposition~\ref{p:Tri}, $M_\Sigma(E)$ is Metzler and
$\lambda_\Sigma(E)$ is an eigenvalue of the full Jacobian. For sufficiently
large $c$, Perron--Frobenius theory gives
$\lambda_\Sigma(E)=\rho(M_\Sigma(E)+cI)-c$, with a nonnegative eigenvector.
It measures asymptotic growth in the transverse linearized system:
$\lambda_\Sigma(E)>0$ permits invasion, whereas
$\lambda_\Sigma(E)<0$ gives exponential decay of all transverse perturbations;
at equality, linearization does not decide invasion.

\beP[transverse factorization on siphon faces][][traFac]\label{t:Inv}
For a reaction network  as defined in Proposition \ref{p:zeroset}, write $\nu^-_\ell\to\nu^+_\ell$, and   assume
each rate vanishes exactly when at least one of its reactant
concentrations vanishes (this  holds for mass-action
kinetics).

For any siphon $\Sigma$
\[
 f_\Sigma(z,y)=M_\Sigma(z,y)z,\qquad
 M_\Sigma(z,y)=\int_0^1D_zf_\Sigma(tz,y)\,dt,
 \qquad M_\Sigma(0,y)=\JtrS{\Sigma}(y).
\]
The matrix $M_\Sigma$ need not be Metzler away from the face. If the ODE
preserves the nonnegative orthant, the block-triangular Jacobian on
$\mathcal F_\Sigma$, its Metzler transverse block, and the stability test at a face
equilibrium are those of Proposition~\ref{p:Tri}.
\eeP
\begin{proof}
On the face all consumption rates of absent species vanish. If the
reaction condition holds, all production rates vanish too. Conversely,
at a relative-interior point of the face, a reaction producing an
absent species without an absent reactant has positive rate. It cannot
be cancelled by consumption of that species, so the face is not
invariant. The integral identity is the fundamental theorem of
calculus applied along $t\mapsto(tz,y)$; at $z=0$ it gives
$M_\Sigma(0,y)=D_zf_\Sigma(0,y)=\JtrS{\Sigma}(y)$.
\end{proof}

\paragraph{Siphons in epidemiology.}
In the multi-strain models considered here, \emph{epidemiologic strains} often
correspond to \emph{minimal siphons}. By a \emph{mutation} we mean only a coupling
between strains that changes the Perron--Frobenius form of the transversal block
from block diagonal to block triangular. The totality of ``infected variables'' is
the total siphon, the union of all minimal siphons, and a \DFE\ (DFE) lies on its
face; a larger siphon may be needed when the face of the total siphon does not
contain a unique fixed point. For a splitting $\Jtr=F-V$ with
$F\ge0$ and $V$ a nonsingular $M$-matrix \cite{Diek,Van},
the Perron spectral radius $\rho(FV^{-1})$ of the next generation matrix on a
siphon, expressed as a function of the resident variables and called a
reproduction function in \cite{AAN,AAHJ,AABH25}, characterizes invasion of boundary
fixed points, and often their existence domains; other eigenvalues of the next
generation matrix play no role in the linear invasion threshold.

\beR[historical notes][][hisS]
The invariant face-siphons correspondence and its use in persistence are developed in
\cite{AdLS,ShiuStu}.

The Jacobian and transversal factorizations are  elementary \cite{AAHK}, \cite{AAH26} and also supply
the idea of the role of Perron data in face-repulsion.

The exact-zero-set assertion of Proposition~\ref{p:zeroset} is already contained
in \cite[Lemma~5.3 and Proposition~5.4]{AdLS}, for real-analytic rates; the proof above gives a
short argument using a local trajectory through the omega-limit point, and needs only $C^1$ rates.
The invariant-face equivalence is \cite[Proposition~5.5]{AdLS}.
\eeR
\ifjournalversion\else
\section{A basic toolbox for siphon calculus}
\label{ss:siphonToolbox}

Throughout this toolbox, $\dot x=f(x)$ is a positive $C^1$ system on a
neighbourhood of the orthant; local Lipschitz regularity suffices for
invariance and mass balance. We use the face, block and spectral notation
of Section~\ref{ss:siphons}. Reaction-network statements retain the
exact-zero rate condition of Proposition~\ref{p:zeroset}.

\paragraph{Conservation classes and inaccessible extinction faces.}
Let $S$ be the stoichiometric matrix. A nonnegative conservation vector
means $q\geq0$, $q\ne0$, $q^\top S=0$; its associated conserved quantity
is $q^\top x$. Conservation laws with mixed signs are not substitutes
for this condition.

\beP[Conservation excludes a siphon face]\label{p:toolConservation}
Suppose $q\geq0$, $q\ne0$, $q^\top S=0$, and
$\operatorname{supp}q\subseteq\Sigma$. On a trajectory with
$q^\top x(0)=c>0$,
\[
 q_\Sigma^\top x_\Sigma(t)=c,
 \qquad \|x_\Sigma(t)\|_1\geq c/\|q\|_\infty.
\]
Consequently that trajectory cannot approach $\mathcal F_\Sigma$.
For a reaction network satisfying the exact-zero condition, if every
nonempty siphon contains such a conservation support, then every bounded
trajectory starting strictly positive satisfies
$\liminf_{t\to\infty}x_i(t)>0$ for every $i$.
\eeP
\begin{proof}
The first assertion follows from $q^\top\dot x=0$ and
$q_\Sigma^\top x_\Sigma\leq\|q\|_\infty\|x_\Sigma\|_1$.
For the last assertion, a bounded trajectory has a compact omega-limit
set. By Proposition~\ref{p:zeroset}, the nonempty exact zero set of a
boundary limit point would be a siphon. Its conservation support has
positive conserved mass along the trajectory but zero mass at that point,
a contradiction. Compactness gives the stated coordinate lower bounds.
\end{proof}
This is the familiar noncritical-siphon obstruction to extinction
\cite{AdLS,ShiuStu}. It yields trajectory-dependent lower bounds;
it does not assert one uniform bound across conservation classes whose
conserved masses may approach zero. A conserved quantity can also
produce zero eigenvalues in the ambient Jacobian: stability must then
be tested within the stoichiometric compatibility class.

\paragraph{Positive weights: one test for stability, decay and comparison.}

\beP[Collatz--Wielandt bounds for Metzler matrices][][ColWie]
For a Metzler matrix $A$, the following are equivalent:
\[
 s(A)<0;\qquad \exists v\gg0:\ Av\ll0;
 \qquad \exists q\gg0:\ q^\top A\ll0.
\]
For every $v\gg0$,
\[
 \min_i\frac{(Av)_i}{v_i}\leq s(A)
 \leq\max_i\frac{(Av)_i}{v_i}.
\]
If a nonnegative solution satisfies
$\dot z\leq Az$ componentwise, and $q^\top A\leq-\alpha q^\top$
for some $\alpha>0$, then
\[
 q^\top z(t)\leq e^{-\alpha t}q^\top z(0).
\]
\eeP
\begin{proof}
If $A$ is Hurwitz, its positive semigroup gives
$(-A)^{-1}=\int_0^\infty e^{At}\,dt\geq0$.
Choose $v=(-A)^{-1}\mathbf1$ and
$q^\top=\mathbf1^\top(-A)^{-1}$; their components are strictly
positive and the strict inequalities follow. Conversely, a positive
left or right weight with strict decrease makes the associated
weighted column or row sums strictly negative, and bounds $s(A)<0$.
The two-sided bounds follow by diagonally scaling with $v$ and
comparing the row sums of a shifted nonnegative matrix.
The final assertion follows from
$\frac{d}{dt}(q^\top z)\leq-\alpha q^\top z$ and Gronwall's inequality.
\end{proof}
These are classical nonsingular $M$-matrix and positive-systems
criteria \cite{Berman}, \cite{rantzer2015scalable}. A uniform comparison on a forward-invariant
region upgrades a local invasion test to an extinction estimate on
that region. Without the comparison hypothesis, a stable boundary
equilibrium alone says nothing about large introductions.

\paragraph{Comparison matrices and explicit extinction rates.}
The preceding positive-weight criterion becomes an explicit decay test
when the comparison matrix is written as $C-V$. The next theorem applies
to any nonnegative differential subsolution; a siphon supplies one
important application.

\beT[A Metzler comparison inequality gives explicit exponential decay]
[comparison decay][Cmp]
\emph{Scope: any nonnegative absolutely continuous function satisfying the
differential inequality; no condition is placed on the other variables.}

Let $V\in\mathbb R^{n\times n}$ be a nonsingular $M$-matrix and $C\ge0$ with $\rho(CV^{-1})<1$.
Then $V-C$ is a nonsingular $M$-matrix. Put
\[
 p=(V-C)^{-1}\mathbf1\gg0,\qquad c=\frac1{\max_ip_i},\qquad
 w^\top=\mathbf1^\top(V-C)^{-1}\gg0,\qquad c'=\frac1{\max_iw_i}.
\]
If $z:[0,\infty)\to\mathbb R_+^n$ is locally absolutely continuous and
$\dot z\le(C-V)z$ almost everywhere, componentwise, then
\[
 \max_i\frac{z_i(t)}{p_i}\le e^{-ct}\max_i\frac{z_i(0)}{p_i},\qquad
 L=w^\top z\ \text{ satisfies }\ \dot L\le-c'L\quad\text{a.e.}
\]
In particular $z(t)\to0$ exponentially.
\eeT
\begin{proof}
$(V-C)^{-1}=V^{-1}\sum_{k\ge0}(CV^{-1})^k\ge V^{-1}\ge0$ and $V-C$ is a
$Z$-matrix, so it is a nonsingular $M$-matrix; $V^{-1}$ has no zero row or
column, hence $p,w\gg0$. The matrix $C-V$ is Metzler and $(C-V)p=-\mathbf1\le-cp$.
The function $r(t)=\max_i z_i(t)/p_i$ is locally absolutely continuous.
At almost every time, its derivative is the maximum of
$\dot z_i(t)/p_i$ over the indices attaining the maximum. For each such
index, $z_i=rp_i$, $z_j\le rp_j$, and the Metzler property gives
\[
 \dot z_i\le((C-V)z)_i\le r((C-V)p)_i=-r\le-crp_i.
\]
Thus $\dot r\le-cr$ almost everywhere, and integration proves the
weighted maximum estimate. For the functional, $w\ge0$ gives
$\dot L\le w^\top(C-V)z=-\mathbf1^\top z\le-L/\max_iw_i$.
\end{proof}

For a siphon, apply the theorem to $z=x_\Sigma$ after verifying the
comparison inequality throughout the stated invariant region. No
cooperativity or convergence assumption on the resident variables is
needed. The inequality may also depend on time, provided the same
constant matrices $C,V$ bound it almost everywhere. A stable transverse
Jacobian on the face alone does not provide this global comparison.

\paragraph{What the siphon lattice does, and does not, tell the linear graph.}
At a relative-interior point of $\mathcal F_\Sigma$, let $G_\Sigma$ have the
directed edge $j\to i$ for $i\ne j$ when
$(\JtrS{\Sigma})_{ij}>0$.

\beT[Siphons constrain the linear influx graph][siphons and the influx graph][Graph]\label{p:toolGraph}\label{t:Sip}
For mass-action kinetics with positive rate constants, if
$T\subseteq\Sigma$ is a siphon, then $T$ is predecessor-closed in
$G_\Sigma$: $i\in T$ and $j\to i$ imply $j\in T$.
If every reaction producing a species of $\Sigma$ has exactly one
reactant molecule in $\Sigma$, the converse also holds.
Thus, in this first-order setting, minimal nonempty siphons contained
in $\Sigma$ are the source strongly connected components of
$G_\Sigma$. They need not be all its strongly connected components.
\eeT
\begin{proof}
An off-diagonal edge arises from a reaction of order one in $\Sigma$:
its unique reactant in $\Sigma$ is $j$, and it produces $i\ne j$.
If $i\in T$, the siphon property forces that reactant $j$ into $T$.
For the converse, take a reaction producing $i\in T$. Its unique
$\Sigma$ reactant is $j$. If $j=i$, it consumes $T$; otherwise it
gives the edge $j\to i$ and predecessor closure forces $j\in T$.
The last assertion is the source-component characterization of
minimal nonempty predecessor-closed sets.
\end{proof}
This is the elementary graph reduction behind monomolecular siphon
analysis \cite{AdLS,ShiuStu,Berman}. \label{pg:2A2B}The added first-order hypothesis is
essential. The networks
\[
 A\to B,\quad B\to A,
 \qquad\text{and}\qquad
 2A\to B,\quad 2B\to A
\]
have the same sole nonempty siphon $\{A,B\}$, but at its zero face the
first transverse graph is strongly connected and the second has no
edges. Adding $A\to\varnothing$, $B\to\varnothing$ to both preserves
the example and gives negative diagonal entries in the second.
Hence the siphon lattice alone cannot determine the linear influx
graph without information about reaction order. Moreover, when the
transverse graph is reducible,
$s(\JtrS{\Sigma})=\max_B s(J_B)$ in Frobenius form, and a
supercritical block may fail to invade from introductions confined
to components which cannot reach it.

\paragraph{From local invasion to persistence.}
The appropriate resident object is an equilibrium, periodic orbit, or
invariant measure, not necessarily a fixed point. Conservation can
exclude some boundary faces outright; positive invasion weights can
repel others. Uniform persistence additionally needs a compact
absorbing setting and control of all relevant boundary recurrent
sets, as provided by an acyclic Morse decomposition or a suitable
average-Lyapunov criterion \cite{FRT94,HS10,HofSch,SmithThieme2011}. A list of
positive invasion numbers at selected equilibria is not a substitute
for checking those hypotheses.

\subsection{Positive variational equations and converse contraction weights}
\label{s:PositiveWeights}
Let $\varphi_t$ be the flow of $\dot y=g(y)$ on a compact forward-invariant
set $K$. For a continuous matrix field $A$, write $\Psi(t,y)$ for the
fundamental matrix of
\[
 \dot\xi=A(\varphi_t(y))\xi,\qquad \Psi(0,y)=I.
\]
For a $C^1$ field $H$, put $D_hH=(\partial_yH)g$, componentwise.
The matrix $A$ must be specified: estimates for the full, transverse,
and exterior-power variational equations have different consequences.

\beP[Common and transported certificates for linear equations along a flow]
[][VarCert]
Suppose $g$ is $C^1$ near $K$ and $A:K\to\mathbb R^{d\times d}$ is
continuous. Each of the following conditions implies uniform exponential decay.
\begin{enumerate}
\item A symmetric $C^1$ field $P$ satisfies, for $0<m\le M$ and $c>0$,
\begin{equation}\label{HG:transportLMI}
 mI\preceq P\preceq MI,\qquad
 D_hP+A^\top P+PA\preceq-2cP\qquad(y\in K).
\end{equation}
Then $\|\Psi(t,y)\|_2\le\sqrt{M/m}\,e^{-ct}$.
\item A fixed diagonal signature matrix $\Sigma$, with entries $\pm1$,
makes $C=\Sigma A\Sigma^{-1}$ Metzler on $K$, and a $C^1$ field $v$ satisfies
\begin{equation}\label{HG:transportLP}
 m\mathbf1\le v\le M\mathbf1,\qquad
 D_hv+C^\top v\le-cv\qquad(y\in K).
\end{equation}
Then $\|\Psi(t,y)\|_1\le(M/m)e^{-ct}$.
This assertion also holds for a continuous flow and continuous $v$
with a continuous forward orbital derivative, used in place of $D_hv$.
\end{enumerate}
The norms $(\xi^\top P(y)\xi)^{1/2}$ and $v(y)^\top|\Sigma\xi|$
respectively decay at rate $c$ along trajectories.

For constant certificates the derivative terms vanish. If
$A(K)\subseteq\operatorname{conv}\{A_1,\ldots,A_q\}$, it suffices to check
one of the following systems with the same certificate at every vertex:
\begin{equation}\label{HG:LMI}
 P\succ0,\qquad A_j^\top P+PA_j\preceq-2cP\quad(1\le j\le q),
\end{equation}
\begin{equation}\label{HG:LP}
 v\gg0,\qquad C_j^\top v\le-cv\quad(1\le j\le q),
 \qquad C_j=\Sigma A_j\Sigma^{-1}\text{ Metzler}.
\end{equation}
For fixed $c$, these are LMIs in $P$ and linear inequalities in $v$,
respectively. They express logarithmic-norm bounds in the constant norms
with transformations $\Theta=P^{1/2}$ and
$\Theta=\operatorname{Diag}(v)\Sigma$.
\eeP
\begin{proof}
Along a solution,
\[
 \frac{d}{dt}(\xi^\top P\xi)
 =\xi^\top(D_hP+A^\top P+PA)\xi
 \le-2c\,\xi^\top P\xi,
\]
which proves the first estimate. For the second, positivity of the
Metzler equation for $\Sigma\xi$ gives, almost everywhere,
\[
 \frac{d}{dt}\bigl(v^\top|\Sigma\xi|\bigr)
 \le(D_hv+C^\top v)^\top|\Sigma\xi|
 \le-cv^\top|\Sigma\xi|.
\]
Integration and the bounds on $v$ prove the claim, including its
orbital-derivative version. Convex combinations preserve the constant
inequalities and the Metzler property. The logarithmic-norm statements
follow from the symmetric part in the Euclidean case and the weighted
column sums in the $\ell^1$ case.
\end{proof}

\beD[Lur'e systems are absolutely contracting when one norm contracts them for a whole class of nonlinearities]
[absolute Lur'e contractivity][LurCon]
A \emph{Lur'e system} is $\dot x=-Vx+B\Phi(Wx,t)$, with $V\in\mathbb R^{n_S\times n_S}$,
$B\in\mathbb R_+^{n_S\times m}$, $W\in\mathbb R_+^{m\times n_S}$. Let $P\succ0$, $\|x\|_P^2=x^\top Px$ and
$\eta>0$. The right-hand side $f$ is \emph{strongly infinitesimally contracting} with respect to
$\|\cdot\|_P$ with rate $\eta$ if
\[
 (f(x_1,t)-f(x_2,t))^\top P(x_1-x_2)\le-\eta\|x_1-x_2\|_P^2
 \qquad\text{for all }x_1,x_2,t.
\]
The system is \emph{absolutely contracting} over a class of nonlinearities if this holds with
the same $P$ and $\eta$ for every $\Phi$ in the class. The class of
\cite{davydov2024exponential} consists of the \emph{cocoercive} maps,
\[
 (\Phi(y_1,t)-\Phi(y_2,t))^\top(y_1-y_2)\ge\rho\|\Phi(y_1,t)-\Phi(y_2,t)\|_2^2,
 \qquad\rho>0.
\]
\eeD

\beR[Absolute Lur'e contractivity and the test $\bar R_0<1/L$]
By \cite[Theorem~1]{davydov2024exponential}, absolute contractivity over the cocoercive class
holds if and only if some $\lambda\ge0$ satisfies
\[
 \begin{pmatrix}
 -V^\top P-PV+2\eta P & PB+\lambda W^\top\\
 B^\top P+\lambda W & -2\lambda\rho I
 \end{pmatrix}\preceq0.
\]
A coordinatewise response with $0\le \Phi'\le L$ is cocoercive with $\rho=1/L$, so
$\dot x=-Vx+\Phi(\beta Wx)$ is the case $B=I$, with $\beta W$ in place of $W$. Put
$\bar R_0=\beta\rho(V^{-1}W)$. For $V=I$, $\beta>0$ and symmetric $W\ge0$ with $r=\rho(W)>0$, the
matrix inequality with $P=I$, $\lambda=1/(\beta r)$ reduces to $\beta Lr\le1-\eta$, so it agrees
with the test $\bar R_0<1/L$.

In general the cocoercive class is strictly larger and the two conditions differ, as the
following illustration, which is not taken from \cite{davydov2024exponential}, shows. Take $V=I$,
$W=\begin{pmatrix}0&a\\0&0\end{pmatrix}$ with $\beta aL>2$. Then $\bar R_0=0$, and every coordinatewise
response gives a triangular system, contracting in the norm $\max\{|z_1|/p_1,|z_2|\}$ for any
$p_1>\beta aL$. The linear map
$\Phi(y)=\frac L2\begin{pmatrix}1&-1\\1&1\end{pmatrix}y$ is cocoercive with $\rho=1/L$, but
$-I+\beta\frac L2\begin{pmatrix}1&-1\\1&1\end{pmatrix}W$ has the eigenvalue $-1+\beta aL/2>0$, so no
$P$ gives absolute contractivity.

With $V$ a nonsingular $M$-matrix and $W\ge0$, absolute contractivity over the cocoercive class
implies $\bar R_0<1/L$: the admissible map $\Phi(y)=Ly$ requires $-V+\beta LW$ to be Hurwitz. The
example above shows that the converse fails. Thus the distinction goes beyond the choice of
norm: the cocoercive class also allows coupling between the coordinates of the response.
For firing-rate networks with monotone nonexpansive coordinatewise activations,
\cite{Gokhale2026} give a sufficient condition with a positive diagonal multiplier. An elementary
deduction from it recovers the spectral bound exactly for $V=I$, including nonsymmetric and
reducible $W\ge0$. Put $M=\beta LW$ and normalize the response to have slopes in $[0,1]$. The
criterion is
\[
 \begin{pmatrix}
 -2(1-\eta)P & P+M^\top Q\\
 P+QM & -2Q
 \end{pmatrix}\preceq0,
 \qquad P\succ0,\quad Q\succ0\ \text{diagonal},\quad\eta>0.
\]
Testing the identity activation shows that its feasibility implies $\rho(M)<1$. Conversely, if
$\rho(M)<1$, the nonnegative matrix $G=(I+M)/2$ is Schur stable and admits a positive diagonal
matrix $P$ with $G^\top PG\prec P$ \cite{rantzer2015scalable}. Taking $Q=P$, the Schur complement is
\[
 2\bigl(G^\top PG-(1-\eta)P\bigr),
\]
which is negative definite for sufficiently small $\eta>0$. The triangular example above satisfies
this criterion: it separates the coupled cocoercive class from the coordinatewise class, not
the test $\bar R_0<1/L$ from the coordinatewise criterion.
These LMIs use a quadratic function of the difference of two states. The classical
Lur'e--Postnikov quadratic-plus-integral construction, associated with Popov criteria
\cite{Khalil2002}, instead works about an equilibrium. For $\dot x=-Vx+\Phi(Wx+h)$ with $V$ positive
diagonal, slopes in $[0,1]$ and an equilibrium $x^*$, put $e=x-x^*$ and
$\psi_i(s)=\Phi_i((Wx^*+h)_i+s)-\Phi_i((Wx^*+h)_i)$, so that $\dot e=-Ve+\psi(We)$. With $X\succ0$
and $N=\operatorname{Diag}(\nu_i)\succeq0$, $T\succeq0$, both diagonal, the function
\[
 L_{\mathrm{LP}}(e)=e^\top Xe+2\sum_i\nu_i\int_0^{(We)_i}\psi_i(s)\,ds
\]
is bounded above and below by quadratics. Writing $z=\psi(We)$, differentiation gives
\[
 \dot L_{\mathrm{LP}}=2e^\top X(-Ve+z)+2z^\top NW(-Ve+z).
\]
Since $T$ is diagonal, the sector inequalities $\psi_i(s)(s-\psi_i(s))\ge0$ give
$z^\top T(We-z)\ge0$. Adding twice this quantity bounds $\dot L_{\mathrm{LP}}$ above by the quadratic
form in $(e,z)$ of the matrix in the \emph{Lur'e--Postnikov LMI}
\begin{equation}\label{eq:LPLMI}
 \begin{pmatrix}
 -XV-VX & X-VW^\top N+W^\top T\\
 X-NWV+TW & NW+W^\top N-2T
 \end{pmatrix}\prec0.
\end{equation}
Feasibility of \eqref{eq:LPLMI} therefore implies global exponential stability of $\bar x$. It
certifies stability of that equilibrium instead of contraction, and this function is distinct
from the Fenchel--Bregman divergence.
The subcritical contraction statement of Theorem~\ref{t:KBsharp} belongs to
the established theory of coordinatewise slope-restricted Lur'e systems \cite{Davydov}.
Its constructive conclusion is different: with $V,W,L$ fixed, every prescribed finite number of
stable equilibria can be realized by one common scalar response whenever $\bar R_0>1/L$. No conclusion
at $\bar R_0=1/L$ is asserted.
\eeR

For monotone systems and positive linear time-varying systems the converse principle is
\cite[Theorem~4.3, Corollary~4.8]{KawanoBesselinkCao2020}. The next theorem formulates it over a
compact base flow: the base need not be cooperative, $A$ need not be its Jacobian, and only
forward evolution is used.

\beT[Cooperative variational equations admit state-dependent contraction weights exactly when they are uniformly exponentially stable]
[converse contraction weights][Wgt]\label{t:SCM-converse}
Let $K$ be compact and forward invariant under $\dot y=g(y)$. Let $g$ be continuous and uniquely integrable near $K$, with continuous
flow $\varphi_t$, and let $A:K\to\mathbb R^{n\times n}$ be continuous and
Metzler. Let $\Psi(t,y)$ be the fundamental matrix of
$\dot\xi=A(\varphi_t(y))\xi$, with $\Psi(0,y)=I$.
The following are equivalent:
\begin{enumerate}
\item There are $M,\alpha>0$ such that
$\|\Psi(t,y)\|\le Me^{-\alpha t}$ for all $t\ge0$ and $y\in K$.
\item There are $c>0$ and a continuous, uniformly positive weight
$v:K\to\mathbb R^n$ with continuous forward orbital derivative
\[
 \dot v(y)=\left.\frac{d}{dt}v(\varphi_t(y))\right|_{t=0+},
 \qquad \dot v+A^\top v\le-cv\quad\text{on }K.
\]
\end{enumerate}
Under these conditions, every fixed $q\gg0$ gives the weight
\[
 v(y)=\int_0^\infty\Psi(t,y)^\top q\,dt,
 \qquad \dot v+A^\top v=-q.
\]
\eeT
\begin{proof}
The implication (2)$\Rightarrow$(1) follows from
Proposition~\ref{p:VarCert} with $\Sigma=I$ and the minimum and maximum
of the components of $v$ on $K$.
For (1)$\Rightarrow$(2), the displayed integral converges uniformly and
is continuous. Since $A$ is bounded, $\Psi(t,y)=I+O(t)$ uniformly near
$t=0$; positivity gives $v_i(y)\ge\delta q_i/2$ for some $\delta>0$.
The cocycle identity yields, for $s\ge0$,
\[
 v(\varphi_s(y))=
 \Psi(s,y)^{-\top}\left[v(y)-\int_0^s\Psi(t,y)^\top q\,dt\right].
\]
Differentiation at $s=0$ gives the asserted continuous orbital derivative.
Boundedness of $v$ then gives $q\ge cv$ for some $c>0$.
\end{proof}

\beR[Smooth variational equations admit smooth contraction weights]
[][WgtSmooth]
Under the equivalent conditions of Theorem~\ref{t:Wgt}, if $g$ and $A$
extend as $C^1$ fields near $K$, the weight can be chosen $C^1$ near $K$.
Indeed, fix $q\gg0$ and choose $T>0$ so that
$\Psi(T,y)^\top q\le q/2$ on $K$. The finite-horizon weight satisfies
\[
 v_T(y)=\int_0^T\Psi(t,y)^\top q\,dt,\qquad
 D_hv_T+A^\top v_T=\Psi(T,y)^\top q-q\le-q/2.
\]
This is the construction of \cite[Appendix~B, (61)]{KawanoBesselinkCao2020}, there with
$q=\mathbf1$. Finite-time smooth dependence gives $C^1$ regularity; uniform positivity
and boundedness give $D_hv_T+A^\top v_T\le-cv_T$ for some $c>0$.
If $K$ is a rectangular box, $C^1$ tensor Bernstein approximation
provides a polynomial weight: sufficiently close approximants preserve
positivity and the strictly negative residual. Polynomial $g,A$ are
not required for existence, but permit polynomial certificate searches.
\eeR

The next theorem is a compact-domain consequence of the converse principle. Its hypothesis is
global convergence to a Hurwitz equilibrium on a forward invariant set, whereas
\cite[Theorem~4.3]{KawanoBesselinkCao2020} starts from incremental exponential stability on a
forward and backward invariant set.

\beT[On a compact invariant convex set, global convergence to a Hurwitz equilibrium of a cooperative system is equivalent to the existence of contraction weights]
[][EqWgt]
Let $K\subset\mathbb R^n$ be compact, convex, with nonempty interior,
and forward invariant under a $C^2$ vector field $f$ defined near $K$.
Suppose $J(x)=Df(x)$ is Metzler on $K$, $f(e)=0$ for some $e\in K$,
and $J(e)$ is Hurwitz. The following are equivalent:
\begin{enumerate}
\item Every solution starting in $K$ converges to $e$.
\item The state variational matrix $\Psi(t,x)=D_x\varphi_t(x)$ satisfies
$\|\Psi(t,x)\|\le Me^{-\alpha t}$ for some $M,\alpha>0$, all $t\ge0$,
and all $x\in K$.
\item There are $c>0$ and a $C^1$ weight $v\gg0$ near $K$ such that
\[
 (\partial_xv)f+J^\top v\le-cv\qquad\text{on }K.
\]
\end{enumerate}
In this case $e$ is globally exponentially stable relative to $K$.
If $K$ is a rectangular box, $v$ may be chosen polynomial.
\eeT
\begin{proof}
Theorem~\ref{t:Wgt} and Remark~\ref{r:WgtSmooth} give (2)$\Leftrightarrow$(3)
and the polynomial conclusion. Convexity gives
\[
 \varphi_t(x)-e
 =\int_0^1\Psi(t,e+\theta(x-e))(x-e)\,d\theta,
\]
so (2) implies exponential convergence and hence (1).

Assume (1). Since $J(e)$ is Metzler and Hurwitz,
\[
 w=(-J(e))^{-\top}\mathbf1\gg0,\qquad P_w=\operatorname{Diag}(w).
\]
For $\mu_1(H)=\max_j(H_{jj}+\sum_{i\ne j}|H_{ij}|)$, continuity gives
$\lambda>0$ and a neighbourhood $U$ of $e$ such that
$\mu_1(P_wJ(x)P_w^{-1})\le-\lambda$ on $U\cap K$.
Local asymptotic stability gives a relatively open forward-invariant
neighbourhood $V$ of $e$ in $K$, contained in $U$.
Every $x\in K$ enters $V$ at some time $t_x$; continuity gives a
relatively open neighbourhood of $x$ entering $V$ at the same time.
A finite subcover and forward invariance give $T<\infty$ such that
$\varphi_t(K)\subset V$ for all $t\ge T$.
With $L_w=\sup_{x\in K}\|P_wJ(x)P_w^{-1}\|_1$, the initial bound
$e^{L_wt}$ and subsequent decay at rate $\lambda$ yield
\[
 \|P_w\Psi(t,x)P_w^{-1}\|_1
 \le e^{(L_w+\lambda)T}e^{-\lambda t}
 \qquad(t\ge0,\ x\in K).
\]
Equivalence of norms proves (2).
\end{proof}

In the setting of Theorem~\ref{t:EqWgt}, the weight $v$ of item (3) defines the intrinsic distance
\[
 d_v(x,y)=\inf_\gamma\int_0^1v(\gamma(\theta))^\top|\gamma'(\theta)|\,d\theta,
\]
over $C^1$ paths in $K$ from $x$ to $y$. Item (3) and the Metzler property give
$d_v(\varphi_t(x),\varphi_t(y))\le e^{-ct}d_v(x,y)$, so $L(x)=d_v(x,e)$ is an equilibrium Lyapunov
function, possibly nonsmooth, with $L(\varphi_t(x))\le e^{-ct}L(x)$.

\emph{Transverse application.}
For a resident flow, taking $A=J^\perp$ in
Proposition~\ref{p:VarCert} or Theorem~\ref{t:Wgt} certifies the transverse
variational equation without requiring a cooperative resident field.
Theorem~\ref{t:EqWgt}, in contrast, concerns the full Jacobian; its
convergence conclusion cannot be inferred from transverse decay alone.

\beP[Composite Lyapunov functions combine resident storage with transverse decay][][ComLya]
Write the equations as $\dot z=f_z(z,y)$, $\dot y=f_y(z,y)$ with $z\ge0$ and $f_z(0,y)=0$, and let
$y^*$ be a resident equilibrium. On a compact forward-invariant set $K\subseteq\{z\ge0\}$ containing
$(0,y^*)$, suppose there are $\ell\gg0$, $a,c>0$, $b\ge0$, and a continuous resident storage
$U\ge0$, vanishing only at $y^*$ and differentiable near the resident projection of $K$, such that
\[
 \ell^\top f_z(z,y)\le-c\ell^\top z,
\]
\[
 \nabla U(y)^\top f_y(0,y)\le-a\|y-y^*\|^2,
 \qquad
 \nabla U(y)^\top\bigl(f_y(z,y)-f_y(0,y)\bigr)\le b\ell^\top z
 \quad ((z,y)\in K).
\]
Then, for any $\kappa>b/c$, the function $H(z,y)=U(y)+\kappa\ell^\top z$ satisfies
\[
 \dot H\le-a\|y-y^*\|^2-(\kappa c-b)\ell^\top z,
\]
and $(0,y^*)$ is globally asymptotically stable relative to $K$.
\eeP
\begin{proof}
Adding the three inequalities gives the estimate for $\dot H$. Since $z\ge0$ and $\ell\gg0$, $H\ge0$ on
$K$ and vanishes only at $(0,y^*)$, and the right-hand side is negative elsewhere. Thus $H$ is a strict
Lyapunov function on the compact forward-invariant set $K$.
\end{proof}
This is a standard composite-Lyapunov argument. The first inequality concerns the full nonlinear
equations on $K$, beyond their transverse linearization on the siphon face.
The linear term is the zero-reference limit of Subsection~\ref{s:FenBre}.
Proposition~\ref{p:ComBas} gives an application to a threshold network.

\subsection{Volterra functions are Fenchel--Bregman divergences}\label{s:FenBre}
A \emph{Lyapunov--Volterra function} is a weighted sum of terms
$x_i-x_i^*-x_i^*\log(x_i/x_i^*)$ about an equilibrium $x^*$. These functions and the term
$\int\Phi^{-1}$ in the energy of symmetric firing-rate networks \cite{Betteti2025} come from one
construction, a convex generator $h$ with $h'=\Phi^{-1}$; the Lyapunov--Volterra case is $\Phi(v)=e^v$.
The relevant objects are the Fenchel coupling, also called the Fenchel--Young gap, and its
representation as a Bregman divergence; see \cite[Appendix~C]{MertikopoulosSandholm2016}. We call
this quantity the \emph{Fenchel--Bregman divergence}. It gives a nonnegative
candidate from a response function; decrease along trajectories requires a separate calculation.

\paragraph{An increasing response determines a convex storage.}
Let $\Phi:\mathbb R\to I$ be a strictly increasing $C^1$ bijection onto an open interval $I$, with
$\Phi'(v)>0$ for every $v$. Choose $x_0\in I$ and define
\[
 h(x)=\int_{x_0}^{x}\Phi^{-1}(s)\,ds,\qquad x\in I.
\]
Then
\[
 h'(x)=\Phi^{-1}(x),\qquad
 h''(x)=(\Phi^{-1})'(x)=\frac{1}{\Phi'(\Phi^{-1}(x))}>0,
\]
so $h$ is strictly convex. Its convex conjugate is $h^*(v)=\sup_{z\in I}\{vz-h(z)\}$. The maximum is
attained uniquely at $z=\Phi(v)$, so
\[
 h^*(v)=v\Phi(v)-h(\Phi(v)),\qquad (h^*)'(v)=\Phi(v).
\]
The Fenchel coupling is
\[
 \mathcal S_h(p,v)=h(p)+h^*(v)-pv,\qquad p\in I,\quad v\in\mathbb R.
\]
The Fenchel--Young inequality gives $\mathcal S_h(p,v)\ge0$, with equality exactly when $p=\Phi(v)$.
Adding a constant to $h$ does not change this coupling. Define the Bregman divergence by
\[
 D_h(p,x)=h(p)-h(x)-h'(x)(p-x).
\]
Substituting $v=\Phi^{-1}(x)=h'(x)$ gives the exact identity
\[
 \mathcal S_h(p,\Phi^{-1}(x))=D_h(p,x).
\]
Thus fixing $p=x^*$ produces an equilibrium-relative candidate that is nonnegative and vanishes only
at $x=x^*$. A response with plateaus has no inverse, and needs a separate convex-subdifferential
argument.

\beP[Symmetric threshold networks admit a response-integral energy][][SymEne]
Consider $\dot x=-x+\Phi(\beta Wx)$ with $W=W^\top$ and $\beta\ge0$, where the same $C^1$ response
$\Phi:\mathbb R\to(0,1)$ is applied coordinatewise, $\Phi$ is onto, and $\Phi'(q)>0$ for every $q$. Choose
$s_0\in(0,1)$ and put
\[
 h(s)=\int_{s_0}^{s}\Phi^{-1}(r)\,dr,\qquad
 \mathcal E(x)=\sum_{i=1}^{n_S}h(x_i)-\frac{\beta}{2}x^\top Wx,\qquad x\in(0,1)^{n_S}.
\]
Then
\[
 \dot{\mathcal E}(x)=-\sum_{i=1}^{n_S}\bigl(\Phi^{-1}(x_i)-\beta(Wx)_i\bigr)\bigl(x_i-\Phi(\beta(Wx)_i)\bigr)\le0,
\]
with equality exactly at equilibria. If $e\in(0,1)^{n_S}$ is an equilibrium, then
\[
 \mathcal E(x)-\mathcal E(e)=\sum_{i=1}^{n_S}D_h(x_i,e_i)-\frac{\beta}{2}(x-e)^\top W(x-e),
\]
\[
 \nabla^2\mathcal E(e)=\operatorname{Diag}\Bigl(\frac{1}{\Phi'(\beta(We)_i)}\Bigr)-\beta W.
\]
This Hessian is positive definite if and only if the equilibrium Jacobian is Hurwitz, and then
$\mathcal E(x)-\mathcal E(e)$ is a strict local Lyapunov function. On any compact forward-invariant
subset of $(0,1)^{n_S}$, every trajectory approaches the equilibrium set; if the equilibria in that
set are isolated, each trajectory converges to one equilibrium.
\eeP
\begin{proof}
Symmetry of $W$ gives $\partial_{x_i}\mathcal E=\Phi^{-1}(x_i)-\beta(Wx)_i$, and the displayed derivative
follows from the equations. Since $\Phi$ is strictly increasing, the two factors in each summand have the
same sign, and their product vanishes exactly when $x_i=\Phi(\beta(Wx)_i)$. At an equilibrium,
$h'(e_i)=\beta(We)_i$; expanding the quadratic term and using the definition of $D_h$ gives the
energy-difference identity. With $S=\operatorname{Diag}(\Phi'(\beta(We)_i))\succ0$, the equilibrium
Jacobian is $J(e)=-I+\beta SW=-S\nabla^2\mathcal E(e)$, which is similar to the symmetric matrix
$-S^{1/2}\nabla^2\mathcal E(e)S^{1/2}$. The convergence statements follow from LaSalle's invariance
principle and, for isolated equilibria, connectedness of the omega-limit set.
\end{proof}

This is the classical symmetric-network energy calculation \cite{Hopfield1984,Betteti2025}, written
with Bregman divergences. Here the divergence is $D_h(x,e)$, whereas the Lyapunov--Volterra function
below uses $D_h(e,x)$. The energy gives convergence to the equilibrium set; in a multistable network
global attraction to one equilibrium does not follow.

\paragraph{The exponential response gives exactly the Lyapunov--Volterra function.}
For $\Phi(v)=e^v$ and $I=(0,\infty)$, choose $h(x)=x\log x-x$, so that $h^*(v)=e^v$. For $x,x^*>0$,
\[
 \mathcal S_h(x^*,\log x)=h(x^*)+h^*(\log x)-x^*\log x=x-x^*-x^*\log\frac{x}{x^*}.
\]
Consequently, for $c_i>0$ and $x_i^*>0$,
\[
 L(x;x^*)=\sum_{i=1}^{n_S}c_i\,\mathcal S_h(x_i^*,\log x_i)
 =\sum_{i=1}^{n_S}c_i\left(x_i-x_i^*-x_i^*\log\frac{x_i}{x_i^*}\right)
\]
is exactly the weighted Lyapunov--Volterra function. The order of the arguments matters:
\[
 D_h(x^*,x)=x-x^*-x^*\log\frac{x}{x^*},\qquad
 D_h(x,x^*)=x\log\frac{x}{x^*}-x+x^*.
\]
The second expression is the entropy term familiar from complex-balanced reaction networks. The two
use the same convex generator with opposite argument order; their dissipation properties must be
checked for the respective dynamics. The exponential response is unbounded and not globally
Lipschitz, so results requiring a global slope bound do not apply to it.

\paragraph{The derivative identifies the compatible equations.}
For coordinatewise responses $\Phi_i$ satisfying the assumptions above, let $h_i'=\Phi_i^{-1}$, fix
$x^*\in\prod_iI_i$, and set
\[
 L_\Phi(x;x^*)=\sum_{i=1}^{n_S}c_i\,\mathcal S_{h_i}(x_i^*,\Phi_i^{-1}(x_i))
 =\sum_{i=1}^{n_S}c_iD_{h_i}(x_i^*,x_i).
\]
Differentiation gives, along any differentiable trajectory,
\[
 \dot L_\Phi=\sum_{i=1}^{n_S}c_i(x_i-x_i^*)(\Phi_i^{-1})'(x_i)\dot x_i.
\]
In particular, suppose the equations have the form
\[
 \dot x_i=\frac{g_i(x)}{(\Phi_i^{-1})'(x_i)},\qquad g(x^*)=0.
\]
Then $\dot L_\Phi=\sum_ic_i(x_i-x_i^*)g_i(x)$, so the response determines the scalar factors that
cancel in the derivative. For $\Phi_i(v)=e^v$ this form is $\dot x_i=x_ig_i(x)$, the multiplicative
structure of generalized Lotka--Volterra equations. Not every threshold or production system has
this form; it must be checked, or the displayed derivative calculated directly.

For $g(x)=A(x-x^*)$ and $C=\operatorname{Diag}(c_1,\ldots,c_{n_S})$,
\[
 \dot L_\Phi=\frac12(x-x^*)^\top(CA+A^\top C)(x-x^*).
\]
Therefore $CA+A^\top C\preceq0$ makes $L_\Phi$ nonincreasing, while $CA+A^\top C\prec0$ makes its
derivative strictly negative away from $x^*$. For exponential responses this is the classical
Volterra--Lyapunov diagonal-stability calculation. Global convergence additionally requires the
appropriate existence and compactness or coercivity hypotheses.

\paragraph{The zero-reference limit produces linear terms for absent species.}
For fixed $x>0$, $\lim_{p\downarrow0}(x-p-p\log(x/p))=x$. Equivalently, extending $h(x)=x\log x-x$
by $h(0)=0$ gives $\mathcal S_h(0,\log x)=x$. Thus a boundary equilibrium naturally leads to
\[
 L(x;x^*)=\sum_{\{i:x_i^*>0\}}c_i\left(x_i-x_i^*-x_i^*\log\frac{x_i}{x_i^*}\right)
 +\sum_{\{i:x_i^*=0\}}c_ix_i.
\]
The resident terms are Volterra terms and the absent-species terms are linear. This identifies a
common convex origin for the two parts of a Perron--Volterra
ansatz. The choice of their weights and the sign of their
combined derivative remain model-dependent, and dissipation of the absent terms still requires a
transverse estimate; Proposition~\ref{p:LVboundary} carries this out for the three-population network.
\section{Dependencies between the results}\label{s:dep}
Figure~\ref{fig:dep} records use of a proved conclusion or criterion in
another result, including use without an explicit internal citation.
An arrow $A\to B$ means that the proof of $B$ uses that conclusion of $A$.
Panels group related questions, not connected components of the proof graph;
placing two results together does not assert a proof dependency.

\begin{figure}[H]
\centering
\begin{minipage}{\linewidth}
\centering\textbf{Arbitrary multistability: gain bounds and fixed siphons}\par
\begin{tikzpicture}[x=1cm,y=-1.25cm,font=\footnotesize,
 T/.style={draw,rectangle,inner sep=4pt,fill=blue!6,align=center},
 P/.style={draw,rounded corners=3pt,inner sep=3pt,fill=green!8},
 C/.style={draw,rounded corners=3pt,inner sep=3pt},
 arr/.style={-{Stealth[length=4pt]},thin,gray!70!black}]
\node[C] (cEnt) at (-2.5,2) {\hyperref[c:Ent]{C\,\ref*{c:Ent}}};
\node[P] (pColWie) at (-2.5,0) {\hyperref[p:ColWie]{P\,\ref*{p:ColWie}}};
\node[T] (tKBsharp) at (-2.5,1) {\hyperref[t:KBsharp]{T\,\ref*{t:KBsharp}}\\Sharp uniform NGM bound};
\node[T] (tsiphMulti) at (2.5,1) {\hyperref[t:siphMulti]{T\,\ref*{t:siphMulti}}\\Fixed siphons and transverse matrices};
\draw[arr] (pColWie.south) -- (tKBsharp.north);
\draw[arr] (tKBsharp.south) -- (cEnt.north);
\end{tikzpicture}\par
{\footnotesize Independent constructions; neither theorem uses the other.}
\end{minipage}\par\bigskip

\begin{minipage}{\linewidth}
\centering\textbf{Diffusion and spreading}\par
\begin{tikzpicture}[x=1.2cm,y=-1.5cm,font=\footnotesize,
 T/.style={draw,rectangle,inner sep=3pt,fill=blue!6},
 P/.style={draw,rounded corners=3pt,inner sep=3pt,fill=green!8},
 C/.style={draw,rounded corners=3pt,inner sep=3pt},
 arr/.style={-{Stealth[length=4pt]},thin,gray!70!black}]
\node[C] (cDifSta) at (0,0) {\hyperref[c:DifSta]{C\,\ref*{c:DifSta}}};
\node[T] (tKBsharp) at (1.2,0) {\hyperref[t:KBsharp]{T\,\ref*{t:KBsharp}}};
\node[P] (tspread) at (3,0) {\hyperref[p:spread]{P\,\ref*{p:spread}}};
\node[P] (pOdd) at (4.8,0) {\hyperref[p:Odd]{P\,\ref*{p:Odd}}};
\node[T] (tSpe) at (6,0) {\hyperref[t:Spe]{T\,\ref*{t:Spe}}};
\node[T] (tCmp) at (7,0) {\hyperref[t:Cmp]{T\,\ref*{t:Cmp}}};
\node[T] (tEqWgt) at (8,0) {\hyperref[t:EqWgt]{T\,\ref*{t:EqWgt}}};
\node[P] (pComBas) at (8.9,0) {\hyperref[p:ComBas]{P\,\ref*{p:ComBas}}};
\node[P] (pLVboundary) at (9.9,0) {\hyperref[p:LVboundary]{P\,\ref*{p:LVboundary}}};
\node[T] (tSpdBd) at (0.6,1) {\hyperref[t:SpdBd]{T\,\ref*{t:SpdBd}}};
\node[T] (tSpdQC) at (2.2,1) {\hyperref[t:SpdQC]{T\,\ref*{t:SpdQC}}};
\node[P] (pWeiSpe) at (3.4,1) {\hyperref[p:WeiSpe]{P\,\ref*{p:WeiSpe}}};
\node[C] (cOddSpd) at (4.8,1) {\hyperref[c:OddSpd]{C\,\ref*{c:OddSpd}}};
\node[C] (cDifAll) at (7,1) {\hyperref[c:DifAll]{C\,\ref*{c:DifAll}}};
\node[P] (pDifLya) at (9.4,1) {\hyperref[p:DifLya]{P\,\ref*{p:DifLya}}};
\draw[arr] (cDifSta) -- (tSpdBd);
\draw[arr] (tKBsharp) -- (tSpdBd);
\draw[arr] (tspread) -- (tSpdBd);
\draw[arr] (tspread) -- (tSpdQC);
\draw[arr] (tspread) -- (pWeiSpe);
\draw[arr] (tspread) -- (cOddSpd);
\draw[arr] (pOdd) -- (cOddSpd);
\draw[arr] (tSpe) -- (cDifAll);
\draw[arr] (tCmp) -- (cDifAll);
\draw[arr] (tEqWgt) -- (cDifAll);
\draw[arr] (pComBas) -- (pDifLya);
\draw[arr] (pLVboundary) -- (pDifLya);
\draw[arr] (cDifSta) to[out=25,in=155,looseness=.45] (cDifAll);
\end{tikzpicture}\par
{\footnotesize C\,\ref{c:DifSta} has its own comparison proof and gives item 2 of T\,\ref{t:SpdBd}; P\,\ref{p:spread} gives item 3, and item 1 is proved directly. Its statement applies to the equilibria constructed in T\,\ref{t:siphMulti}, T\,\ref{t:KBsharp}, T\,\ref{t:Net} and T\,\ref{t:KBoverlapCusp}, which are not drawn as arrows. Items of C\,\ref{c:DifAll} use T\,\ref{t:Spe}, T\,\ref{t:Cmp} and T\,\ref{t:EqWgt} in turn; P\,\ref{p:ComBas} gives item 3 of P\,\ref{p:DifLya}, and P\,\ref{p:LVboundary} gives item 4.}
\end{minipage}\par\bigskip

\begin{minipage}[t]{0.48\linewidth}
\centering\textbf{Reduction tools and their applications}\par
\begin{tikzpicture}[x=1.65cm,y=-1.12cm,font=\footnotesize,
 T/.style={draw,rectangle,inner sep=3pt,fill=blue!6},
 P/.style={draw,rounded corners=3pt,inner sep=3pt,fill=green!8},
 C/.style={draw,rounded corners=3pt,inner sep=3pt},
 arr/.style={-{Stealth[length=4pt]},thin,gray!70!black}]
\node[P] (pKBpositive) at (1,0) {\hyperref[p:KBpositive]{P\,\ref*{p:KBpositive}}};
\node[T] (tSpe) at (1,1) {\hyperref[t:Spe]{T\,\ref*{t:Spe}}};
\node[P] (pColWie) at (0,1) {\hyperref[p:ColWie]{P\,\ref*{p:ColWie}}};
\node[P] (pKBshared) at (0,2) {\hyperref[p:KBshared]{P\,\ref*{p:KBshared}}};
\node[T] (tKBhom) at (1,2) {\hyperref[t:KBhom]{T\,\ref*{t:KBhom}}};
\node[T] (tKBexchange) at (2,2) {\hyperref[t:KBexchange]{T\,\ref*{t:KBexchange}}};
\draw[arr] (pKBpositive) -- (tSpe);
\draw[arr] (tSpe) -- (pKBshared);
\draw[arr] (pColWie) -- (pKBshared);
\draw[arr] (tSpe) -- (tKBhom);
\draw[arr] (tSpe) -- (tKBexchange);
\end{tikzpicture}
\end{minipage}\hfill
\begin{minipage}[t]{0.48\linewidth}
\centering\textbf{Cusp criteria and explicit constructions}\par
\begin{tikzpicture}[x=1.65cm,y=-1.12cm,font=\footnotesize,
 T/.style={draw,rectangle,inner sep=3pt,fill=blue!6},
 P/.style={draw,rounded corners=3pt,inner sep=3pt,fill=green!8},
 C/.style={draw,rounded corners=3pt,inner sep=3pt},
 arr/.style={-{Stealth[length=4pt]},thin,gray!70!black}]
\node[P] (psadCus) at (0,0) {\hyperref[p:sadCus]{P\,\ref*{p:sadCus}}};
\node[P] (pKBpositive) at (1,0) {\hyperref[p:KBpositive]{P\,\ref*{p:KBpositive}}};
\node[T] (pQua) at (2,1) {\hyperref[t:Qua]{T\,\ref*{t:Qua}}};
\node[T] (tActF) at (0,1) {\hyperref[t:ActF]{T\,\ref*{t:ActF}}};
\node[T] (tcuspAsym) at (0,2) {\hyperref[t:cuspAsym]{T\,\ref*{t:cuspAsym}}};
\node[T] (tKBhom) at (0,3) {\hyperref[t:KBhom]{T\,\ref*{t:KBhom}}};
\node[T] (tNet) at (1,3) {\hyperref[t:Net]{T\,\ref*{t:Net}}};
\node[T] (tKBoverlapCusp) at (2,3) {\hyperref[t:KBoverlapCusp]{T\,\ref*{t:KBoverlapCusp}}};
\node[C] (cNin) at (2,4) {\hyperref[c:Nin]{C\,\ref*{c:Nin}}};
\draw[arr] (psadCus) -- (tActF);
\draw[arr] (pKBpositive) -- (tActF);
\draw[arr] (psadCus) to[bend right=35] (tcuspAsym);
\draw[arr] (psadCus) to[bend right=35] (tKBhom);
\draw[arr] (psadCus) -- (pQua);
\draw[arr] (psadCus) -- (tNet);
\draw[arr] (pQua) -- (tNet);
\draw[arr] (pQua) to[bend right=35] (tKBoverlapCusp);
\draw[arr] (psadCus) -- (tKBoverlapCusp);
\draw[arr] (tKBoverlapCusp) -- (cNin);
\end{tikzpicture}
\end{minipage}\par\bigskip
\begin{minipage}[t]{0.48\linewidth}
\centering\textbf{Production representation, siphon structure and boundary stability}\par
\begin{tikzpicture}[x=1.65cm,y=-1.12cm,font=\footnotesize,
 T/.style={draw,rectangle,inner sep=3pt,fill=blue!6},
 P/.style={draw,rounded corners=3pt,inner sep=3pt,fill=green!8},
 C/.style={draw,rounded corners=3pt,inner sep=3pt},
 arr/.style={-{Stealth[length=4pt]},thin,gray!70!black}]
\node[P] (pzeroset) at (0,0) {\hyperref[p:zeroset]{P\,\ref*{p:zeroset}}};
\node[P] (pTri) at (1,0) {\hyperref[p:Tri]{P\,\ref*{p:Tri}}};
\node[P] (pKBpositive) at (2,0) {\hyperref[p:KBpositive]{P\,\ref*{p:KBpositive}}};
\node[P] (tCon) at (0,1) {\hyperref[p:toolConservation]{P\,\ref*{p:toolConservation}}};
\node[T] (pInvF) at (1,1) {\hyperref[t:InvF]{T\,\ref*{t:InvF}}};
\node[T] (tActF) at (2,1) {\hyperref[t:ActF]{T\,\ref*{t:ActF}}};
\node[P] (pKBshared) at (2,2) {\hyperref[p:KBshared]{P\,\ref*{p:KBshared}}};
\node[P] (ptraFac) at (0,3) {\hyperref[p:traFac]{P\,\ref*{p:traFac}}};
\node[T] (tGraph) at (1,3) {\hyperref[t:Graph]{T\,\ref*{t:Graph}}};
\node[P] (pOdd) at (2,3) {\hyperref[p:Odd]{P\,\ref*{p:Odd}}};
\node[T] (tNet) at (1,4) {\hyperref[t:Net]{T\,\ref*{t:Net}}};
\node[T] (tKBoverlapCusp) at (2,4) {\hyperref[t:KBoverlapCusp]{T\,\ref*{t:KBoverlapCusp}}};
\node[P] (pKBrepresentation) at (0,4) {\hyperref[p:KBrepresentation]{P\,\ref*{p:KBrepresentation}}};
\draw[arr] (pzeroset) -- (tCon);
\draw[arr] (pTri) -- (pInvF);
\draw[arr] (pTri) -- (ptraFac);
\draw[arr] (pTri) to[bend right=35] (tGraph);
\draw[arr] (pTri) to[bend right=35] (tNet);
\draw[arr] (pTri) -- (tKBoverlapCusp);
\draw[arr] (pKBpositive) -- (tActF);
\draw[arr] (pKBpositive) to[bend right=35] (pKBshared);
\draw[arr] (pInvF) -- (pOdd);
\end{tikzpicture}
\end{minipage}\hfill
\begin{minipage}[t]{0.48\linewidth}
\centering\textbf{Positive weights and contraction}\par
\begin{tikzpicture}[x=1.65cm,y=-1.12cm,font=\footnotesize,
 T/.style={draw,rectangle,inner sep=3pt,fill=blue!6},
 P/.style={draw,rounded corners=3pt,inner sep=3pt,fill=green!8},
 C/.style={draw,rounded corners=3pt,inner sep=3pt},
 arr/.style={-{Stealth[length=4pt]},thin,gray!70!black}]
\node[P] (pColWie) at (0,0) {\hyperref[p:ColWie]{P\,\ref*{p:ColWie}}};
\node[P] (pVarCert) at (2,0) {\hyperref[p:VarCert]{P\,\ref*{p:VarCert}}};
\node[T] (tCmp) at (0,1) {\hyperref[t:Cmp]{T\,\ref*{t:Cmp}}};
\node[P] (pSymEne) at (0,2) {\hyperref[p:SymEne]{P\,\ref*{p:SymEne}}};
\node[P] (pComLya) at (0,3) {\hyperref[p:ComLya]{P\,\ref*{p:ComLya}}};
\node[P] (pComBas) at (1,3) {\hyperref[p:ComBas]{P\,\ref*{p:ComBas}}};
\node[P] (pLVboundary) at (2,3) {\hyperref[p:LVboundary]{P\,\ref*{p:LVboundary}}};
\draw[arr] (pComLya) -- (pComBas);
\draw[arr] (pComBas) -- (pLVboundary);
\node[T] (tWgt) at (2,1) {\hyperref[t:Wgt]{T\,\ref*{t:Wgt}}};
\node[T] (tEqWgt) at (2,2) {\hyperref[t:EqWgt]{T\,\ref*{t:EqWgt}}};
\draw[arr] (pColWie) -- (tCmp);
\draw[arr] (pColWie) -- (tEqWgt);
\draw[arr] (pVarCert) -- (tWgt);
\draw[arr] (tWgt) -- (tEqWgt);
\end{tikzpicture}
\end{minipage}\par\bigskip
\caption{Use of proved results. Parentheses in the group headings refer to
the enumeration in the abstract. Repeated nodes
(P\,\ref{p:KBpositive}, P\,\ref{p:KBshared}, P\,\ref{p:ColWie},
T\,\ref{t:ActF}, T\,\ref{t:KBhom}, T\,\ref{t:Net},
T\,\ref{t:KBoverlapCusp}) identify links between panels.
The reduction tools are consolidated in T\,\ref{t:Spe};
the quartic construction and minimum-degree result form the single
node T\,\ref{t:Qua}.
Arrows indicate proof usage only; panels are thematic, not connected components.
Nodes without incident arrows have independent proofs and no recorded downstream use.}
\label{fig:dep}
\end{figure}
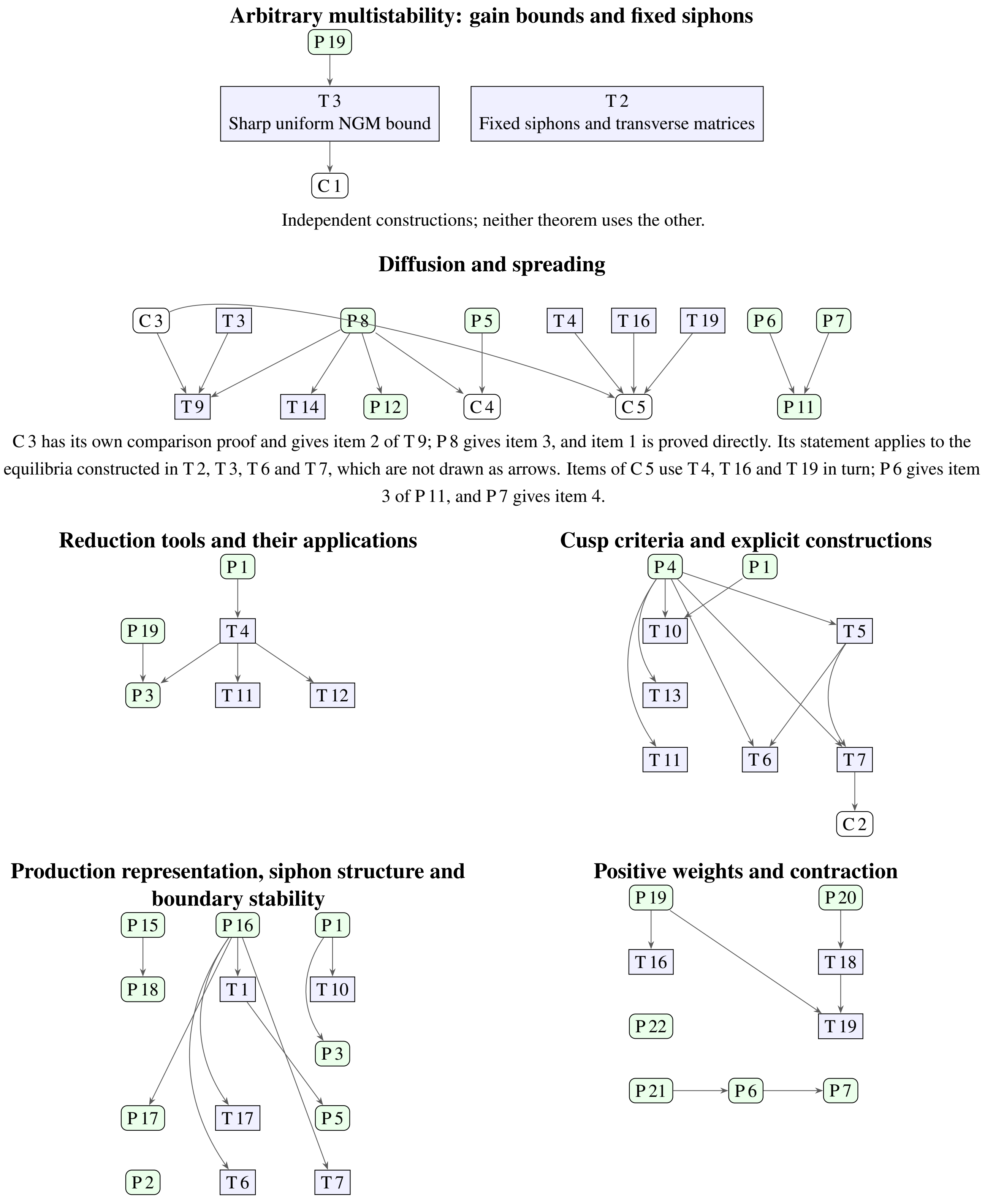

\fi
\funding{This research received no external funding.}

\dataavailability{The numerical parameters are given in the text. The accompanying scripts \texttt{verify\_cusp\_tables.py} and \texttt{verify\_125.py} reproduce the main cusp and equilibrium-count computations. Additional wave simulations, included only in the arXiv version, use separate front scripts.}

\useofartificialintelligence{A large language model (Claude, Anthropic) was used to help check
mathematical statements against the cited sources, to revise the text and to assist with \LaTeX. The
author reviewed all content and takes full responsibility for it.}

\conflictsofinterest{The author declares no conflicts of interest.}

\clearpage
\begin{adjustwidth}{-\extralength}{0cm}
\reftitle{References}
\bibliography{SipArxivMDPI}
\PublishersNote{}
\end{adjustwidth}
\end{document}